\documentclass[10pt]{article}
\usepackage[utf8]{inputenc}
\usepackage{enumitem}
\usepackage{fullpage}
\usepackage{amsmath}
\usepackage{amssymb}
\usepackage{authblk}
\usepackage{natbib}
\usepackage{xcolor}
\usepackage{verbatim}
\RequirePackage{graphicx}
\usepackage{bbm}
\usepackage{tikz, subcaption}
\usetikzlibrary{bayesnet, shapes, arrows, positioning}
\usepackage{subcaption}
\usepackage{url}
\usepackage{booktabs}

\usepackage{amsthm}
\theoremstyle{plain}
\newtheorem*{lemma*}{Lemma}
\newtheorem{theorem}{Theorem}
\newtheorem{lemma}{Lemma}

\newtheorem{proposition}{Proposition}

\newtheorem{assumption}{Assumption}
\title{\bf Multiply Robust Causal Mediation Analysis with\\ Continuous Treatments}

\author{Yizhen Xu$^{1}$, Numair Sani$^{2}$, Ilya Shpitser$^{2}$, AmirEmad Ghassami$^{*,3}$ \\
\small
1. \textit{Division of Biostatistics, Department of Population Health Sciences, University of Utah, USA}
\\
\small
2. \textit{Department of Computer Science, Johns Hopkins University, USA}
\\
\small
3. \textit{Department of Mathematics and Statistics, Boston University, USA}
}

\date{First Version: May 19, 2021; Current Version: November 16, 2026}

\begin{document}
\maketitle

\begin{abstract}
In many applications, researchers are interested in the direct and indirect causal effects of a treatment or exposure on an outcome of interest. Mediation analysis offers a rigorous framework for identifying and estimating these causal effects. For binary treatments, efficient estimators for the direct and indirect effects are presented by Tchetgen Tchetgen and Shpitser (2012) based on the influence function of the parameter of interest. These estimators possess desirable properties such as multiple-robustness and asymptotic normality while allowing for slower than root-n rates of convergence for the nuisance parameters. However, in settings involving continuous treatments, these influence function-based estimators are not readily applicable without making strong parametric assumptions. In this work, utilizing a kernel smoothing approach, we propose an estimator suitable for settings with continuous treatments inspired by the influence function-based estimation strategy. Our proposed approach employs cross-fitting, relaxing the smoothness requirements on the nuisance functions and allowing them to be estimated at slower rates than the target parameter. Additionally, similar to influence function-based estimators, our proposed estimator is multiply robust and asymptotically normal, allowing for inference in settings where parametric assumptions may not be justified.\\

\end{abstract}

\let\thefootnote\relax\footnotetext{~\\[-3mm]Y. Xu, A. Ghassami, and N. Sani contributed equally and are co–first authors.\\$*$ Corresponding author: AmirEmad Ghassami (Email: ghassami@bu.edu)}

\section{Introduction}
\label{sec:intro}

Estimating the causal effect of a treatment, policy, or intervention on an outcome of interest is a fundamental task in various fields such as epidemiology, economics, medicine, and sociology. A common parameter of interest is the average causal effect (ACE), which has been extensively studied \citep{hernan2020causal, imbens2015causal}. However, in addition to estimating the ACE, one may also be interested in the pathways and mechanisms through which the treatment affects the outcome of interest. Causal mediation analysis offers a precise and rigorous mathematical framework to answer such questions \citep{robins1992identifiability, tchetgen2012semiparametric,
pearl2001direct,
vanderweele2009marginal, goetgeluk2008estimation, imai2010identification, van2008direct, lange2011direct, lange2012simple}.

Much of the literature on mediation analysis assumes that the treatment of interest is binary. However, interventions involving the dosage of a drug, and the duration or frequency of an activity are better described as continuous variables. In such cases, mediation effects are naturally represented by a multi-dimensional surface rather than a scalar parameter. 
This learning task is challenging if a priori shape constraints are not imposed on the surface.
Additionally, the presence of continuous treatments complicates the estimation of nuisance parameters, making the estimation of causal parameters more challenging.

The challenges related to estimating ACE in the continuous treatment setting have been addressed in multiple works \citep{kennedy2017nonparametric, AI2021, imbens-continuous, kreif2015evaluation, 10.2307/2673642, su2019non,kallus2018policy, colangelo2020double, hill2011bayesian}. A common method
is based on outcome regression, which requires the correct specification of the relevant models, and hence machine learning methods such as Bayesian additive regression trees (BART) \citep{hill2011bayesian} are often used. However, this inherits the rate of the outcome regression estimation, and complex machine learning methods tend to have a slower convergence rate than simple parametric methods \citep{wasserman2006all, tsybakov2009nonparametric}. An alternative approach involves specifying a parametric form for the dose-response curve or projecting the true curve onto a parametric model, as presented in \cite{robins2000marginal, van1998locally, neugebauer2007nonparametric}. However, these methods may suffer from bias when the dose-response curve is misspecified. In contrast to approaches involving parametric assumptions on the dose-response curve, \cite{kennedy2017nonparametric}  leverage semiparametric theory by utilizing a two-stage estimator that first constructs a doubly robust pseudo-outcome in the first stage, and then regresses the pseudo-outcome on the treatment in the second stage using non-parametric regression methods. \citet{colangelo2020double} utilize double machine learning along with applying kernel smoothing to the augmented inverse propensity weighted (AIPW) score \citep{robins1995semiparametric}. This results in a slower convergence rate of nuisance parameters, while still guaranteeing fast rates for the target parameter. 
However, these approaches are not investigated for mediation analysis in the presence of continuous treatments.

In this paper, we propose a kernel smoothing approach inspired by influence function-based estimators \citep{tsiatis2007semiparametric,newey1994asymptotic, bickel1993efficient,ichimura2015influence, tchetgen2012semiparametric} to deal with continuous treatments for causal mediation analysis. We propose an estimator that, under mild regularity conditions, is consistent and asymptotically normal. Our work aims to extend the results for the continuous treatment ACE to the case of mediation analysis involving continuous treatments in the presence of complex nuisance functions. \citet{huber2020direct} tackle this problem by weighting the observations with a generalized propensity score that involves two nuisance functions, which are the conditional density of treatment given covariates and the conditional density of treatment given mediators and covariates. 
In their proposed approach, the nuisance functions can be estimated either parametrically or non-parametrically.
However, their estimator for the causal parameter is not robust with respect to the misspecification of the two nuisance functions and also inherits the rate of the nuisance function estimators, which could be slow. In contrast, we propose an approach motivated by influence function-based estimation strategy and hence obtain many of the desirable properties of influence functions, namely allowing for slower estimation of nuisance functions, as well as robustness properties. Our work draws from the existing causal mediation literature that discusses the identification and estimation of causal mediation effects \citep{pearl2001direct, imai2010identification, tchetgen2012semiparametric}. Additionally, we utilize the cross-fitting strategy to relax the smoothness assumptions on the nuisance functions \citep{chernozhukov2018double}. 
In a related line of work, \citet{munoz2012population} and \citet{diaz2020causal} study interventions on the treatment mechanism, treating the intervention and the resulting post-intervention exposure separately. Their approach requires pre-specifying a known function that maps a proposed intervention to the induced change in the continuous exposure. Building on this line of work, \citet{liu2024general} recently proposed general one-step machine learning estimators for a range of mediation parameters and extended them to continuous or multivariate exposures through modified treatment policies, which define effects as contrasts between user-specified shifts of the observed exposure. In contrast, in our setting, there is no separate intervention and exposure variables, the intervention variable itself is continuous (e.g., dose or intensity), and the target is the natural effect surface indexed by fixed treatment levels \((a,a')\), which is not pathwise differentiable. This distinction results in a different statistical parameter of interest and subsequently different needed analysis methodology. 

The remainder of this paper is organized as follows. Section \ref{sec:med} introduces the formal mediation analysis framework, describes its identifying assumptions, and discusses an influence function-based estimator of mediation effects for binary treatments. Section \ref{sec:c-med} extends the influence function-based approach to continuous treatment settings and describes the sample-splitting and smoothing procedures. In Section \ref{sec:asyp}, we provide our main results along with the required regularity conditions. Section \ref{simu} presents simulation
results, and in Section \ref{appli}, we apply our proposal to study the effect of the duration of Job Corp training on criminal status mediated by employment.

\section{Mediation Analysis}
\label{sec:med}

\begin{figure}[t!]
\centering
		\tikzstyle{block} = [draw, circle, inner sep=2.5pt, fill=lightgray]
		\tikzstyle{input} = [coordinate]
		\tikzstyle{output} = [coordinate]
        \begin{tikzpicture}[scale=0.75, transform shape]
            \tikzset{edge/.style = {->,> = latex'}}
            \node[] (a) at  (-2,0) {$A$};
            \node[] (x) at  (-3,2) {$X$};
            \node[] (m) at  (0,1) {$M$};
            \node[] (y) at  (2,0) {$Y$};    
            \node[] (la) at  (0,-1) {$(a)$};    
            \draw[-stealth] (a) to (m);                                                         
            \draw[-stealth,dotted] (a) to (y);
            \draw[-stealth,dotted] (x) to (a);                                                         
            \draw[-stealth,dotted] (x) to (m);
            \draw[-stealth,dotted,bend left=45] (x) to (y);
            \draw[-stealth] (m) to (y);
        \end{tikzpicture}
        \begin{tikzpicture}[scale=0.75, transform shape]
            \tikzset{edge/.style = {->,> = latex'}}
            \node[] (a) at  (-2,0) {$A$};
            \node[] (x) at  (-3,2) {$X$};
            \node[] (m) at  (0,1) {$M$};
            \node[] (y) at  (2,0) {$Y$}; 
            \node[] (lb) at  (0,-1) {$(b)$};
            \draw[-stealth,dotted] (a) to (m);                                                         
            \draw[-stealth] (a) to (y);
            \draw[-stealth,dotted] (x) to (a);                                                         
            \draw[-stealth,dotted] (x) to (m);
            \draw[-stealth,dotted,bend left=45] (x) to (y);
            \draw[-stealth,dotted] (m) to (y);
        \end{tikzpicture}
        \begin{tikzpicture}[scale=0.75, transform shape]
            \tikzset{edge/.style = {->,> = latex'}}
            \node[] (a) at  (-2,0) {$A$};
            \node[] (x) at  (-3,2) {$X$};
            \node[] (m) at  (0,1) {$M$};
            \node[] (y) at  (2,0) {$Y$};
            \node[] (lc) at  (0,-1) {$(c)$};
            \draw[-stealth] (a) to (m);                                                         
            \draw[-stealth] (a) to (y);
            \draw[-stealth,dotted] (x) to (a);                                                         
            \draw[-stealth,dotted] (x) to (m);
            \draw[-stealth,dotted,bend left=45] (x) to (y);
            \draw[-stealth] (m) to (y);
        \end{tikzpicture}
        \caption{A graphical representation of the decomposition of total effect into direct and indirect effects. Part $(a)$ represents the indirect effect, part $(b)$ represents the direct effect, and part $(c)$ represents the total effect.}
        \label{fig:MemMod}
\end{figure}

Let $A\in\mathcal{A}$ be the continuous treatment variable, $Y\in\mathcal{Y}$ be the outcome variable, and $M\in\mathcal{M}$ be a mediator variable that relays part of the causal effect of $A$ on $Y$. In addition, let $X\in\mathcal{X}$ denote the observed pre-treatment covariates in the setting. See Figure \ref{fig:MemMod}$(c)$ for a graphical representation of the causal relationships between the variables. To describe the causal effect of the treatment on the outcome, we use the potential outcome framework \citep{rubin1974estimating}. Let $Y^{(A=a)}$ be the random variable representing the potential outcome when the treatment is set to value $a$. We are interested in comparing the treatment values of $a$ and $a'$. A popular way to measure the causal effect of this change in treatment is to use the average causal effect (ACE), which captures the difference in the expected value of the potential outcome variables, that is
\[
ACE(a,a^\prime)=\mathbb{E}[Y^{(a)}-Y^{(a^\prime)}],
\]
where $\mathbb{E}[\cdot]$ denotes the population expectation operator. When no confusion arises, we occasionally suppress the dependence of indexed causal quantities on the treatment levels \(a\) and \(a'\) to simplify notation.

The total ACE of the treatment $A$ on the outcome $Y$ can be partitioned into the part mediated by variable $M$ and the part directly affecting outcome $Y$ (see Figure \ref{fig:MemMod}). 
To formally define this partitioning, let $Y^{(a,m)}$ denote the potential outcome variable corresponding to the outcome when the treatment is set to value $a$ and the mediator is set to value $m$, and $M^{(a)}$ denote the mediator variable when the treatment is set to value $a$. \citet{robins1992identifiability} and \citet{pearl2001direct} proposed the following partitioning of the ACE into the natural direct and indirect effects:
\begin{equation}
\label{eq:NDENIE}
\begin{aligned}
ACE(a,a^\prime)
&=\overbrace{\mathbb{E}[Y^{(a)}-Y^{(a^\prime)}]}^{\text{total effect}}\\
&=\mathbb{E}[Y^{(a,M^{(a)})}-Y^{(a',M^{(a^\prime)})}]\\
&=\underbrace{\mathbb{E}[Y^{(a,M^{(a)})}-Y^{(a,M^{(a^\prime)})}]}_{\text{natural indirect effect}}
+\underbrace{\mathbb{E}[Y^{(a,M^{(a^\prime)})}-Y^{(a',M^{(a^\prime)})}]}_{\text{natural direct effect}}.
\end{aligned}
\end{equation}
The two terms in Equation~\eqref{eq:NDENIE} define the natural indirect effect \(NIE(a,a')\) and the natural direct effect \(NDE(a,a')\), respectively. The natural direct effect (NDE) and natural indirect effect (NIE) can be described as follows. NDE captures the change in the expectation of the outcome if the value of the treatment variable is switched between the two arms of the experiment, while the mediator behaves as if the treatment has not changed. NIE captures the change in the expectation of the outcome if the value of the treatment variable is fixed, while the mediator behaves as if the treatment has been switched between the two arms of the experiment.
In the following subsection, we discuss the estimation of NDE and NIE from observational data.

\subsection{Estimating Natural Direct and Indirect Effects}

To estimate the natural direct and indirect effects, from the partitioning in Equation \eqref{eq:NDENIE}, it suffices to focus on estimating the parameter 
\[
\psi_0(a,a')=\mathbb{E}[Y^{(a,M^{(a^\prime)})}],
\]
for $a,a'\in\mathcal{A}$. Suppose i.i.d. data from a distribution $f$ on variables $O=\{A,X,M,Y\}$ are given. In general, the estimand \(\psi_0(a,a')\) is not identified from observational data, and identification assumptions are required to relate the distribution of the observational data to that of counterfactual variables. We require the following assumptions for the identification of \(\psi_0(a,a')\) from the observed distribution on variables, \(f(O)\).

\begin{assumption}[Identification Assumptions]
\label{assumption:id}
Let $X_1 \perp X_2 \mid X_3$ indicate that the random variables $X_1$ and $X_2$ are conditionally independent given the random variable $X_3$.
 \begin{enumerate}
 \item 
{\bf Consistency.} For all $a\in\mathcal{A}$ and $m\in\mathcal{M}$,
\begin{align*}
Y^{(a,m)}=Y&\quad\text{if }A=a\text{ and } M=m,\\
M^{(a)}=M&\quad\text{if }A=a.  
\end{align*}
\item 
{\bf Sequential Exchangeability.} For all $a,a'\in\mathcal{A}$, and $m\in\mathcal{M}$,
\begin{align*}
	& Y^{(a,m)}\perp \{A,M\}\mid X,\\
	& M^{(a)}\perp A\mid X,\\
	& Y^{(a,m)}\perp M^{(a^\prime)}\mid X.
\end{align*}
\item 
{\bf Positivity.} For all $a\in\mathcal{A}$, $m\in\mathcal{M}$ and $x\in\mathcal{X}$,
\begin{align*}
&f_{M|A,X}(m|A = a,X = x)>0,\\
&f_{A|X}(a|X = x)>0,
\end{align*}
where $f_{M|A,X}$ and $f_{A|X}$ are the conditional density of $M$ given $A$ and $X$, and the conditional density of $A$ given $X$, respectively.
\end{enumerate}
\end{assumption}

Assumption~\ref{assumption:id} is the continuous-treatment analogue of the consistency, sequential ignorability, and positivity conditions used to identify natural direct and indirect effects under a binary treatment \citep{pearl2001direct,imai2010identification,tchetgen2012semiparametric}. In particular, together with consistency, the sequential ignorability assumption of \citet{imai2010identification} adopted by \citet{tchetgen2012semiparametric} implies our sequential exchangeability conditions. The interpretation of these conditions, including the untestable cross-world exchangeability condition \(Y^{(a,m)}\perp M^{(a')}\mid X\), is the same as in the binary case and is discussed in Section~1 of the Supplementary Material. There we also note that when the reference level $a'$ is fixed in advance, as in our application, the cross-world condition need only hold for that single reference level rather than for all $a'\in\mathcal A$. The main difference from the binary setting is the target: it is no longer only 4 target values corresponding to $a,a'\in\{0,1\}$. Instead the target is the function $(a,a')\mapsto\psi_0(a,a')$ defined with  $a$ and $a'$ that can take values in a continuum. Hence, the NDE and NIE are effect curves or surfaces indexed by the pair of treatment levels, and Assumption~\ref{assumption:id} must hold over a continuum of treatment values. The consequences for estimation and inference are summarized in the discussion section.

Under Assumption \ref{assumption:id}, the estimand $\psi_0(a,a')$ can be identified from the observed distribution $f(O)$ using the following expression called the mediation formula \citep{robins1992identifiability, pearl2001direct, imai2010identification}:
\begin{equation}
\label{eq:MCF}
\psi_0(a,a')=\int_{\mathcal{X}}\int_{\mathcal{M}}\mathbb{E}[Y|A=a,M=m,X=x]f_{M|A,X}(m|A=a',X=x)f_X(x)dm dx,
\end{equation}
where $f_X$ is the marginal distribution of $X$.

Using Equation \eqref{eq:MCF}, one can estimate the parameter of interest $\psi_0(a,a')$ by first estimating the nuisance functions $\mathbb{E}[Y|A,M,X]$ and $f_{M|A,X}$, and then using a plug-in estimator to estimate \(\psi_0(a,a')\) as follows
\[
\frac{1}{n}\sum_{i=1}^n
\int_{\mathcal{M}}\hat{\mathbb{E}}[Y_i|A=a,M=m,X_i]\hat{f}_{M|A,X}(m|A=a',X_i)dm.
\]
Unfortunately, this estimator is sensitive to bias in the estimation of nuisance functions. That is, misspecifying either of the nuisance functions induces bias in the estimation of the parameter of interest.

As an alternative approach, in the case of binary treatment, that is, $\mathcal{A}=\{0,1\}$, \citet{tchetgen2012semiparametric} developed a semiparametric approach to inference for mediation analysis. 
They derived the efficient influence function for $\psi_0(a,a')$ as
\begin{equation}
\label{eq:IFb}
\begin{aligned}
IF_{\psi_0}(O) 
&= I(A = a)\lambda(a,X)\frac{\alpha(a',M,X)}{\alpha(a,M,X)}\{Y - \gamma(X,M,a)\}\\
&\quad+ I(A = a')\lambda(a',X)\{\gamma(X,M,a) - \eta(a, a', X)\} + \eta(a, a', X) - \psi_0(a,a'),
\end{aligned}
\end{equation}
where $\lambda(a,X) := 1/{f(a|X)}$, $\alpha(a,M,X) := f(M| a, X)$, and $\gamma(X,M,a) := \mathbb{E}[Y|A=a,M,X]$ are the nuisance functions, $a,a'\in\{0,1\}$, $I(\cdot)$ denotes the indicator function, and
\[
\eta(a, a', X) = \int_{\mathcal{M}} \gamma(X, m, a) \alpha(a', m,X) dm.
\]
Note that $IF_{\psi_0}$ is comprised of three nuisance functions: $\lambda$, $
\alpha$, and $\gamma$. \cite{tchetgen2012semiparametric} showed that the estimator based on this influence function has the multiple robustness property, that is, it is consistent even if the model for one (but not more than one) nuisance function is misspecified. Formally, let 
\begin{itemize}
\item $\mathfrak{M}_{ym}$ be the sub-model in which the model for $\gamma$ and $\alpha$ are correctly specified.
\item $\mathfrak{M}_{ya}$ be the sub-model in which the model for $\gamma$ and $\lambda$ are correctly specified.
\item $\mathfrak{M}_{ma}$ be the sub-model in which the model for $\alpha$ and $\lambda$ are correctly specified.
\end{itemize}
 The estimator for \(\psi_0(a,a')\) based on the influence function $IF_{\psi_0}$ defined as
\begin{align*}
\hat{\psi}^{TTS}(a,a')
= \frac{1}{n}\sum^n_{i=1}&\bigg\{ I(A_i = a)\hat{\lambda}(a,X_i)\frac{\hat{\alpha}(a',M_i,X_i)}{\hat{\alpha}(a,M_i,X_i)}\{Y_i - \hat{\gamma}(X_i,M_i,a)\}\\
&\quad+ I(A_i = a')\hat{\lambda}(a',X_i)\{\hat{\gamma}(X_i,M_i,a) - \hat{\eta}(a, a', X_i)\} + \hat{\eta}(a, a', X_i)\bigg\},
\end{align*}
is consistent when the truth lies in the submodel union $\mathfrak{M}_{ym}\cup\mathfrak{M}_{ya}\cup\mathfrak{M}_{ma}$, all estimators of nuisance functions converge to some functions in probability,  and the estimators of nuisance functions in the correctly specified submodels are consistent, where $$\hat{\eta}(a, a', X) = \int_{\mathcal{M}} \hat{\gamma}(X,m,a) \hat{\alpha}(a',m, X)dm.$$

Inspired by this result, in the following section, we propose a kernel-based estimator for mediation effects in settings with continuous treatment variables, while preserving multiple robustness and allowing for the nuisance parameters to be estimated at a slower rate than the parameter of interest.

\section{Continuous-Treatment Mediation Analysis}
\label{sec:c-med}
 
In the case of continuous treatments, the parameter of interest, $\psi_0(a,a')$, is no longer regular \citep{bickel1993efficient,ichimura2015influence}. Therefore, the method of \cite{tchetgen2012semiparametric} cannot be applied directly. 
However, their estimator can be modified to be suitable for inference in the case of continuous treatments, while still obtaining desirable properties such as asymptotic normality, robustness to misspecification of nuisance functions, and valid inference while allowing for 
the nuisance parameters to be estimated at a slower rate than the parameter of interest.
Specifically, we modify $\hat{\psi}^{TTS}(a,a')$ by employing a kernel smoothing technique, wherein the indicators in the calculation of $\hat{\psi}^{TTS}(a, a')$ 
 are replaced by kernel-based weights. The weights are functions of treatment values falling within a neighborhood (defined by the bandwidth parameter $h$) of $a$ and $a'$. This modification introduces several challenges in the inference, which we will present and address in Section \ref{sec:asyp}.

Let $d_A$ denote the dimension of the treatment variable, and let
\[
K_h(a):=\frac{1}{h^{d_A}}\prod_{j=1}^{d_A}k\Big(\frac{a_j}{h}\Big),
\]
where $k(\cdot)$ is a kernel function, and $h$ denotes the bandwidth parameter.
We propose to use the following modification of the efficient influence function in Equation \eqref{eq:IFb} for any $a$ and $a' \in \mathcal{A}$:

\begin{align}\label{eq:IF}
m(O;\alpha,\lambda,\gamma,\psi(a,a'))
&= K_h(A-a)\,\lambda(a,X)\,\frac{\alpha(a',M,X)}{\alpha(a,M,X)}\{Y-\gamma(X,M,a)\}\nonumber\\
&\quad + K_h(A-a')\,\lambda(a',X)\{\gamma(X,M,a)-\eta(a,a',X)\}\nonumber\\
&\quad + \eta(a,a',X) - \psi(a,a') .
\end{align}

The moment function in Equation~\eqref{eq:IF} is parametrized by the nuisance functions \((\alpha,\lambda,\gamma)\), and in our implementation \(\eta(a,a',X)\) is computed by Monte Carlo integration over the fitted conditional mediator distribution. Alternative parametrizations and computations are also possible, such as rewriting the mediator density ratio in terms of the conditional treatment density \(f_{A\mid M,X}\) following Bayes-rule \citep{diaz2020causal,farbmacher2022causal} or estimating \(\eta(a,a',X)\) by nested regression; these involve different nuisance functions with their own rate conditions. Details are discussed in the Supplementary Material.

Note that in Equation \eqref{eq:IF}, $m(O;\alpha,\lambda,\gamma,\psi(a,a'))$ also depends on the choice of kernel function and its bandwidth $h$. For simplicity, and with a slight abuse of notation, we omit $K_h$ from the notation for $m(\cdot)$. To derive desired results on consistency, asymptotic normality, and multiple robustness, we require the kernel $k(\cdot)$ to satisfy the following conditions.
\begin{assumption}[Kernel \& Bandwidth Assumptions]
\label{assumption:kernel}
The kernel function $k(\cdot)$ satisfies
\begin{enumerate}
    \item $\int k(u) du = 1$
    \item $\int uk(u) du = 0$
    \item $0 < \int u^6 k(u) du < \infty$ 
    \item $\int k^2(u) du < \infty$
    \item $0 < \int u^2 k^2(u) du < \infty$
\end{enumerate}

Additionally, the kernel bandwidth $h$ is assumed to be a function of the sample size $n$ and satisfies $h \rightarrow 0$, $nh^{d_A} \rightarrow \infty$ and $nh^{d_A + 4} \rightarrow C_h$, for a constant $C_h$, as $n\rightarrow \infty$.
\end{assumption}

These assumptions are satisfied by common kernels such as the Gaussian kernel and Epanechnikov kernel. Note that in the moment function in Equation \eqref{eq:IF}, the nuisances are not functions of the parameter of interest $\psi_0(a, a')$. Therefore, having estimators for nuisance functions suffices for obtaining an estimator for $\psi_0(a, a')$. Next, we describe the estimation procedure for utilizing Equation \eqref{eq:IF} to estimate $\psi_0(a, a')$.

In applications, the treatment values \(a\) and \(a^\prime\) define the causal contrast of interest. We recommend choosing \(a^\prime\) as a substantively meaningful reference exposure, such as a low, standard, baseline, or policy-relevant treatment level, and choosing \(a\) values that correspond to feasible and interpretable alternative exposure levels. Both \(a\) and \(a^\prime\) should lie in regions with adequate empirical support, so that the positivity assumption is plausible and the nuisance estimations do not lead to unstable weights dominated by a small number of observations. In some applications, no single value for \(a\) and/or \(a^\prime\) may be naturally preferred. In such settings, one may prespecify a grid of scientifically meaningful values for both \(a\) and \(a^\prime\), and report the NDE and NIE over the resulting set of treatment contrasts. This provides a more complete summary of the mediation surface and avoids selecting a single contrast post hoc. The source code in our GitHub repository (\url{https://github.com/yizhenxu/Continuous-Treatment-Mediation.git}) implements parallelization over treatment contrasts, making repeated evaluation over a prespecified grid computationally feasible in practice.

\vspace{1em}
\noindent
\textbf{Estimation Procedure.} We use the cross-fitting estimation approach of \cite{chernozhukov2018double} for separating the estimation of the nuisance functions from the parameter of interest.
This approach is beneficial since weaker smoothness requirements are needed for the estimation of nuisance functions. In the cross-fitting approach, we partition the sample indices into $L$ folds $\{I_1,...,I_L\}$ of roughly equal size. Data from the $\ell$-th fold is denoted by $O_{I_\ell}$, and the data in the rest of the folds is denoted by $O^c_{I_\ell}$.
For $\ell\in\{1,...,L\}$, we estimate the nuisance functions $\alpha$, $
\lambda$, $\gamma$ by $\hat{\alpha}_\ell,\hat{\lambda}_\ell,\hat{\gamma}_\ell$ based on data $O^c_{I_{\ell}}$.
For all $\ell$, let $\hat{\psi}_\ell$ be the estimator for $\psi_0(a,a')$ obtained by solving
\[
\frac{1}{|I_\ell|}\sum_{i\in I_\ell}m(O_i;\hat{\alpha}_\ell,\hat{\lambda}_\ell,\hat{\gamma}_\ell,\hat{\psi}_\ell(a,a'))=0.
\]
Our proposed estimator for $\psi_0(a,a')$ is
\begin{equation}
\label{eq:cf-estimator}
\hat{\psi}^{MR}(a,a')=\frac{1}{L}\sum_{\ell=1}^L \hat{\psi}_\ell(a,a'),
\end{equation}
where MR stands for multiply robust. 
Since \(\psi_0(a,a)=\mathbb{E}[Y^{(a)}]\) and \(\psi_0(a',a')=\mathbb{E}[Y^{(a')}]\), the natural direct and indirect effects in Equation~\eqref{eq:NDENIE} are estimated by \(\widehat{NDE}(a,a')=\hat{\psi}^{MR}(a,a')-\hat{\psi}^{MR}(a',a')\) and \(\widehat{NIE}(a,a')=\hat{\psi}^{MR}(a,a)-\hat{\psi}^{MR}(a,a')\).  In the next section, we present the asymptotic properties of our proposed estimator, along with the required regularity conditions.

\medskip

\noindent
\textbf{Remark on practical nuisance estimation.}
In implementation, $\lambda(a,x)$ is obtained by estimating the conditional treatment density $f_{A|X}(a|x)$ and setting $\hat{\lambda}(a,x)=1/\hat f_{A|X}(a|x)$. For a univariate continuous treatment, common parametric choices include normal, log-normal, gamma, or other generalized propensity score models, possibly with flexible mean and variance functions \citep{10.2307/2673642,imbens-continuous,huber2020direct}. For multivariate treatments, one may instead use a joint conditional density model for $A\mid X$, such as a multivariate normal location-scale model, a mixture model, a copula-based model, or a flexible conditional density estimator. Nonparametric and machine learning approaches to conditional density estimation, including kernel and mixed-data estimators, orthogonal-series conditional density estimators, and least-squares conditional density or density-ratio estimators, may also be used when their convergence rates are compatible with Assumption~\ref{assumption:DR} \citep{li2008nonparametric,izbicki2017converting,sugiyama2010conditional}. 

The nuisance function $\alpha(a,m,x)=f_{M|A,X}(m|a,x)$ can be estimated similarly as a conditional density or probability mass function for the mediator. If $M$ is binary or categorical, $\alpha$ can be estimated using logistic, multinomial, or other probabilistic classification methods. If $M$ is continuous, $\alpha$ can be estimated using parametric conditional density models or nonparametric conditional density estimators. For a multivariate mediator $M=(M_1,\ldots,M_{d_M})$, directly estimating the joint conditional density may be challenging; a useful alternative is the factorization
\[
f_{M|A,X}(m|a,x)
=
\prod_{j=1}^{d_M}
f_{M_j|A,X,M_{1:(j-1)}}(m_j|a,x,m_{1:(j-1)}),
\]
which reduces the problem to a sequence of lower-dimensional conditional density or classification problems. In applications where only the ratio $\alpha(a',M,X)/\alpha(a,M,X)$ is required, direct conditional density-ratio estimation can also be used, provided that the resulting estimator satisfies the corresponding product-rate requirement. Regardless of the particular nuisance-estimation strategy, the theoretical results require the pairwise product rates in Assumption~\ref{assumption:DR}. Specifically, in Section 6, we use reproducing kernel Hilbert space (RKHS) conditional mean embedding estimator for density estimation for both nuisances.
\hfill$\diamond$

\section{Asymptotic Analysis}
\label{sec:asyp}

In this section, we provide asymptotic properties of our proposed estimator $\hat{\psi}^{MR}(a,a')$ in Equation \eqref{eq:cf-estimator}. We start by stating the required regularity conditions.

\begin{assumption}[Regularity Conditions]
\label{assumption:reg}
~\\
\vspace{-5mm}
\begin{enumerate}
\begin{sloppypar}
\item For all $Y$, $M$, and $X$, the functions $f(a \mid Y, M, X)$, $f(a \mid M, X)$, $f(a \mid X)$, $\gamma(X,M,a)$ as a function of $a$ are three times continuously differentiable with respect to each dimension of $a$, and the functions and their first, second, and third derivatives with respect to $a$ are bounded in each dimension.
\end{sloppypar}
\item The nuisance functions $\alpha,\lambda, \gamma$ and the estimators $\hat{\alpha},\hat{\lambda},\hat{\gamma}$ are bounded. Additionally, $\alpha$, $\lambda$ and their estimators $\hat{\alpha}$, $\hat{\lambda}$ are bounded away from zero.
\item $Y$'s conditional variance $var(Y|a,m,x)$ and its first and second derivative with respect to each dimension of $a$ are bounded for any $a\in \mathcal{A}$, $m\in \mathcal{M}$, and $x\in \mathcal{X}$.

\end{enumerate}
\end{assumption}

In addition to the regularity conditions, we require the following conditions regarding the convergence of the estimators of the nuisance functions.

\begin{assumption}[Convergence of Nuisance Estimators]
\label{assumption:converegence}
~\\
For any value $a \in \mathcal{A}$, the estimators $\hat{\alpha}(a,M, X)$, $\hat{\lambda}(a,X)$, and $\hat{\gamma}(X,M,a)$ satisfy the following conditions:
\begin{enumerate}
    \item $\int\left(\hat{\lambda}(a, x) - \lambda(a, x) \right)^2 f_{X}(x) dx \xrightarrow[]{P} 0$,
    \item $\int \left(\hat{\alpha}(a, m, x) - \alpha(a, m, x) \right)^2 f_{M, X}(m, x) dm dx \xrightarrow[]{P} 0$, 
    \item $\int\left(\hat{\gamma}(x,m,a) - \gamma(x,m,a) \right)^2 f_{M, X}(m, x) dm dx \xrightarrow[]{P} 0$,
\end{enumerate}
where $\xrightarrow[]{P}$ indicates convergence in probability.

\end{assumption}

 Similar to influence function-based estimators, in Assumption \ref{assumption:converegence}, we do not require individual nuisance estimators to satisfy convergence rate conditions. However, in our proposed method, we have requirements on the convergence rate of the product for the nuisance estimators as follows. 

\begin{assumption}[Nuisance Convergence Rates]
\label{assumption:DR}
~\\
For any value $a, a' \in \mathcal{A}$, the estimators $\hat{\alpha}(a,M, X)$, $\hat{\lambda}(a,X)$, and $\hat{\gamma}(X,M,a)$ satisfy the following conditions:
{\footnotesize
\begin{enumerate}
    \item
        \begin{align*}
        &\sqrt{nh^{d_A}}\left(\int \left(\hat{\alpha}(a', m, x) - \alpha(a', m, x) \right)^2f_{M, X}(m,x)dm dx\right)^{\frac{1}{2}}\\
        &\hspace{1.2in}\times \left(\int \left(\hat{\gamma}(x,m,a) - \gamma(x,m,a) \right)^2f_{M, X}(m,x) dm dx\right)^{\frac{1}{2}} \xrightarrow[]{P}0,
        \end{align*}
    \item
       \[
        \sqrt{nh^{d_A}}\left(\int \left(\hat{\lambda}(a', x) - \lambda(a', x) \right)^2f_{X}( x)dx\right)^{\frac{1}{2}}
        \left(\int \left(\hat{\gamma}(x,m,a) - \gamma(x,m,a) \right)^2f_{M, X}(m,x) dm dx\right)^{\frac{1}{2}} \xrightarrow[]{P} 0,
        \]
    \item 
        \[
        \sqrt{nh^{d_A}}\left(\int \left(\hat{\lambda}(a',x) - \lambda(a', x) \right)^2f_{X}(x) dx\right)^{\frac{1}{2}}\left(\int \left(\hat{\alpha}(a, m, x) - \alpha(a, m, x) \right)^2f_{M, X}(m,x)dm dx\right)^{\frac{1}{2}}
         \xrightarrow[]{P} 0.
         \]
\end{enumerate}
}
\end{assumption}

As seen in Assumption \ref{assumption:DR},
our requirements on the convergence rate of nuisance function estimators are on the product of the error rates, rather than on the individual nuisance function estimators. Therefore, if one of the estimators converges at a slow rate, the other estimator can compensate. This is a desirable property when working with non-parametric estimators since they typically have slow rates of convergence. Note that Assumption \ref{assumption:DR} includes $\sqrt{h^{d_A}}$ that is not involved in the standard influence function-based approach for binary treatment \citep{tchetgen2012semiparametric}.

\medskip

\noindent
\textbf{Remark on the interpretation of Assumption~\ref{assumption:DR}.}
Let\\
\(
r_\lambda(a)=\{\int(\hat{\lambda}(a,x)-\lambda(a,x))^2 f_X(x)\,dx\}^{1/2},
\)\\
\(
r_\alpha(a)=\{\int(\hat{\alpha}(a,m,x)-\alpha(a,m,x))^2 f_{M,X}(m,x)\,dm\,dx\}^{1/2},
\)
and\\
\(
r_\gamma(a)=\left\{\int(\hat{\gamma}(x,m,a)-\gamma(x,m,a)\right)^2 f_{M,X}(m,x)\,dm\,dx\}^{1/2}.
\)\\
Then Assumption~\ref{assumption:DR} can be summarized as requiring the relevant pairwise products
$r_\alpha(a^\prime)r_\gamma(a)$,
$r_\lambda(a^\prime)r_\gamma(a)$,
$r_\lambda(a^\prime)r_\alpha(a)$
to be $o_P\{(nh^{d_A})^{-1/2}\}$. Hence, the effective sample size for the smoothed estimating equation is $nh^{d_A}$ rather than $n$. This differs from the standard EIF-based requirement for binary or discrete treatments, where the corresponding pairwise product condition is typically $o_P(n^{-1/2})$ \citep{tchetgen2012semiparametric,chernozhukov2018double}. Since $h\to 0$, the target parameter in the continuous-treatment case is estimated at the slower rate $(nh^{d_A})^{-1/2}$, and the nuisance product-rate requirement is correspondingly less stringent than the binary-treatment requirement, although the overall target parameter is also estimated more slowly.

For example, the bandwidth that balances the squared smoothing bias and the variance satisfies
\(
h\asymp n^{-1/(d_A+4)}.
\)
Under this rate-optimal choice, Assumption~\ref{assumption:DR} requires
\[
r_\alpha(a^\prime)r_\gamma(a),\quad
r_\lambda(a^\prime)r_\gamma(a),\quad
r_\lambda(a^\prime)r_\alpha(a)
=
o_P\left(n^{-2/(d_A+4)}\right).
\]
Equivalently, if the two nuisance estimators in a given product have comparable rates, each needs to converge faster than $n^{-1/(d_A+4)}$. When $d_A=1$, this gives the nonparametric rate threshold $n^{-1/5}$ for each nuisance estimator in the equal-rate case, compared with the faster-than-$n^{-1/4}$ requirement in the classical binary-treatment EIF setting. This requirement is satisfied with, e.g., Gaussian RKHS function class we used in data application \citep{wainwright2019high}.

More generally, suppose a nuisance estimator has an $L_2$ convergence rate $n^{-\beta}$, up to logarithmic factors. Then the product condition for a pair of nuisance estimators with exponents $\beta_j$ and $\beta_k$ is
\(
\beta_j+\beta_k>\frac{2}{d_A+4}
\)
under the rate-optimal bandwidth. For standard nonparametric estimators over an $s$-smooth function class with effective dimension $p$, rates of the form $n^{-s/(2s+p)}$ are typical \citep{wasserman2006all,tsybakov2009nonparametric}. Thus, in the equal-smoothness case, the above condition is satisfied when
$\frac{s}{2s+p}>\frac{1}{d_A+4}$,
or equivalently
$s>\frac{p}{d_A+2}$.
Parametric nuisance estimators are more than sufficient for Assumption~\ref{assumption:DR}, and flexible machine learning estimators may also be used when their rates satisfy the pairwise product conditions above. 
\hfill$\diamond$

\medskip

In the case of binary treatment variables, the combination of assumptions similar to Assumptions \ref{assumption:converegence} and \ref{assumption:DR} can lead to asymptotic normality, which is used to construct Wald-style confidence intervals. However, when the treatments are continuous, the Central Limit Theorem (CLT) cannot be directly applied to our proposed method because the bandwidth $h$ varies as a function of the sample size $n$, implying that the distribution of Equation \eqref{eq:IF} changes with $n$. 
Instead, we impose additional assumptions stated below to satisfy the Lyapunov's condition for CLT and achieve asymptotic normality. 

\begin{assumption}[Assumptions for Lyapunov CLT]
\label{assumption:lyapunov-conditions}
~
\begin{enumerate}
    \item $\mathbb{E}\left[|Y-\gamma(X,M,a)|^3 \big| A = a', M = m, X = x\right]$ is bounded for any $(a,a',m,x)\in \mathcal{A}\times \mathcal{A}\times\mathcal{M} \times \mathcal{X}$.
    \item $\int^{\infty}_{-\infty} k(u)^{c_1}k(u+\tilde{c})^{c_2} du < \infty$ and $\int^{\infty}_{-\infty} u^2k(u)^{c_1}k(u+\tilde{c})^{c_2} du < \infty$ for $\tilde{c} \in \mathbb{R}$ and  $c_1,c_2\in\{0,1,2,3\}$ such that $c_1+c_2 \in \{2,3\}$.
\end{enumerate}        
\end{assumption}

In practice, almost all commonly used kernel functions (Gaussian, Epanechnikov, triangular, biweight, etc.) satisfy condition 2 in Assumption \ref{assumption:lyapunov-conditions}. The main exclusion is kernels with too heavy tails (like Cauchy, with $\alpha$=1), which would make the second moment condition diverge.
Having stated the assumptions in our setting, we now provide the following result regarding the asymptotic behavior of the proposed estimator in Equation \eqref{eq:cf-estimator}.

\begin{theorem}
\label{thm:convergence}
Under Assumptions \ref{assumption:id}-\ref{assumption:DR}, for any values of $a,a'\in\mathcal{A}$,
\[
\sqrt{nh^{d_A}}(\hat{\psi}^{MR}(a,a')-\psi_0(a,a'))=
\sqrt{\frac{h^{d_A}}{n}}
\sum_{i=1}^nm(O_i;\alpha,\lambda,\gamma,\psi_0(a,a'))+o_p(1).
\]
Additionally, if Assumption \ref{assumption:lyapunov-conditions} holds, then $\sqrt{nh^{d_A}}(\hat{\psi}^{MR}(a,a')-\psi_0(a,a')-h^2 B(a,a'))$ converges to the Gaussian distribution $\mathcal{N}(0,V(a,a'))$, where $B(a,a')$ and $V(a, a^\prime)$  are defined as
\begin{equation*}
\begin{aligned}
B(a,a')
&= \Bigl[\int u^{2}\,k(u)\,du\Bigr]\,
   \mathbb{E}\Biggl[
      \frac{\alpha(a',M,X)}{\alpha(a,M,X)\,\lambda(a,X)}
      \Biggl\{
        \sum_{j=1}^{d_A} \partial_{a_j}\gamma(X,M,a)\,\partial_{a_j} f(a\mid X,M) \\
&\qquad\qquad\qquad\qquad
        + \frac{1}{2}\Bigl(\sum_{j=1}^{d_A} \partial^{2}_{a_j}\gamma(X,M,a)\Bigr)\,f(a\mid X,M)
      \Biggr\} \\
&\qquad\qquad
      + \bigl\{\gamma(X,M,a) - \eta(a,a',X)\bigr\}\,
        \frac{1}{2}\,
        \frac{\sum_{j=1}^{h_{d_A}} \partial^{2}_{a_j} f(a'\mid X,M)}{\lambda(a',X)}
   \Biggr]
   + O(h)
\end{aligned}
\end{equation*}
and
\begin{equation*}
    \begin{aligned}
    V(a, a^\prime) = &\left[\int k(u)^2  du\right]^{d_A}\mathbb{E}\Bigg\{\frac{\alpha^2(a', M, X) f(a|X,M)}{\alpha^2(a, M, X)} \lambda^2(a|X)var(Y|X,M,a)\\
    &\qquad\qquad+ \lambda(a'|X)var[E(Y|X,M,a)|X,a']\Bigg\}.
\end{aligned}
\end{equation*}

\end{theorem}
All proofs are provided in the Supplementary Material. 
Theorem \ref{thm:convergence} provides results on the point-wise convergence of $\hat{\psi}^{MR}(a,a')$ and establishes the asymptotic normality of our estimator. Additionally, $\hat{\psi}^{MR}(a,a')$ has a multiple robustness property analogous to the $\hat{\psi}^{TTS}(a,a')$, formally stated in Proposition \ref{lem:multiple-robustness}.

\begin{proposition}\label{lem:multiple-robustness}
    Under Assumptions \ref{assumption:id}, \ref{assumption:kernel}, \ref{assumption:reg}, and \ref{assumption:lyapunov-conditions}, the proposed estimator $\hat{\psi}^{MR}(a, a')$ will be a consistent estimator for $\psi_0(a, a')$ as long as all three nuisance function estimators converge in probability to some functions and any two out of the three conditions in Assumption \ref{assumption:converegence} hold.    
\end{proposition}

While Theorem \ref{thm:convergence} and Proposition \ref{lem:multiple-robustness} establish properties of $\hat{\psi}^{MR}(a,a')$ that are desirable for point estimation, uncertainty quantification through the calculation of valid confidence intervals requires the estimation of $V(a, a')$ and $B(a, a')$. However, these are hard to estimate due to their complicated analytical forms. Nevertheless, by choosing an undersmoothing bandwidth $h$ that satisfies $\sqrt{nh^{d_A + 4} } \rightarrow 0$, valid confidence intervals can still be constructed without estimating $B(a, a')$. This is due to the fact that the bias of local smoother is of order $h^2$, while the standard deviation is of order $(nh^{d_A})^{-1/2}$.
Choosing a bandwidth $h$ such that $\sqrt{n h^{d_A+4}} \to 0$ ensures that the bias vanishes faster than the standard deviation and thus becomes asymptotically negligible for inference. The dependence on \(d_A\) also provides practical guidance on the dimensionality of the treatment variable for which the proposed approach is feasible. The product kernel assigns non-negligible weight mainly to observations whose treatment values fall in an \(h\)-neighborhood of the target value \(a\). If the treatment density is regular near $a$, the probability that an observation lies in such a neighborhood is of order \(h^{d_A}\), so the effective local sample size around $a$ is of order $n h^{d_A}$. As \(d_A\) increases, this effective sample size decreases rapidly for fixed \(n\) and \(h\), reflecting the usual curse of dimensionality for kernel smoothing. Combining the second-order smoothing bias $O(h^2)$ with the stochastic error $O((n h^{d_A})^{-1/2})$, the leading bias--variance tradeoff is heuristically summarized as $\operatorname{MSE}(h)
\approx
h^4+\frac{1}{n h^{d_A}}$. Balancing the two terms gives the usual kernel smoothing bandwidth order \(h\asymp n^{-1/(d_A+4)}\), and the corresponding pointwise estimation error rate is \(n^{-2/(d_A+4)}\). This rate deteriorates as \(d_A\) increases, showing the impact of the curse of dimensionality on convergence. Therefore, the proposed product-kernel implementation is most practical for scalar or low-dimensional continuous treatments. For moderate- or high-dimensional treatments, substantially larger sample sizes and stronger overlap would be required; otherwise, one may need to prespecify lower-dimensional treatment summaries or impose additional structure on the treatment-response surface.

Given that the bias is asymptotically negligible, we focus on estimating $V(a, a')$. Naturally, an estimator for $V(a, a')$ can be constructed as follows.
\[\widehat{V}(a, a') = h^{d_A}\frac{1}{L}\sum_{\ell = 1}^L \frac{1}{|I_\ell|}\sum_{i \in I_\ell} m^2(O_i;\hat{\alpha},\hat{\lambda},\hat{\gamma},\hat{\psi}_\ell(a,a')).
\]

We present the additional assumptions necessary for the consistency of $\widehat{V}(a, a')$ below. 

\begin{assumption}[Consistency of $\widehat{V}(a, a')$]
\label{assumption:v-hat}
~
\begin{enumerate}
    \item $\mathbb{E}\{[Y-\gamma(X,M,a)]^4 | A = a', M = m, X = x\}$ is bounded for any $(a,a',m,x)\in \mathcal{A}\times \mathcal{A}\times\mathcal{M} \times \mathcal{X}$.

    \item $\int^{\infty}_{-\infty} k(u)^{c_1}k(u+\tilde{c})^{c_2} du < \infty$ and $\int^{\infty}_{-\infty} u^2k(u)^{c_1}k(u+\tilde{c})^{c_2} du < \infty$ for $\tilde{c} \in \mathbb{R}$ and $c_1,c_2\in\{0,1,2,3,4\}$ such that $c_1+c_2 \in \{2,3,4\}$.
\end{enumerate}    
\end{assumption}

Assumption \ref{assumption:v-hat} is similar to Assumption \ref{assumption:lyapunov-conditions} but of higher order. Assumption 7.2 is satisfied by commonly used non-compact kernels with sufficiently fast-decaying tails in standard kernel smoothing. For example, consider the Gaussian kernel $k(u)=(2\pi)^{-1/2}\exp(-u^2/2).$
For any fixed $\tilde c\in\mathbb R$ and $c_1,c_2\in\{0,1,2,3,4\}$ with $c_1+c_2\in\{2,3,4\}$, there is $k(u)^{c_1}k(u+\tilde c)^{c_2}
=(2\pi)^{-(c_1+c_2)/2}
\exp\{-\frac{1}{2}[c_1 u^2+c_2 (u+\tilde{c})^2]\}\propto C \exp\{-\frac{c_1+c_2}{2}(u+\frac{c_2\tilde c}{c_1+c_2})^2
\}$. The product $k(u)^{c_1}k(u+\tilde c)^{c_2}$ is proportional to a Gaussian density in $u$. Therefore both
$\int_{-\infty}^{\infty} k(u)^{c_1}k(u+\tilde c)^{c_2} du$ and $\int_{-\infty}^{\infty} u^2 k(u)^{c_1}k(u+\tilde c)^{c_2}du
$
are finite. Hence the Gaussian kernel satisfies Assumption 7.2. Bounded compactly supported kernels, such as the uniform, triangular, and Epanechnikov kernels, also satisfy the assumption immediately because the relevant integrands are bounded and have compact support. We have the following result regarding the consistency of $\widehat{V}(a, a')$.

\begin{proposition}\label{lem:v-hat-estimator}
    Under Assumptions \ref{assumption:id}-\ref{assumption:v-hat}, for any values of $a,a'\in\mathcal{A}$, $\widehat{V}(a, a')$ is a consistent estimator for $V(a, a')$.
\end{proposition}

Using the result of Proposition \ref{lem:v-hat-estimator}, $\widehat{V}(a, a')$ can be used to construct asymptotically valid confidence intervals as follows. Choose an undersmoothing bandwidth $h$ that satisfies $\sqrt{nh^{d_A + 4}} \rightarrow 0$, so that $\sqrt{nh^{d_A}}h^2 B(a, a')$ is asymptotically negligible. Then, the $(1 - \alpha)$ confidence interval is given as 
\begin{equation}
\label{eq:cie}
\Bigg[\hat{\psi}^{MR}(a, a') \pm \Phi^{-1}(1 - \alpha/2)\sqrt{\frac{\widehat{V}(a, a')}{nh^{d_A}}}\Bigg],
\end{equation}
where $\Phi$ is the CDF of $\mathcal{N}(0,1)$.  However, there is little practical guidance on how to implement undersmoothing, and in most applications it functions primarily as a technical device to simplify asymptotic derivations \citep{kennedy2017nonparametric}.  In particular, undersmoothing sequences are not unique and hence choosing some bandwidth satisfying the undersmoothing condition may induce the impression of arbitrary tuning. In our simulations and applications, we first follow a more principled and widely used data-driven approach by selecting the bandwidth according to Silverman’s rule of thumb \citep{van2000asymptotic, silverman2018density}, and we report pointwise confidence intervals based on Equation~\eqref{eq:cie}. Then, to evaluate the robustness of our conclusions, we also conduct sensitivity analyses under different bandwidths (undersmoothing and oversmoothing). As suggested by \citet{wasserman2006all}, adopting such a practical rule avoids the need to artificially eliminate asymptotic bias; instead, one acknowledges the presence of residual bias and addresses it through reporting and sensitivity analysis. Multiple robustness and cross-fitting reduce sensitivity to nuisance-model misspecification and nuisance-estimation error, but they do not remove the irregularity of $\psi_0(a,a')$; the practical implications of this irregularity are discussed in Section~7 of the Supplementary Material.

\section{Simulation Study}\label{simu}

We conducted a simulation study to demonstrate that the proposed estimator is consistent and multiply robust, and a sensitivity analysis to assess the uncertainty of the proposed estimator under different bandwidths and sample sizes. The data-generating process is as follows:
\begin{align*}
&{\bf X} = (X_1, X_2, X_3) \sim \mathcal{N}(0, \text{diag}\{0.25,0.1,0.8\}),\\
& A \sim \mathcal{N}(5+X_1+0.2 X_1^2,1),\\
&\delta(A, X) = \text{sigmoid}(-5+5A+2X_2+10AX_3),\\
&M \sim \text{Bernoulli}(\delta(A, X)), \\
& Y \sim \mathcal{N}( -A + 20M + 5MX_1 + X_2 ,1).
\end{align*}

The parameter of interest is $\psi_0(a, a')$ at $a = 4.5$ and $a' = 6$. Under the described simulation setting, the true parameter value is $9.1$, calculated based on Monte Carlo approximation of $\psi_0(a, a') = \int_{\mathcal{X}}\eta(a,a',X=x)dx$. To demonstrate the multiple robustness property, we considered various types of model misspecification in Table \ref{fig:sim-Table1}, where we also compared the proposed estimator $\hat{\psi}^{MR}(a, a')$ to the estimator in \citet{huber2020direct}, $\hat{\psi}^H(a, a')$, and the estimator without bias correction, $\hat{\psi}^{\eta}(a,a')$. To ensure comparability, we calculated the three estimators under the cross-fitting approach and defined $\hat{\psi}^H(a, a')$ and $\hat{\psi}^{\eta}(a, a')$ as
$$\hat{\psi}^H(a, a') = \frac{1}{L}\sum^L_{\ell = 1}\hat{\psi}^H_\ell(a, a')\quad\text{and}\quad  \hat{\psi}^\eta(a, a') = \frac{1}{L}\sum^L_{\ell = 1}\hat{\psi}^\eta_\ell(a, a'),$$
where
\begin{align*}
    \hat{\psi}^H_\ell(a, a') &=  \frac{1}{|I_\ell|}\frac{\sum_{i\in I_\ell} K_h(A_i - a)\hat{\lambda}(a,X_i)\frac{\hat{\alpha}(a',M_i,X_i)}{\hat{\alpha}(a,M_i,X_i)}Y_i}{\sum_{j\in I_\ell} K_h(A_j - a)\hat{\lambda}(a,X_j)\frac{\hat{\alpha}(a',M_j,X_j)}{\hat{\alpha}(a,M_j,X_j)}},\\
    \hat{\psi}^{\eta}_\ell(a,a') &= 
    \frac{1}{|I_\ell|}\sum_{i\in I_\ell} \hat{\eta}_\ell(a,a',X_i) = \frac{1}{|I_\ell|}\sum_{i\in I_\ell}\int_{\mathcal{M}} \hat{\gamma}_\ell(X_i, m, a) \hat{\alpha}_\ell(a',m,X_i)dm,
\end{align*}
 $\hat{\eta}(a,a',X) = \int_{\mathcal{M}} \hat{\gamma}(X, m, a) \hat{\alpha}(a',m,X)dm$, setting $L=3$. We used 1000 simulation replicates for each of the sample sizes 2000, 5000, and 8000, and chose the kernel bandwidth using the Silverman rule of thumb \citep{silverman2018density} under Gaussian kernels. 
The types of model misspecification considered include the scenario where all three models, $\mathbb{E}[Y|A,M,X]$, $f(M|A,X)$, and $f(A | X)$, are correctly specified, scenarios where only two out of the three models are correctly specified, and the scenario where all three models are misspecified.  The explicit nuisance model specifications corresponding to each column of Table~1 are provided in the Supplementary Material. As shown in Table \ref{fig:sim-Table1}, our proposed estimator has minimal or close to minimal bias for all scenarios except when all models are misspecified, demonstrating its theoretically proven multiple robustness property. The regression-only estimator \(\widehat{\psi}^{\eta}\) is sensitive to misspecification of the outcome model, and the generalized propensity-score weighted estimator \(\widehat{\psi}^{H}\) exhibits larger bias and RMSE in several misspecified settings. Additionally, bias and the root mean square error (RMSE) across simulation replicates reduce as the sample size gets larger, showing the consistency of our estimator. The bias becomes significant when all models are misspecified for all sample sizes and considered estimators. 

To further evaluate the finite-sample behavior of the estimators under varying degrees of treatment overlap, we conducted an additional overlap sensitivity analysis, reported in Supplementary Material Tables 2--9. In this analysis, we kept the reference exposure fixed at \(a'=6\) and evaluated the estimators over $a \in \{2,3,4,5,7,8,9,10\}$. In this analysis, as in Table 1, all estimators use the model-based conditional treatment density without the H\'ajek-type stabilization or clipping introduced for the application in Section~6.
Under the data-generating mechanism above, the continuous treatment distribution is centered near 5. Thus, values of \(a\) close to the center of the treatment distribution, such as \(a=4,5,7,8\), correspond to relatively strong or moderate overlap, whereas values farther in the tails, such as \(a=2,3,9,10\), represent increasingly poor overlap and more severe finite-sample positivity stress. The results show the expected deterioration in estimator stability as \(a\) moves into low-density regions of the treatment distribution. In the moderate-overlap region, estimators behavior is similar to that in Table 1. In the poor-overlap regions, all estimators become less stable, as reflected by inflated RMSEs, particularly for values of \(a\) deep in the tails. This behavior is expected because kernel smoothing around treatment values with limited empirical support relies on fewer effective observations and can amplify the effect of estimated inverse density weights. The proposed augmented estimator reduces bias relative to non-augmented alternatives in many settings where its nuisance-model requirements are satisfied, but it does not remove the fundamental finite-sample information loss caused by weak overlap.

\begin{table}[t]
\centering
\begin{tabular}{ccccccc}
\hline
\multicolumn{7}{c}{Absolute average bias (RMSE) when correct models are:}             \\
n    & Estimator & $(Y, M, A)$         & $(Y, M)$          & $(M, A)$          & $(Y, A)$          & None        \\ \hline
2000 & $\hat{\psi}^{MR}$ & 0.05 (0.23) & 0.04 (0.23) & 0.02 (0.43) & 0.04 (0.47) & 0.39 (0.54) \\ 
     & $\hat{\psi}^{\eta}$ & 0 (0.22) & 0 (0.22) & 0.54 (0.58) & 0.17 (0.28) & 0.63 (0.67) \\ 
     & $\hat{\psi}^{H}$ & 0.05 (0.51) & 0.39 (0.54) & 0.05 (0.51) & 0.19 (0.54) & 0.39 (0.54) \\ \hline
5000 & $\hat{\psi}^{MR}$ & 0.03 (0.15) & 0.02 (0.15) & 0.03 (0.16) & 0.04 (0.32) & 0.41 (0.5) \\ 
     & $\hat{\psi}^{\eta}$ & 0 (0.14) & 0 (0.14) & 0.54 (0.56) & 0.18 (0.22) & 0.63 (0.65) \\ 
     & $\hat{\psi}^{H}$ & 0 (0.34) & 0.41 (0.49) & 0 (0.34) & 0.23 (0.41) & 0.41 (0.49) \\ \hline
8000 & $\hat{\psi}^{MR}$ & 0.02 (0.12) & 0.02 (0.11) & 0.02 (0.12) & 0.01 (0.27) & 0.4 (0.45) \\ 
     & $\hat{\psi}^{\eta}$ & 0.01 (0.11) & 0.01 (0.11) & 0.53 (0.54) & 0.17 (0.2) & 0.62 (0.63) \\ 
     & $\hat{\psi}^{H}$ & 0.01 (0.26) & 0.4 (0.45) & 0.01 (0.26) & 0.22 (0.34) & 0.4 (0.45) \\  \hline
\end{tabular}
  \caption{Estimated absolute average bias (RMSE) of different estimators at $a = 4.5$ and $a' = 6$ averaged across 1000 simulation replicates under Silverman smoothing bandwidth, given sample size $n = 2000, 5000, $ and 8000.}
\label{fig:sim-Table1}
\end{table}


\begin{table}[t]
\centering
\footnotesize
\begin{tabular}{c| @{\hspace{0.7\tabcolsep}} c @{\hspace{0.7\tabcolsep}} c @{\hspace{0.7\tabcolsep}} c @{\hspace{0.7\tabcolsep}} c @{\hspace{0.7\tabcolsep}} c @{\hspace{0.7\tabcolsep}} c @{\hspace{0.7\tabcolsep}} c}
  \hline
&\multicolumn{7}{c}{Bandwidth}              \\
 & n & 0.1 & 0.2 & 0.3 & 0.4 & 0.5 & 0.6 \\ 
  \hline
Absolute Average & 2000 & 0.01 (0.26) & 0.03 (0.24) & 0.06 (0.24) & 0.1 (0.24) & 0.14 (0.26) & 0.19 (0.29) \\ 
 Bias (RMSE)  & 5000 & 0.01 (0.17) & 0.02 (0.15) & 0.06 (0.16) & 0.1 (0.17) & 0.14 (0.2) & 0.19 (0.24) \\ 
   & 8000 & 0 (0.13) & 0.02 (0.12) & 0.05 (0.12) & 0.09 (0.14) & 0.13 (0.17) & 0.18 (0.21) \\ \hline
 Mean (SD) of  & 2000 & 11.63 (2.32) & 10.63 (1.27) & 10.26 (0.8) & 10.06 (0.56) & 9.94 (0.41) & 9.86 (0.44) \\ 
  \scriptsize{$\sqrt{\hat{V}(a,a')}$}  & 5000 & 11.61 (0.98) & 10.58 (0.38) & 10.22 (0.25) & 10.03 (0.21) & 9.92 (0.18) & 9.84 (0.15) \\ 
   & 8000 & 11.56 (0.72) & 10.56 (0.29) & 10.21 (0.16) & 10.02 (0.11) & 9.91 (0.09) & 9.83 (0.09) \\ \hline
  Coverage & 2000 & 0.95 & 0.94 & 0.94 & 0.93 & 0.9 & 0.86 \\ 
  & 5000 & 0.94 & 0.94 & 0.93 & 0.9 & 0.82 & 0.72 \\ 
   & 8000 & 0.95 & 0.95 & 0.94 & 0.89 & 0.76 & 0.62 \\ 
   \hline
\end{tabular}
\caption{Sensitivity analysis over sample sizes and bandwidths under correct model specifications and $L=3$: absolute average bias (RMSE), average of $\hat{V}^{1/2}$, and coverage across 1000 simulation replicates under sample sizes n = 2000, 5000, 8000 and pre-specified bandwidths $h \in \{0.1, 0.2, \ldots, 0.6\}$.  }
\label{fig:sim-Table3}
\end{table}

Table \ref{fig:sim-Table3} summarizes the sensitivity of the estimator $\hat{\psi}^{MR}(a,a')$ under correct model specifications to different sample sizes and kernel smoothing bandwidths, by reporting absolute average bias, average of $\hat{V}(a,a')^{1/2}$, and Monte Carlo coverage. We define coverage as the proportion of simulation replicates that include the true value $\psi_0(a, a')$ in the estimated 95$\%$ confidence interval.  Recall that $\hat{V}(a, a')$ is the empirical variance of the estimating functions; we construct pointwise 95\% confidence intervals at significance level $\alpha = 0.05$ via a Wald approach as in Equation \eqref{eq:cie}. The Silverman kernel bandwidths for n = 2000, 5000, 8000, and under $L=3$ are approximately 0.28, 0.23, and 0.20. For the estimated mediation function $\hat{\psi}^{MR}(a,a')$, the presence of bias becomes apparent with reduced length of confidence intervals and decreasing coverage when the kernel bandwidths are larger than the Silverman-suggested optimal bandwidths, i.e., when the bandwidths are greater or equal to $0.3$. This pattern persists across different sample sizes. Coverage is theoretically guaranteed when the sample size $n$ goes to infinity and $\sqrt{nh^{d_A + 4}} \rightarrow 0$. We can see from Table \ref{fig:sim-Table3} that when bandwidth equals 0.1 (undersmoothed), coverage of the proposed pointwise 95\% confidence interval is indeed at least 0.95 for all sample sizes. In practice, choosing an undersmoothed kernel bandwidth can guarantee a relatively smaller bias in finite sample settings (for the price of conservative coverage). However, as previously discussed, there is no guidance on the choice of an undersmoothing bandwidth in practice.

\section{Application}\label{appli}

We applied the proposed approach to the Job Corps study \citep{huber2020direct,schochet2008does,schochet2001national}. Study participants were enrolled between 16 and 24 years old and from low-income households. The program provides eight months or approximately 1,200 hours of training on average. We aim to study the effect of the duration of Job Corp training $(a)$ on the binary outcome of the occurrence of any criminal arrests in the fourth year following program participation $(Y)$, with the proportion of weeks employed in the second year being the mediator $(M)$. Our study design follows \citet{huber2020direct}, who considered a similar causal mechanism and focused on the actual number of arrests in the fourth year as the outcome.

We consider a rich set of time-invariant socioeconomic variables as pre-treatment confounders $X$, similar to the study in \citep{huber2020direct}. Table 10 in Supplementary Material presents summary statistics of the following variables: outcome $Y$, mediator $M$, treatment $A$, and confounders $X$. Missing values in confounders are addressed by including the indicators of missingness as covariates. Moreover, following previous work on this dataset \citep{huber2020direct,flores2012estimating}, we applied our evaluation to the 
4,000 individuals in the dataset who received training, i.e., with a training duration in the program strictly greater than zero. Table 1 in Supplementary Material shows that on average,  5.1$\%$ of the participants had a history of imprisonment, and 23.75$\%$ had been arrested at least once before joining the study. Additionally, 8.7$\%$ of the individuals included in the study were arrested for criminal activities during the fourth year after study participation. 

We investigate whether longer Job Corps training reduces criminal behavior through employment or mechanisms beyond employment, by evaluating the natural direct and indirect effect of treatment durations at $a\in\{100,200,\ldots,2000\}$ hours versus just $a'=60$ hours, which corresponds to two weeks of training.  

Following \citet{huber2020direct}, we model the log-transformed treatment; the nuisance-model specifications are described below. To improve finite-sample stability of the localized inverse-propensity weights near treatment values with limited overlap, we used two modifications. First, we used a H\'ajek-type stabilized propensity density estimator \citep{hernan2020causal}, defined within each cross-fitting split as
\[
\hat f_H(a|X_i)
=
\hat f_0(a|X_i)
\left\{
\frac{1}{|I_{-\ell}|}
\sum_{j\in I_{-\ell}}
\frac{K_h(A_j-a)}{\hat f_0(a|X_j)}
\right\},
\]
where $\hat f_0(a|X)$ denotes the model-based conditional treatment density estimator and $I_{-\ell}$ denotes the training folds used to estimate the nuisance functions for observations in fold $\ell$. The estimator $\hat f_H$ is consistent whenever $\hat f_0$ is. Second, we used the clipped density estimator $\hat f_{H,c}(a|X_i)=\max\{\hat f_H(a|X_i),c\}$ with $c=0.01$ in place of $\hat f_H(a|X_i)$. Clipping is a finite-sample bias--variance trade-off. It can reduce the influence of observations with extremely large inverse-propensity weights, but it is not a substitute for the positivity assumption. The consistency of $\hat f_H$ and the properties of clipping are discussed in Section~8 of the Supplementary Material.

\begin{figure}[t!]
    \begin{center}
       \centering
    \includegraphics[scale = 0.5] {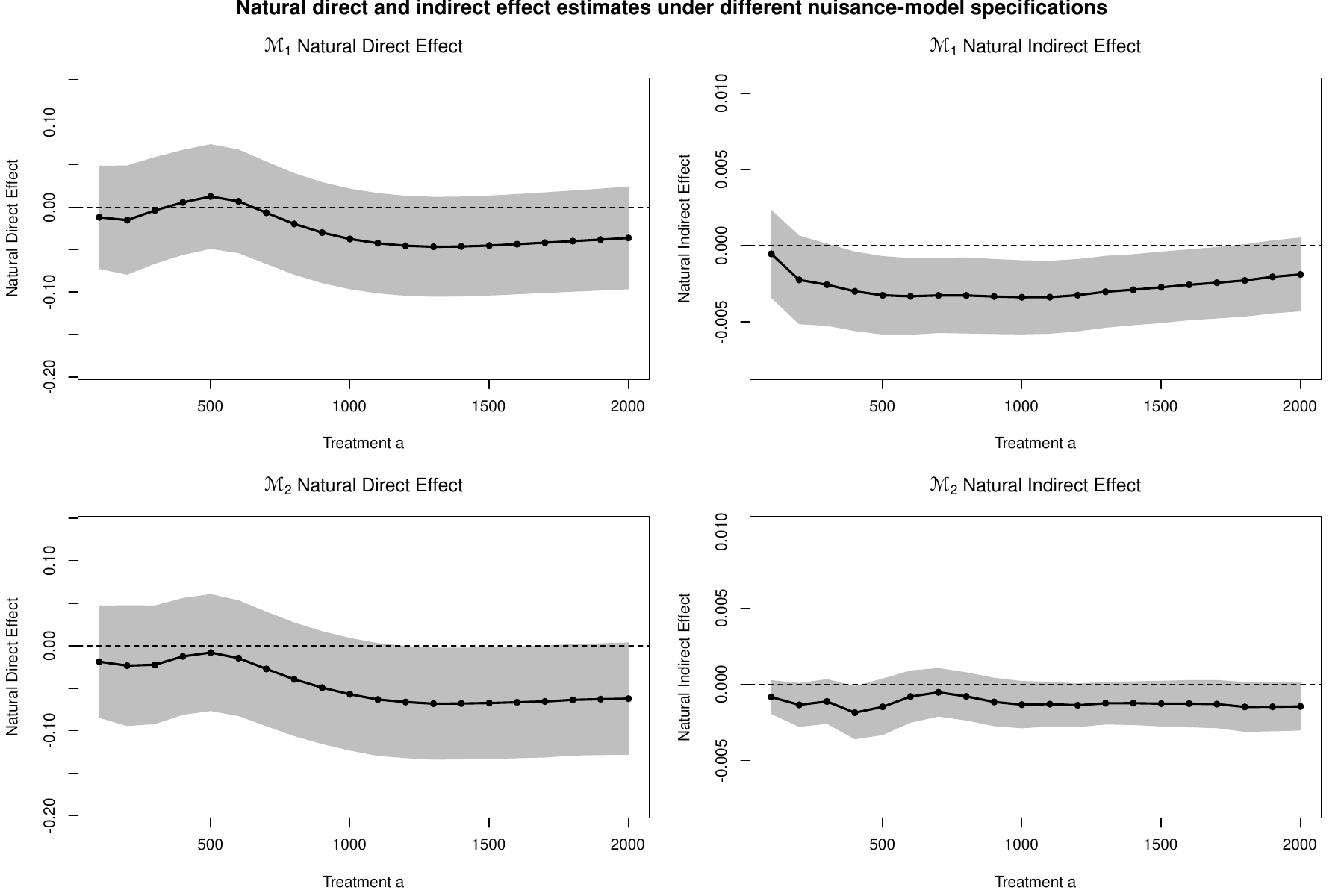}
   \end{center}
    \caption{Direct effect $\hat{\psi}^{MR}(a, a')-\hat{\psi}^{MR}(a', a')$ and indirect effect $\hat{\psi}^{MR}(a, a)-\hat{\psi}^{MR}(a, a')$ for $a'=60$ and $a\in \{100, 200, \ldots, 2000\}$ under Silverman bandwidth and clipping of the H\'ajek propensity at 0.01. $\mathcal M_1$ specifies parametric nuisance models: a logistic generalized linear model for the outcome regression, a Gaussian generalized linear model for the treatment density $f(A\mid X)$, and a beta regression model for the mediator density $f(M\mid A,X)$. $\mathcal M_2$ specifies RKHS Gaussian-kernel working models: kernel support vector machines for the outcome regression and working mean models, together with conditional mean embedding (CME) estimators for $f(A\mid X)$ and $f(\operatorname{logit}(M)\mid A,X)$, with the CME regularization parameter selected by internal three-fold cross-validation. For $\mathcal M_2$, mediator density ratios are evaluated on the logit-transformed mediator scale, where the logit Jacobian cancels in the ratio. The black line connects point estimates at the evaluated treatment values as a visual aid, and the grey shaded region connects the corresponding pointwise 95\% confidence interval limits.}
\label{fig:app-demo}
\end{figure}

\begin{figure}[t!]
    \begin{center}
   \includegraphics[scale = 0.5] {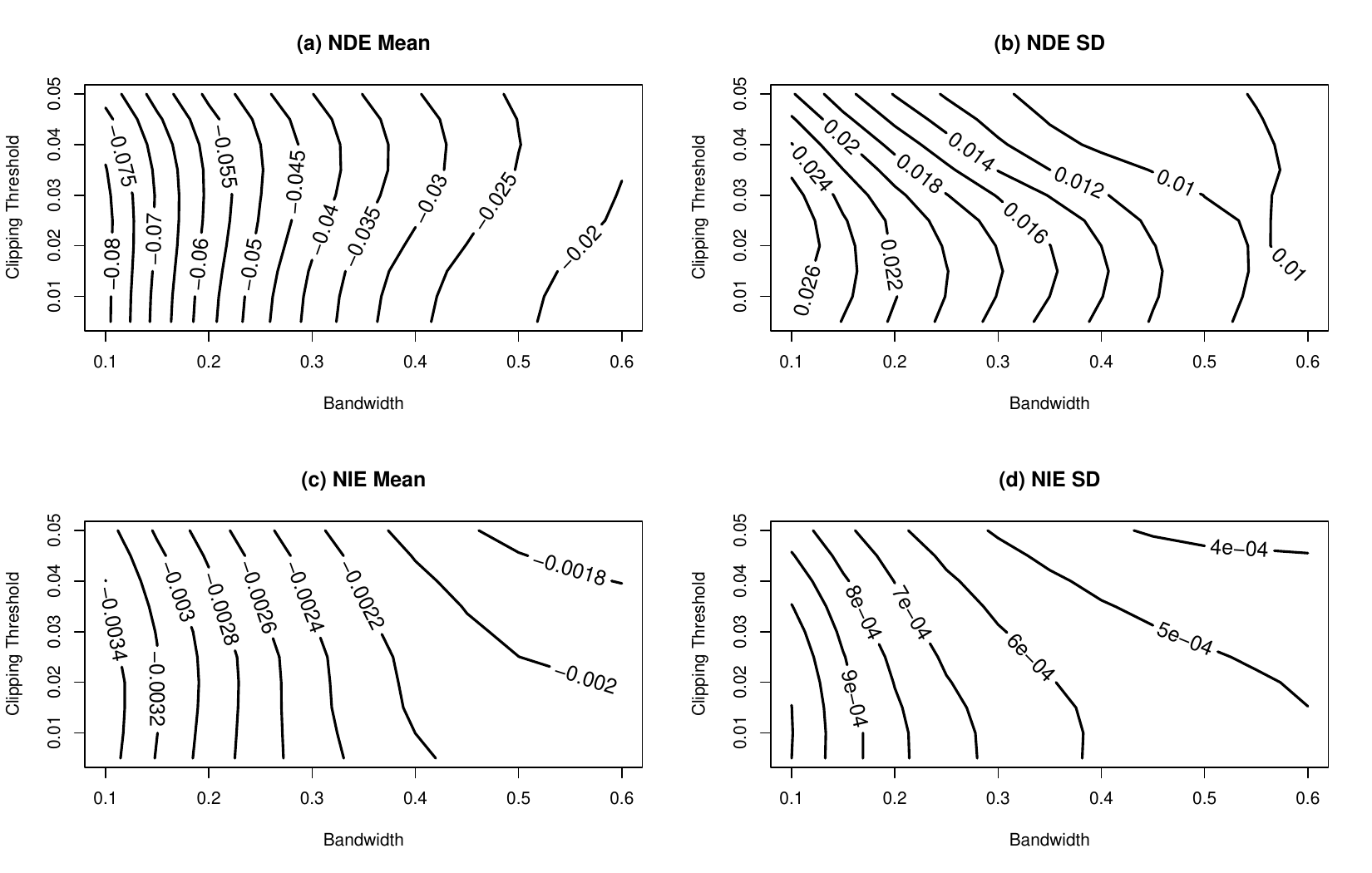}
   \end{center}
   \caption{ 
Sensitivity analysis for comparing treatments $a = 1500$ and $a' = 60$ using the proposed approach under different values of clipping threshold and bandwidth. In the application, the fold-specific Silverman rule-of-thumb bandwidths fall in the range \(0.23\)--\(0.25\), which is the reference bandwidth range for the main analysis. Bandwidths below approximately \(0.23\) correspond to undersmoothing, while bandwidths above approximately \(0.25\) correspond to oversmoothing. Panels show: (a) contour plot of the estimated mean of the natural direct effect; (b) contour plot of the estimated standard deviation of the natural direct effect; (c) contour plot of the estimated mean of the natural indirect effect; and (d) contour plot of the estimated standard deviation of the natural indirect effect.
   }
\label{fig:app-sensitivity}
\end{figure}

In practice, we treat clipping and bandwidth selection as sensitivity parameters. Useful diagnostics include summaries of the localized weights
\(
w_i(a)={K_h(A_i-a)}/{\hat f_{H,c}(a|X_i)}
\)
and the corresponding effective sample size
\(
n_{\mathrm{eff}}(a)
=
{\{\sum_i w_i(a)\}^2}/{\sum_i w_i^2(a)}.
\)
Large maximum weights, small effective sample sizes, or conclusions that change substantially across reasonable clipping thresholds should be interpreted as evidence of limited overlap for the target treatment contrast \citep{cole2008constructing,petersen2012diagnosing,ionides2008truncated}. Accordingly, Figure~\ref{fig:app-sensitivity} reports a sensitivity analysis over both smoothing bandwidths and clipping thresholds.

Figure \ref{fig:app-demo} displays the mean and 95\% confidence interval of natural direct and indirect effects over the range of values for $a$, under H\'ajek-type stabilized weighted propensities clipped at 0.01 \citep{ionides2008truncated} and Gaussian kernels with bandwidth chosen using the Silverman-type rule of thumb \citep{silverman2018density}. The confidence interval is obtained via Equation \eqref{eq:cie}. 
We fix the reference treatment level at $a'=60$ hours and compared it with evaluated treatment levels $a\in\{100,200,\ldots,2000\}$ hours. The first row reports results under $\mathcal M_1$, which uses parametric nuisance models consisting of a logistic generalized linear model for the outcome regression, a Gaussian generalized linear model for the log-transformed treatment density $f(\widetilde A\mid X)$, and a beta regression model for the shifted mediator density $f(M^*\mid \widetilde A,X)$, where $\widetilde A=\log(A)$, $A$ denotes raw treatment hours, and $M^*={[(n-1)(M/100)+0.5]}/n$ denotes the shifted mediator, with $M$ being the proportion of weeks employed in the second year recorded on a 0--100 scale and $n$ being the sample size. The transformation first rescales $M$ to $M/100$ and then maps it into the open unit interval so that beta-regression and logit-scale density estimation can be applied. The second row reports results under $\mathcal M_2$, which uses RKHS Gaussian-kernel working models, with kernel support vector machines for the outcome regression and working mean models, and RKHS conditional mean embedding (CME) estimators for $f(\widetilde A\mid X)$ and $f(\operatorname{logit}(M^*)\mid \widetilde A,X)$. The CME regularization parameter is selected by internal three-fold cross-validation.
 
 Figure \ref{fig:app-sensitivity} reports the sensitivity analysis and demonstrates how the estimated mean and empirical standard deviation vary by different clipping thresholds and smoothing bandwidths. For the Job Corps application, the fold-specific Silverman rule-of-thumb bandwidths are $(0.25, 0.24, 0.23)$. Therefore, in Figure 3, bandwidths smaller than approximately 0.23 represent undersmoothing relative to the data-driven choice, bandwidths between 0.23 and 0.25 represent the Silverman bandwidth region used in the main analysis, and bandwidths larger than approximately 0.25 represent oversmoothing. This sensitivity analysis shows how the estimated effects and standard deviations change when moving from undersmoothed to oversmoothed specifications. The source code is available at
\url{https://github.com/yizhenxu/Continuous-Treatment-Mediation.git}.

Both \citet{schochet2008does} and \citet{huber2020direct} identified significant effects of the training program in reducing criminal arrests, especially when the training duration is over 1000 hours. Similar to \citet{huber2020direct}, we observed nonlinear NDE and NIE of the training duration on the occurrence of arrests, demonstrated in Figure~\ref{fig:app-demo}. We magnified the NIE axis by a factor of 10 to highlight its presence, since the NIE is smaller in magnitude compared to the NDE. Under $\mathcal M_1$, we observed negative NIE estimates of the Job Corps training at $a'=60$ and $a$ between 400 and 1700 hours. The pointwise 95\% confidence intervals are below, or very close to, zero for most evaluated treatment levels, suggesting that part of the training effect may operate indirectly through employment under the parametric nuisance specification. Regarding NDE, although the mean is negative, the 95\% confidence interval is not entirely below zero. Under $\mathcal M_2$, the estimated NIEs are also mostly negative but closer to zero, and the corresponding confidence intervals generally overlap or lie very close to zero. Thus, the evidence for a nonzero indirect effect is attenuated under the more flexible RKHS CME nuisance specification. The NDE estimates under $\mathcal M_2$ are again mostly negative over moderate-to-large values of $a$, with the 95\% confidence interval falling below zero for $a$ between 1200 and 1600 but otherwise lying on or very close to the boundary of zero for larger $a$, providing only limited evidence of a statistically significant direct effect. Overall, the two specifications yield qualitatively similar negative point estimates for both NDE and NIE, while the indirect-effect signal is stronger under $\mathcal M_1$ than under $\mathcal M_2$.

Our sensitivity analysis in Figure~\ref{fig:app-sensitivity} indicates that while the overall results are fairly robust to the choice of bandwidth and clipping threshold, certain parameter settings can yield a statistically significant effect. For example, with oversmoothing bandwidth 0.3, the NDE mean and standard deviation are roughly -0.04 and 0.018, respectively, leading to an NDE that is significantly negative  for comparing $a=1500$ and $a'=60$. As the bandwidth increases from the undersmoothing region to the oversmoothing region, both the NDE and NIE means are attenuated toward zero, and their estimated standard deviations decrease. This reflects the usual bias–variance tradeoff in kernel smoothing: smaller bandwidths yield more localized but more variable estimates, whereas larger bandwidths yield smoother estimates with smaller variance but greater attenuation of the estimated effects. In addition, across the range of clipping thresholds considered, the qualitative pattern remains stable. Overall, these results suggest that extended training may reduce criminal behavior, with evidence of both direct effects and a smaller indirect pathway through employment, although the strength of the direct-effect evidence varies across bandwidth choices.

\section{Discussion}
In this paper, we proposed a multiply robust approach for estimating natural direct and indirect effects when the treatment is continuous. By replacing the treatment indicators that arise in the binary treatment influence function-based estimator of the mediation functional with kernel weights, the estimator targets pointwise cross-world potential outcome means at prespecified treatment values while retaining key advantages of influence-function based estimation, including robustness to certain nuisance-model misspecifications. The simulation studies support the multiple robustness and show stable performance when the required nuisance components are correctly specified in the relevant multiply robust submodels. The data application illustrates how the method can be used to summarize direct and indirect effect curves over scientifically meaningful treatment contrasts.

The comparison below summarizes how the reviewed binary-treatment framework changes when the treatment is continuous. With a continuous treatment, $P(A=a)=0$, so the identification formula in Equation~\eqref{eq:MCF} evaluates its components, $\gamma$ and $\alpha$, at treatment values of probability zero. The values of $\gamma$ and $\alpha$ at $a$ and $a'$ are therefore well defined only when these functions are continuous in the treatment argument. In addition, positivity requires the treatment density to be bounded away from zero near $a$ and $a'$. Both continuity and positivity are guaranteed by Assumption~\ref{assumption:reg}. Because $\psi_0(a,a')$ depends on the data distribution only through its behavior at the fixed points $a$ and $a'$, small perturbations of the distribution near these points can change $\psi_0(a,a')$ by arbitrarily large amounts; hence $\psi_0(a,a')$ is not pathwise differentiable in the nonparametric model and has no influence function. In Equation~\eqref{eq:IFb}, this is reflected in the indicators $I(A=a)$ and $I(A=a')$, which are zero almost surely when $A$ is continuous. The kernel weights used in Section~\ref{sec:c-med} replace these indicators by local averages over a window of width $h$, at the cost of a variance of order $h^{-d_A}$ and the slower $\sqrt{nh^{d_A}}$ rate established in Section~\ref{sec:asyp}.
\begin{center}
\small
\renewcommand{\arraystretch}{1.25}
\begin{tabular}{p{0.2\textwidth}p{0.35\textwidth}p{0.35\textwidth}}
 & Binary treatment, \(\mathcal A=\{0,1\}\) & Continuous treatment \\
\midrule
Target & A single contrast, e.g., \(\mathbb{E}[Y^{(1,M^{(0)})}]\) & The function \(\psi_0(a,a')\), \((a,a')\in\mathcal A\times\mathcal A\) \\
Positivity & \(P(A=a\mid X)>0\) & \(f_{A\mid X}(a\mid X)>0\), bounded away from zero near \(a\) (Assumption~\ref{assumption:reg}) \\
Identification & Formula evaluated at \(A=a\), with \(P(A=a)>0\) & Formula evaluated at fixed points \(a,a'\), with \(P(A=a)=0\); requires continuity in \(a\) (Assumption~\ref{assumption:reg}) \\
Influence function & Efficient influence function \eqref{eq:IFb} exists & Does not exist; \(\psi_0(a,a')\) is irregular \\
Estimation & Indicator weights \(I(A=a)\), \(I(A=a')\) & Kernel weights \(K_h(A-a)\), \(K_h(A-a')\) \\
Rate and bias & \(\sqrt{n}\); no smoothing bias & \(\sqrt{nh^{d_A}}\); smoothing bias \(h^2B(a,a')\) \\
Robustness & Consistent in \(\mathfrak{M}_{ym}\cup\mathfrak{M}_{ya}\cup\mathfrak{M}_{ma}\) & Same union model; rate conditions on pairwise products of nuisance errors \\
\end{tabular}
\end{center}

The continuous-treatment setting also highlights several practical limitations. First, pointwise mediation functionals are irregular in nonparametric models, so inference depends on the local effective sample size near the evaluated treatment values rather than on the full sample size. Consequently, bandwidth selection, treatment overlap, and stabilization of inverse-density weights play a central role in finite-sample performance. Second, although cross-fitting and multiple robustness reduce sensitivity to nuisance estimation error, they do not eliminate the need for adequate support near the target treatment values. In practice, the proposed estimator should therefore be accompanied by overlap diagnostics, sensitivity analyses over bandwidths and weight-stabilization choices, and cautious interpretation when the treatment values of interest lie in low-density regions of the observed treatment distribution.

Several extensions are important directions for future work. One natural extension is to develop theory and implementation for multivariate mediators. While the current framework formally allows general mediator spaces, direct estimation of joint conditional mediator densities may become difficult when the mediator is continuous or multivariate. Structured factorizations, conditional density-ratio estimation, representation learning, and dimension-reduction strategies may provide more scalable alternatives. Another direction is to study settings with high-dimensional outcomes or mediators, where the target may be a vector-valued mediation surface, a projection of such a surface, or a low-dimensional summary indexed by scientific priorities. 

Finally, another promising direction is a curve-estimation formulation based on pseudo-outcomes, analogous to Kennedy et al. \cite{kennedy2017nonparametric}. For a fixed reference treatment level \(a'\), one could in principle construct cross-fitted debiased pseudo-outcomes whose conditional mean, as a function of the observed treatment value \(A\), equals the mediation curve \(a \mapsto \psi_0(a,a')\). A second-stage nonparametric regression of such pseudo-outcomes on \(A\) would then estimate the full curve directly, rather than evaluating the pointwise kernel estimating equation developed here over a prespecified grid of treatment values. This is indeed a feasible and useful possibility, and it would provide a complementary route to estimating smooth mediation surfaces, facilitate data-adaptive smoothing over treatment values, and potentially improve scalability.
At the same time, this pseudo-outcome regression strategy represents a substantially different point of view from the pointwise estimation approach developed in the present paper. In particular, the mediation functional involves both an outcome-regression component and a mediator-distribution component evaluated at the reference level \(a'\). Therefore, preserving the multiple-robustness structure would require constructing a pseudo-outcome that appropriately debiases both components. In addition, such a formulation would require separate theoretical analysis of the second-stage regression, smoothing-parameter selection, and potentially simultaneous inference over the estimated treatment curve. We therefore leave the full development of this pseudo-outcome curve-estimation approach to future work.

\section*{Acknowledgments}
We are grateful to Professor Eric Tchetgen Tchetgen for insightful discussions and constructive comments that helped improve this work.

\bibliography{paper-ref.bib}
\bibliographystyle{apalike}

\clearpage

\begin{center}
{\LARGE \bfseries Supplementary Material}
\end{center}

\vspace{1em}

\setcounter{section}{0}
\setcounter{subsection}{0}
\setcounter{table}{0}

This supplement provides additional details and results supporting the main text. Section \textcolor{red}{1} discusses the identification assumptions in Assumption 1 of the main manuscript and their relation to the assumptions used for binary treatments. Section \textcolor{red}{2} describes the computation of the integral term \(\eta(a,a',X)\) for both binary and continuous mediators. Section \textcolor{red}{3} presents the nuisance model specifications used in the simulation study, including the correctly specified and misspecified models corresponding to Table 1 of the main manuscript. Section \textcolor{red}{4} reports additional simulation tables, analogous to Table 1 in the main manuscript, for different causal contrasts to illustrate estimator performance under varying degrees of treatment overlap. Section \textcolor{red}{5} contains the proofs of the main theoretical results in the manuscript. Section \textcolor{red}{6} includes a proof for the consistency of H\'ajek-type propensity estimator in cross validation. Section \textcolor{red}{7} shows a descriptive summary of the data used in our application.

\section{Discussion of the identification assumptions}
\label{sec:supp-id}

Assumption 1 in the main manuscript is the continuous-treatment analogue of the consistency, sequential ignorability, and positivity conditions used to identify natural direct and indirect effects under a binary treatment \citep{robins1992identifiability,pearl2001direct,imai2010identification,tchetgen2012semiparametric}. The discussion in this section applies equally to binary and continuous treatments; the aspects specific to continuous treatments are discussed in Section 2.1 of the main manuscript.

\textbf{Interpretation.} Assumption 1 gives a sufficient set of causal identification conditions for the natural direct and indirect effects \citep{robins1992identifiability,pearl2001direct,imai2010identification}. The consistency condition links the observed data to the relevant potential variables: if the observed treatment and mediator are equal to $(a,m)$, then the observed outcome is the corresponding potential outcome $Y^{(a,m)}$, and if the observed treatment is equal to $a$, then the observed mediator is the corresponding potential mediator $M^{(a)}$. The sequential exchangeability conditions require that the measured baseline covariates $X$ are rich enough to account for the confounding needed to identify the treatment--mediator, treatment--outcome, and mediator--outcome components of the natural-effect functional. The positivity condition requires sufficient overlap in the observed data: the treatment and mediator values appearing in the target estimand must occur with positive density, conditional on the relevant covariates. Without this support condition, the causal contrasts at the specified values of $a$, $a'$, and $m$ cannot be learned from the observed data without extrapolation.

\textbf{Relation to the binary-treatment assumptions.} \citet{tchetgen2012semiparametric} adopt the sequential ignorability assumption of \citet{imai2010identification}: for all treatment levels $e,e'$ and all $m\in\mathcal M$,
\[
\text{(i) } \{Y^{(e',m)},M^{(e)}\}\perp A\mid X,
\qquad
\text{(ii) } Y^{(e',m)}\perp M\mid A=e,X.
\]
With $e,e'$ ranging over $\mathcal A$, conditions (i) and (ii), together with consistency, imply the sequential exchangeability conditions in Assumption 1. Indeed, (i) implies $Y^{(a,m)}\perp A\mid X$ and $M^{(a)}\perp A\mid X$, and combining the former with (ii) gives $Y^{(a,m)}\perp \{A,M\}\mid X$. For the cross-world condition, (ii) with $e=a'$ and consistency, $M=M^{(a')}$ on the event $\{A=a'\}$, give $Y^{(a,m)}\perp M^{(a')}\mid A=a',X$; by (i), the joint conditional distribution of $(Y^{(a,m)},M^{(a')})$ given $(A=a',X)$ equals that given $X$, so that $Y^{(a,m)}\perp M^{(a')}\mid X$. Our formulation, which corresponds to the no-unmeasured-confounding conditions of \citet{pearl2001direct} and \citet{vanderweele2015explanation}, is therefore no stronger than the assumptions used in the binary-treatment setting. For continuous treatments it has the additional advantage that it does not condition on the null events $\{A=e\}$, as condition (ii) does.

\textbf{Cross-world exchangeability.} The condition
\(
Y^{(a,m)} \perp M^{(a^\prime)} \mid X
\)
is the cross-world exchangeability condition. It is called ``cross-world'' because it relates potential variables under two different hypothetical intervention regimes: the outcome that would be observed if treatment and mediator were set to $(a,m)$, and the mediator that would be observed if treatment were set to $a'$. This condition is not empirically testable from the observed data distribution alone and has been the subject of substantial discussion in the mediation literature \citep{robins1992identifiability,pearl2001direct,avin2005identifiability,vanderweele2015explanation}. Substantively, it rules out residual dependence, conditional on $X$, between the potential mediator under $a'$ and the potential outcome under the joint intervention $(a,m)$. In graphical terms, it is closely related to the absence of treatment-induced mediator--outcome confounders not included in $X$. In applications where the cross-world exchangeability condition is not substantively defensible, alternative estimands such as interventional or randomized interventional direct and indirect effects may be more appropriate \citep{vansteelandt2017interventional, diaz2020causal}; however, those estimands are outside the scope of the present work.

\textbf{Target-specific weakening.} Although the cross-world exchangeability is stated uniformly over all treatment contrasts
\((a,a')\in\mathcal A\times\mathcal A\), this is stronger than necessary
for a prespecified subset of contrasts. For example, if the reference
treatment level is fixed at \(a^\circ\) and interest is restricted to
\(\psi_0(a,a^\circ)\) over a prespecified set of values \(a\in\mathcal G\), the
cross-world condition need only hold as $Y^{(a,m)} \perp M^{(a^\circ)} \mid X,
 a\in\mathcal G,\; m\in\mathcal M$, together with the corresponding target-specific exchangeability and
positivity conditions. We state the condition uniformly to support the
general estimation framework for arbitrary treatment contrasts.

\section{Computation of integral $\hat\eta(a,a',X)$}
Recall that the nuisance function appearing in the estimating equation is
\[
\eta(a,a',X)
=
\int_{\mathcal M} \gamma(X,m,a)\alpha(a',m,X)\,dm,
\]
where $\gamma(X,M,a)=E(Y\mid A=a,M,X),
\alpha(a,M,X)=f_{M\mid A,X}(M\mid a,X)$. Thus, $\eta(a,a',X)$ is the conditional mean of the outcome regression evaluated at treatment level $a$,
after averaging over the conditional distribution of the mediator under treatment level $a'$.

Our implementation computes $\hat\eta(a,a',X)$ using the same fold-specific nuisance estimates used in
the construction of $\hat\psi^{MR}(a,a')$. Let $\{I_1,\ldots,I_L\}$ denote the cross-fitting folds. For
$i\in I_\ell$, the nuisance functions are estimated using only observations outside the $\ell$th fold, yielding
$\hat\gamma_\ell$ and $\hat\alpha_\ell$. We then compute
\[
\hat\eta_\ell(a,a',X_i)
=
\int_{\mathcal M}
\hat\gamma_\ell(X_i,m,a)
\hat\alpha_\ell(a',m,X_i)\,dm.
\]
Equivalently, if $M_i^*$ denotes a draw from the estimated conditional mediator distribution
\[
M_i^* \sim \hat\alpha_\ell(a',\cdot,X_i),
\]
then
\[
\hat\eta_\ell(a,a',X_i)
=
\mathbb E_{\hat\alpha_\ell(a',\cdot,X_i)}
\left\{
\hat\gamma_\ell(X_i,M_i^*,a)
\right\}.
\]
The numerical implementation differs according to whether the mediator is binary or continuous, as
described below.

\subsection*{Binary mediator}
In the simulation study, $M\in\{0,1\}$. Therefore, the integral defining $\eta(a,a',X)$ reduces to a finite
sum and is evaluated exactly. Let
\[
\hat p_{\ell,a'}(X_i)
=
\hat P_\ell(M=1\mid A=a',X=X_i).
\]
Then
\[
\hat\alpha_\ell(a',1,X_i)=\hat p_{\ell,a'}(X_i),
\qquad
\hat\alpha_\ell(a',0,X_i)=1-\hat p_{\ell,a'}(X_i),
\]
and hence, for $i\in I_\ell$,
\[
\hat\eta_\ell(a,a',X_i)
=
\hat\gamma_\ell(X_i,1,a)\hat p_{\ell,a'}(X_i)
+
\hat\gamma_\ell(X_i,0,a)\{1-\hat p_{\ell,a'}(X_i)\}.
\]
Thus, with a binary mediator, no Monte Carlo approximation is needed: the fitted outcome
regression is evaluated twice, once at $M=1$ and once at $M=0$, and these two fitted values are averaged
using the fitted mediator probability under $A=a'$.

\subsection*{Continuous mediator}
In the data application, the mediator is continuous. The raw treatment value is the number of Job Corps
training hours, and the nuisance models are fitted on the log-treatment scale $\tilde A=\log(A)$. Thus, when evaluating the nuisance functions at raw treatment levels $a$ and $a'$, the corresponding
values used in the fitted nuisance models are $\tilde a=\log(a)$ and $\tilde a'=\log(a')$.

The mediator is the proportion of weeks employed in the second year, recorded on a $0$--$100$ scale.
We first transform it to $M^*
=
\frac{N-1}{N}\frac{M}{100}+\frac{0.5}{N}$, so that $M^*\in(0,1)$. The integral defining $\eta$ is then evaluated on this shifted mediator scale:
\[
\hat\eta_\ell(a,a',X_i)
=
\int_0^1
\hat\gamma_\ell(X_i,m^*,\tilde a)
\hat\alpha_\ell(\tilde a',m^*,X_i)dm^*= \mathbb E_{\hat\alpha_\ell}
\left[
\hat\gamma_\ell(X_i,M^*,\tilde a)
\mid \tilde A=\tilde a',X_i
\right].
\]
Because this integral is not available in closed form for the continuous-mediator nuisance estimators, we
approximate it by Monte Carlo integration. Specifically, for each held-out observation $i\in I_\ell$, we
generate $R$ independent samples from the fitted conditional mediator distribution, $M_{ir}^*
\sim
\hat\alpha_\ell(\tilde a',\cdot,X_i), r=1,\ldots,R$, and compute the Monte
Carlo average
\[
\hat\eta^{MC}_\ell(a,a',X_i)
=
\frac{1}{R}
\sum_{r=1}^R
\hat\gamma_\ell(X_i,M_{ir}^*,\tilde a).
\]
Then by the law of large numbers, $\hat\eta^{MC}_\ell(a,a',X_i)$ asymptotically approaches $\hat\eta_\ell(a,a',X_i)$.  In our implementation, we use $R=1000$ Monte Carlo draws. We therefore use the Monte
Carlo average as the estimate of the integral.

\section{Alternative nuisance parametrizations}
\label{sec:supp-alt-param}

The moment function in Equation (4) of the main manuscript,
\[
\begin{aligned}
m(O;\alpha,\lambda,\gamma,\psi(a,a'))
&= K_h(A-a)\,\lambda(a,X)\,\frac{\alpha(a',M,X)}{\alpha(a,M,X)}\{Y-\gamma(X,M,a)\}\\
&\quad + K_h(A-a')\,\lambda(a',X)\{\gamma(X,M,a)-\eta(a,a',X)\}
+ \eta(a,a',X) - \psi(a,a'),
\end{aligned}
\]
involves the conditional mediator density
\(\alpha(a,m,x)=f_{M\mid A,X}(m\mid a,x)\), both through the density ratio
\(\alpha(a',M,X)/\alpha(a,M,X)\) and through the integral defining
\(\eta(a,a',X)\). In this section we describe two alternative parametrizations that avoid part or all of the direct estimation of \(\alpha\).

\textbf{Treatment-density ratio.} Following the same Bayes-rule argument used in
\citet{diaz2020causal} and \citet{farbmacher2022causal}, the mediator density ratio may be
rewritten in terms of conditional treatment densities. Whenever the relevant densities are well-defined and bounded away from
zero, $\alpha(a,m,x)
=
f_{A\mid M,X}(a\mid m,x) f_{M\mid X}(m\mid x)/f_{A\mid X}(a\mid x)$,
and hence
\[
\lambda(a,x)\frac{\alpha(a',m,x)}{\alpha(a,m,x)}
=
\frac{1}{f_{A\mid X}(a\mid x)}
\frac{
f_{A\mid M,X}(a'\mid m,x) f_{M\mid X}(m\mid x) / f_{A\mid X}(a'\mid x)
}{
f_{A\mid M,X}(a\mid m,x) f_{M\mid X}(m\mid x) / f_{A\mid X}(a\mid x)
}
=
\lambda(a',x)\frac{f_{A\mid M,X}(a'\mid m,x)}{f_{A\mid M,X}(a\mid m,x)} .
\]
Thus, the first term of the moment function can be equivalently expressed without
\(\alpha\) as
\[
K_h(A-a)\lambda(a',X)
\frac{f_{A\mid M,X}(a'\mid M,X)}{f_{A\mid M,X}(a\mid M,X)}
\{Y-\gamma(X,M,a)\}.
\]
This parametrization avoids the mediator density $\alpha$ but introduces a new nuisance component $\pi(a,m,x)=f_{A\mid M,X}(a\mid m,x)$, with corresponding regularity and product-rate conditions.

\textbf{Nested regression for \(\eta\).} The integral defining \(\eta(a,a',X)\) can be written as a nested
conditional mean:
\[
\eta(a,a',X)
=
\int_{\mathcal M} \gamma(X,m,a)\alpha(a',m,X)\,dm
=
E[\gamma(X,M,a)\mid A=a',X]:=\omega(a,a',X).
\]
In our empirical implementation, we estimate $\eta(a,a',X)$ using the original integral representation and approximate the integral by Monte Carlo draws from the fitted conditional mediator distribution, as described in Section 2. The $\omega(a,a',X)$ formulation suggests an alternative: form pseudo-outcomes $\hat\gamma(X_i,M_i,a)$ and estimate their conditional
mean given $(A,X)$, evaluated at $A=a'$. This nested-regression approach avoids explicit integration with
respect to the mediator density in the construction of $\eta$ and eliminates Monte Carlo error.

However, the nested-regression approach introduces an additional nuisance function $\omega(a,a',X)$ and
therefore an additional modeling or smoothing step. In particular, because $A$ is continuous, estimating
$\omega(a,a',X)$ requires pointwise prediction at $A=a'$, which may involve smoothing bias, tuning-parameter
selection, and overlap concerns near $a'$. Moreover, unless the first term of the moment function is also rewritten
using the treatment-density ratio above, the mediator density $\alpha$ is still needed for
$\alpha(a',M,X)/\alpha(a,M,X)$. Hence, a formulation that completely avoids direct estimation of mediator-density-related nuisances
requires the enlarged nuisance tuple $(\pi,\lambda,\gamma,\omega)$ rather than the original tuple $(\alpha,\lambda,\gamma)$, with corresponding regularity and product-rate
conditions for these nuisance components.

\section{Nuisance-model specifications in the simulation study}
\label{app:simulation-nuisance}

This section gives the nuisance-model specifications used in the simulation study presented in Section~5 Table 1. For all simulation settings, nuisance functions were estimated within the cross-fitting procedure. Specifically, for each fold, the nuisance models were fitted using the observations outside the fold and evaluated on the held-out fold. The number of folds was \(L=3\). The working models used to generate the columns of Table~1 in the manuscript are summarized in the Table~\ref{tab:simulation-nuisance-models} below. The misspecified versions omit key terms from these true nuisance mechanisms. In particular, the misspecified treatment model uses only an intercept; the misspecified mediator models either use only an intercept or omit \(X_2\); and the misspecified outcome models have more variations, omitting different terms needed to represent \(E(Y\mid A,M,X)\). These specifications were chosen to evaluate the multiple-robustness property by considering cases in which all three nuisance models are correct, exactly two nuisance models are correct, or none of the nuisance models is correct. 
\begin{table}[htbp]
\centering
\scriptsize
\renewcommand{\arraystretch}{1.25}
\caption{Nuisance-model specifications used in the simulation study.}
\label{tab:simulation-nuisance-models}
\begin{tabular}{p{0.18\textwidth}p{0.42\textwidth}p{0.28\textwidth}}
\hline
\multicolumn{3}{l}{\textbf{Outcome and treatment nuisance models}} \\
\hline
Correct models 
& Outcome model for \(\gamma\) 
& Treatment model for \(f(A\mid X)\) \\
\hline
\((Y,M,A)\)
&
\(Y \sim A + M + M:X_1 + X_2 - 1\)
&
\(A \sim X_1 + X_1^2\)
\\

\((Y,M)\)
&
\(Y \sim A + M + M:X_1 + X_2 - 1\)
&
\(A \sim 1\)
\\

\((M,A)\)
&
\(Y \sim M + X_1 + X_2\)
&
\(A \sim X_1 + X_1^2\)
\\

\((Y,A)\)
&
\(Y \sim A + M + M:X_1 + X_2 - 1\)
&
\(A \sim X_1 + X_1^2\)
\\

None
&
\(Y \sim A + X_2\)
&
\(A \sim 1\)
\\
\hline
\end{tabular}

\vspace{0.8em}

\begin{tabular}{p{0.18\textwidth}p{0.72\textwidth}}
\hline
\multicolumn{2}{l}{\textbf{Mediator nuisance model}} \\
\hline
Correct models 
& Mediator model for \(\alpha\) \\
\hline
\((Y,M,A)\)
&
\(\mathrm{logit}\,P(M=1\mid A,X) \sim A + X_2 + A:X_3\)
\\

\((Y,M)\)
&
\(\mathrm{logit}\,P(M=1\mid A,X) \sim A + X_2 + A:X_3\)
\\

\((M,A)\)
&
\(\mathrm{logit}\,P(M=1\mid A,X) \sim A + X_2 + A:X_3\)
\\

\((Y,A)\)
&
\(\mathrm{logit}\,P(M=1\mid A,X) \sim 1\)
\\

None
&
\(\mathrm{logit}\,P(M=1\mid A,X) \sim A + A:X_3\)
\\
\hline
\end{tabular}

\end{table}

\section{Additional simulation results under alternative causal contrasts}
This illustrates estimator performance across regions of poor, moderate, and strong treatment overlap.

\begin{table}[htp]
\centering
\begin{tabular}{ccccccc}
\hline
\multicolumn{7}{c}{Absolute average bias (RMSE) across varying levels of model misspecification} \\
n    & Estimator & YMA & YM & MA & YA & None \\ \hline
2000 & $\hat{\psi}^{MR}$ & 1.1 (22.54) & 0.26 (1.07) & 6.31 (134.37) & 0.32 (1.91) & 0.83 (8.27) \\
     & $\hat{\psi}^{\eta}$ & 0 (0.22) & 0 (0.22) & 3.04 (3.04) & 0.17 (0.28) & 2.76 (2.83) \\
     & $\hat{\psi}^{H}$ & 0.81 (2.96) & 1.84 (2.75) & 0.81 (2.96) & 1.91 (3.54) & 1.85 (2.75) \\ \hline
5000 & $\hat{\psi}^{MR}$ & 0.17 (1.77) & 0.14 (0.32) & 0.76 (5.74) & 0.17 (1.75) & 1.51 (2.31) \\
     & $\hat{\psi}^{\eta}$ & 0 (0.14) & 0 (0.14) & 3.04 (3.04) & 0.18 (0.23) & 2.77 (2.8) \\
     & $\hat{\psi}^{H}$ & 0.64 (2.55) & 1.76 (2.27) & 0.64 (2.55) & 1.79 (3.09) & 1.75 (2.26) \\ \hline
8000 & $\hat{\psi}^{MR}$ & 0.14 (0.6) & 0.11 (0.25) & 0.56 (2.51) & 0.12 (0.6) & 1.68 (2.17) \\
     & $\hat{\psi}^{\eta}$ & 0.01 (0.11) & 0.01 (0.11) & 3.03 (3.03) & 0.17 (0.2) & 2.75 (2.77) \\
     & $\hat{\psi}^{H}$ & 0.48 (2.33) & 1.83 (2.21) & 0.48 (2.33) & 1.6 (2.81) & 1.82 (2.2) \\ \hline
\end{tabular}
  \caption{Estimated absolute average bias (RMSE) of different estimators at $(a, a') = (2,6)$, averaged across 1000 simulation replicates under Silverman smoothing bandwidth, given sample size n = 2000, 5000, and 8000.}
\label{fig:append-Table-a2}
\end{table}

\begin{table}[htp]
\centering
\begin{tabular}{ccccccc}
\hline
\multicolumn{7}{c}{Absolute average bias (RMSE) across varying levels of model misspecification} \\
n    & Estimator & YMA & YM & MA & YA & None \\ \hline
2000 & $\hat{\psi}^{MR}$ & 0.18 (0.8) & 0.15 (0.3) & 0.38 (4.37) & 0.18 (0.57) & 1.07 (2.64) \\
     & $\hat{\psi}^{\eta}$ & 0.01 (0.22) & 0.01 (0.22) & 2.05 (2.06) & 0.18 (0.28) & 1.92 (1.97) \\
     & $\hat{\psi}^{H}$ & 0.12 (1.44) & 1.24 (1.46) & 0.12 (1.44) & 0.79 (1.64) & 1.25 (1.47) \\ \hline
5000 & $\hat{\psi}^{MR}$ & 0.13 (0.31) & 0.1 (0.19) & 0.13 (0.72) & 0.13 (0.41) & 1.2 (1.32) \\
     & $\hat{\psi}^{\eta}$ & 0.01 (0.14) & 0.01 (0.14) & 2.05 (2.05) & 0.18 (0.23) & 1.93 (1.95) \\
     & $\hat{\psi}^{H}$ & 0.04 (1.11) & 1.25 (1.36) & 0.04 (1.11) & 0.72 (1.33) & 1.24 (1.35) \\ \hline
8000 & $\hat{\psi}^{MR}$ & 0.1 (0.19) & 0.08 (0.15) & 0.08 (0.44) & 0.09 (0.3) & 1.25 (1.33) \\
     & $\hat{\psi}^{\eta}$ & 0 (0.11) & 0 (0.11) & 2.04 (2.04) & 0.18 (0.21) & 1.91 (1.92) \\
     & $\hat{\psi}^{H}$ & 0.07 (0.95) & 1.28 (1.35) & 0.07 (0.95) & 0.73 (1.19) & 1.27 (1.35) \\ \hline
\end{tabular}
  \caption{Estimated absolute average bias (RMSE) of different estimators at $(a, a') = (3,6)$, averaged across 1000 simulation replicates under Silverman smoothing bandwidth, given sample size n = 2000, 5000, and 8000.}
\label{fig:append-Table-a3}
\end{table}

\begin{table}[htp]
\centering
\begin{tabular}{ccccccc}
\hline
\multicolumn{7}{c}{Absolute average bias (RMSE) across varying levels of model misspecification} \\
n    & Estimator & YMA & YM & MA & YA & None \\ \hline
2000 & $\hat{\psi}^{MR}$ & 0.1 (0.28) & 0.07 (0.24) & 0.01 (0.93) & 0.09 (0.49) & 0.71 (0.85) \\
     & $\hat{\psi}^{\eta}$ & 0.01 (0.22) & 0.01 (0.22) & 1.05 (1.07) & 0.18 (0.28) & 1.06 (1.1) \\
     & $\hat{\psi}^{H}$ & 0.02 (0.7) & 0.71 (0.84) & 0.02 (0.7) & 0.35 (0.79) & 0.71 (0.84) \\ \hline
5000 & $\hat{\psi}^{MR}$ & 0.06 (0.17) & 0.05 (0.16) & 0.02 (0.35) & 0.07 (0.33) & 0.73 (0.79) \\
     & $\hat{\psi}^{\eta}$ & 0.01 (0.14) & 0.01 (0.14) & 1.05 (1.06) & 0.18 (0.23) & 1.07 (1.09) \\
     & $\hat{\psi}^{H}$ & 0.01 (0.46) & 0.73 (0.79) & 0.01 (0.46) & 0.36 (0.58) & 0.73 (0.79) \\ \hline
8000 & $\hat{\psi}^{MR}$ & 0.05 (0.13) & 0.04 (0.12) & 0.02 (0.16) & 0.04 (0.27) & 0.7 (0.74) \\
     & $\hat{\psi}^{\eta}$ & 0 (0.11) & 0 (0.11) & 1.04 (1.05) & 0.18 (0.21) & 1.06 (1.07) \\
     & $\hat{\psi}^{H}$ & 0.02 (0.39) & 0.7 (0.75) & 0.02 (0.39) & 0.32 (0.5) & 0.7 (0.74) \\ \hline
\end{tabular}
  \caption{Estimated absolute average bias (RMSE) of different estimators at $(a, a') = (4,6)$, averaged across 1000 simulation replicates under Silverman smoothing bandwidth, given sample size n = 2000, 5000, and 8000.}
\label{fig:append-Table-a4}
\end{table}

\begin{table}[htp]
\centering
\begin{tabular}{ccccccc}
\hline
\multicolumn{7}{c}{Absolute average bias (RMSE) across varying levels of model misspecification} \\
n    & Estimator & YMA & YM & MA & YA & None \\ \hline
2000 & $\hat{\psi}^{MR}$ & 0.02 (0.23) & 0.02 (0.23) & 0.02 (0.24) & 0.01 (0.46) & 0.08 (0.37) \\
     & $\hat{\psi}^{\eta}$ & 0.01 (0.22) & 0.01 (0.22) & 0.05 (0.22) & 0.18 (0.28) & 0.21 (0.3) \\
     & $\hat{\psi}^{H}$ & 0.01 (0.39) & 0.08 (0.37) & 0.01 (0.39) & 0.14 (0.42) & 0.08 (0.37) \\ \hline
5000 & $\hat{\psi}^{MR}$ & 0.02 (0.15) & 0.01 (0.15) & 0.01 (0.16) & 0.03 (0.32) & 0.09 (0.27) \\
     & $\hat{\psi}^{\eta}$ & 0.01 (0.14) & 0.01 (0.14) & 0.05 (0.15) & 0.18 (0.23) & 0.22 (0.26) \\
     & $\hat{\psi}^{H}$ & 0 (0.28) & 0.09 (0.27) & 0 (0.28) & 0.14 (0.31) & 0.09 (0.27) \\ \hline
8000 & $\hat{\psi}^{MR}$ & 0.01 (0.12) & 0.01 (0.11) & 0.01 (0.12) & 0 (0.27) & 0.08 (0.22) \\
     & $\hat{\psi}^{\eta}$ & 0 (0.11) & 0 (0.11) & 0.04 (0.11) & 0.18 (0.21) & 0.21 (0.23) \\
     & $\hat{\psi}^{H}$ & 0 (0.22) & 0.08 (0.23) & 0 (0.22) & 0.13 (0.26) & 0.08 (0.22) \\ \hline
\end{tabular}
  \caption{Estimated absolute average bias (RMSE) of different estimators at $(a, a') = (5,6)$, averaged across 1000 simulation replicates under Silverman smoothing bandwidth, given sample size n = 2000, 5000, and 8000.}
\label{fig:append-Table-a5}
\end{table}

\begin{table}[htp]
\centering
\begin{tabular}{ccccccc}
\hline
\multicolumn{7}{c}{Absolute average bias (RMSE) across varying levels of model misspecification} \\
n    & Estimator & YMA & YM & MA & YA & None \\ \hline
2000 & $\hat{\psi}^{MR}$ & 0.14 (0.3) & 0.1 (0.26) & 0.15 (0.48) & 0.15 (0.5) & 1.36 (1.64) \\
     & $\hat{\psi}^{\eta}$ & 0.01 (0.22) & 0.01 (0.22) & 1.95 (1.97) & 0.18 (0.28) & 1.49 (1.56) \\
     & $\hat{\psi}^{H}$ & 0 (1.01) & 1.37 (1.61) & 0 (1.01) & 0.12 (1) & 1.38 (1.62) \\ \hline
5000 & $\hat{\psi}^{MR}$ & 0.1 (0.2) & 0.07 (0.17) & 0.1 (0.32) & 0.09 (0.34) & 1.39 (1.52) \\
     & $\hat{\psi}^{\eta}$ & 0.01 (0.14) & 0.01 (0.14) & 1.95 (1.96) & 0.18 (0.23) & 1.5 (1.52) \\
     & $\hat{\psi}^{H}$ & 0.04 (0.72) & 1.4 (1.51) & 0.04 (0.72) & 0.15 (0.72) & 1.4 (1.51) \\ \hline
8000 & $\hat{\psi}^{MR}$ & 0.08 (0.16) & 0.06 (0.14) & 0.07 (0.26) & 0.09 (0.3) & 1.4 (1.5) \\
     & $\hat{\psi}^{\eta}$ & 0 (0.11) & 0 (0.11) & 1.96 (1.96) & 0.18 (0.21) & 1.5 (1.51) \\
     & $\hat{\psi}^{H}$ & 0.01 (0.61) & 1.4 (1.5) & 0.01 (0.61) & 0.11 (0.62) & 1.4 (1.5) \\ \hline
\end{tabular}
  \caption{Estimated absolute average bias (RMSE) of different estimators at $(a, a') = (7,6)$, averaged across 1000 simulation replicates under Silverman smoothing bandwidth, given sample size n = 2000, 5000, and 8000.}
\label{fig:append-Table-a7}
\end{table}

\begin{table}[htp]
\centering
\begin{tabular}{ccccccc}
\hline
\multicolumn{7}{c}{Absolute average bias (RMSE) across varying levels of model misspecification} \\
n    & Estimator & YMA & YM & MA & YA & None \\ \hline
2000 & $\hat{\psi}^{MR}$ & 0.28 (0.87) & 0.2 (0.4) & 0.99 (3) & 0.27 (0.92) & 2.11 (3.74) \\
     & $\hat{\psi}^{\eta}$ & 0.01 (0.22) & 0.01 (0.22) & 2.95 (2.96) & 0.18 (0.28) & 2.35 (2.43) \\
     & $\hat{\psi}^{H}$ & 0.4 (2.78) & 2.21 (3.13) & 0.4 (2.78) & 0.56 (2.78) & 2.22 (3.13) \\ \hline
5000 & $\hat{\psi}^{MR}$ & 0.19 (0.51) & 0.13 (0.26) & 0.54 (1.71) & 0.18 (0.57) & 2.18 (2.85) \\
     & $\hat{\psi}^{\eta}$ & 0.01 (0.14) & 0.01 (0.14) & 2.95 (2.96) & 0.18 (0.23) & 2.35 (2.39) \\
     & $\hat{\psi}^{H}$ & 0.44 (2.15) & 2.27 (2.72) & 0.44 (2.15) & 0.6 (2.14) & 2.28 (2.72) \\ \hline
8000 & $\hat{\psi}^{MR}$ & 0.14 (0.39) & 0.1 (0.2) & 0.43 (1.44) & 0.15 (0.46) & 2.24 (2.66) \\
     & $\hat{\psi}^{\eta}$ & 0 (0.11) & 0 (0.11) & 2.96 (2.96) & 0.18 (0.21) & 2.35 (2.37) \\
     & $\hat{\psi}^{H}$ & 0.31 (1.82) & 2.29 (2.59) & 0.31 (1.82) & 0.44 (1.84) & 2.29 (2.58) \\ \hline
\end{tabular}
  \caption{Estimated absolute average bias (RMSE) of different estimators at $(a, a') = (8,6)$, averaged across 1000 simulation replicates under Silverman smoothing bandwidth, given sample size n = 2000, 5000, and 8000.}
\label{fig:append-Table-a8}
\end{table}

\begin{table}[htp]
\centering
\begin{tabular}{ccccccc}
\hline
\multicolumn{7}{c}{Absolute average bias (RMSE) across varying levels of model misspecification} \\
n    & Estimator & YMA & YM & MA & YA & None \\ \hline
2000 & $\hat{\psi}^{MR}$ & 0.64 (3.31) & 0.36 (2.28) & 2.5 (15.77) & 0.64 (3.32) & 1.86 (23.77) \\
     & $\hat{\psi}^{\eta}$ & 0.01 (0.22) & 0.01 (0.22) & 3.95 (3.96) & 0.18 (0.29) & 3.2 (3.3) \\
     & $\hat{\psi}^{H}$ & 1.91 (5.37) & 3.03 (6.03) & 1.91 (5.37) & 2.01 (5.39) & 3.03 (6.03) \\ \hline
5000 & $\hat{\psi}^{MR}$ & 0.16 (3.06) & 0.24 (0.84) & 1.85 (12.13) & 0.13 (3.32) & 2.77 (11.05) \\
     & $\hat{\psi}^{\eta}$ & 0.01 (0.14) & 0.01 (0.14) & 3.95 (3.96) & 0.18 (0.23) & 3.21 (3.25) \\
     & $\hat{\psi}^{H}$ & 2.09 (5.24) & 3.36 (5.88) & 2.09 (5.24) & 2.19 (5.26) & 3.35 (5.88) \\ \hline
8000 & $\hat{\psi}^{MR}$ & 0.1 (2.61) & 0.18 (0.72) & 2.25 (25.47) & 0.1 (2.61) & 3.03 (9.14) \\
     & $\hat{\psi}^{\eta}$ & 0 (0.11) & 0 (0.11) & 3.96 (3.96) & 0.18 (0.21) & 3.2 (3.23) \\
     & $\hat{\psi}^{H}$ & 2.03 (4.77) & 3.42 (5.44) & 2.03 (4.77) & 2.13 (4.78) & 3.42 (5.45) \\ \hline
\end{tabular}
  \caption{Estimated absolute average bias (RMSE) of different estimators at $(a, a') = (9,6)$, averaged across 1000 simulation replicates under Silverman smoothing bandwidth, given sample size n = 2000, 5000, and 8000.}
\label{fig:append-Table-a9}
\end{table}

\begin{table}[htp]
\centering
\begin{tabular}{ccccccc}
\hline
\multicolumn{7}{c}{Absolute average bias (RMSE) across varying levels of model misspecification} \\
n    & Estimator & YMA & YM & MA & YA & None \\ \hline
2000 & $\hat{\psi}^{MR}$ & 0.67 (15.77) & 1.36 (9.58) & 13.16 (237.18) & 0.66 (15.75) & 7.16 (122.57) \\
     & $\hat{\psi}^{\eta}$ & 0.01 (0.22) & 0.01 (0.22) & 4.95 (4.96) & 0.18 (0.29) & 4.05 (4.17) \\
     & $\hat{\psi}^{H}$ & 3.6 (6.92) & 4.21 (7.43) & 3.6 (6.92) & 3.66 (6.94) & 4.21 (7.43) \\ \hline
5000 & $\hat{\psi}^{MR}$ & 0.77 (28.75) & 0.66 (6.75) & 12.76 (266.83) & 0.8 (29.65) & 5.14 (78.94) \\
     & $\hat{\psi}^{\eta}$ & 0.01 (0.14) & 0.01 (0.14) & 4.95 (4.95) & 0.18 (0.23) & 4.06 (4.11) \\
     & $\hat{\psi}^{H}$ & 3.98 (7.4) & 4.51 (7.77) & 3.98 (7.4) & 4.02 (7.43) & 4.51 (7.77) \\ \hline
8000 & $\hat{\psi}^{MR}$ & 0.27 (8.03) & 0.09 (5.02) & 0.27 (37.2) & 0.28 (8.04) & 1.18 (71.63) \\
     & $\hat{\psi}^{\eta}$ & 0 (0.11) & 0 (0.11) & 4.96 (4.96) & 0.18 (0.21) & 4.05 (4.08) \\
     & $\hat{\psi}^{H}$ & 3.87 (7.26) & 4.23 (7.59) & 3.87 (7.26) & 3.9 (7.27) & 4.23 (7.59) \\ \hline
\end{tabular}
  \caption{Estimated absolute average bias (RMSE) of different estimators at $(a, a') = (10,6)$, averaged across 1000 simulation replicates under Silverman smoothing bandwidth, given sample size n = 2000, 5000, and 8000.}
\label{fig:append-Table-a10}
\end{table}

\newpage
\section{Proofs}
Before we start with the proofs, we establish some lemmas that will help us with the proofs in the rest of the Supplementary Material. 


\begin{lemma}
\label{lem:bounding-term}

Let $\{X_m\}$ and $\{Y_m\}$ be a sequence of random variables. Then under conditions outlined in Lemma 6.1 in \cite{chernozhukov2018double}, $\mathbb{E}[|X_m| \mid Y_m] = o_p(1)$ implies $X_m = o_p(1)$.
\end{lemma}

\begin{proof}

By the Conditional Markov Inequality, for any $\epsilon > 0$,
\begin{align*}
    p(|X_m| \geq \epsilon \mid Y_m) \leq \frac{\mathbb{E}[|X_m| \mid Y_m]}{\epsilon}
\end{align*}

By $\mathbb{E}[|X_m| \mid Y_m] = o_p(1)$, there is $p(|X_m| \geq \epsilon \mid Y_m) = o_p(1)$. An application of Lemma 6.1 then yields $p(|X_m| > \epsilon) \rightarrow 0$, therefore $X_m = o_p(1)$.
\end{proof}

\begin{lemma}
\label{lem:kernelexpansion}
Under Assumption 2, for a twice continuously differentiable function $f$ with bounded first and second derivative, we have

\[
\int_A K_h(A-a)f(A)dA=f(a)+O(h^2).
\]

\end{lemma}
\begin{proof}

\begin{align*}
&\int_A K_h(A-a)f(A)dA\\
&=\int \left[\prod_{j = 1}^{d_A}k(u_j)\right]f(uh+a)du_1\dots du_{d_A}\\
&=\int \left[\prod_{j = 1}^{d_A}k(u_j)\right]\Bigg\{f(a)+\sum_{j = 1}^{d_A} u_jh \partial_{a_j}f(a)+ 
\frac{1}{2}\sum_{j = 1}^{d_A} \sum_{j^\prime = 1}^{d_A} u_ju_{j^\prime}h^2\partial_{a_j}\partial_{a_{j^\prime}}\partial_{a_{j^{\prime}}}f(a)|_{\Bar{a}}\Bigg\}\\
&\qquad\qquad du_1 \dots du_{d_A}\\
&=f(a)+O(h^2),
\end{align*}
where $\Bar{a}$ is in between $A$ and $a$. 
\end{proof}

Remark: We assume the second derivative is bounded over the support of the function $f(a)$, which is a stronger assumption than $O(1)$ since the bound holds everywhere as opposed to only for $a \geq c$ where $c$ is a constant. If $\nu(x)$ and $\omega(x)$ are two arbitrary functions, then $\int \nu(x) \omega(x) dx = O(1) \int |\omega(x)|dx$ is true when $\nu(x)$ is bounded, but not when $\nu(x) = O(1)$, e.g. when $\nu(x) = 1/x$ and $\omega(x) = \mathbb{I}\{0\le x \le 1\}$ .

\subsection{Proof of Theorem 1}

We follow a similar outline as \cite{colangelo2020double} and \cite{chernozhukov2018double}. The proof for this theorem is split into two parts. The first part establishes that the proposed estimator satisfies
\begin{align*}
    \sqrt{\frac{h^{d_A}}{n}}\sum^L_{\ell=1} \sum_{i \epsilon I_\ell} \left\{m(O_i;\hat{\alpha}_\ell,\hat{\lambda}_\ell,\hat{\gamma}_\ell,\psi_0(a,a')) - m(O_i;\alpha,\lambda,\gamma,\psi_0(a,a')) \right\} = o_p(1),
\end{align*} and the second part establishes that $\sqrt{nh^{d_A}}(\hat{\psi}^{MR}(a,a')-\psi_0(a,a')-B(a,a'))$ converges to the Gaussian distribution $\mathcal{N}(0,V(a,a'))$. 

Starting with the first part of the proof, note that
\begin{align*}
    &\sqrt{nh^{d_A}} \frac{1}{n}\sum^L_{\ell=1}\sum_{i\in I_\ell} \Big\{m(O_i; \hat{\alpha}_\ell,\hat{\lambda}_\ell, \hat{\gamma}_\ell, \hat{\psi}_\ell(a,a'))- m(O_i;\alpha,\lambda,\gamma,\psi_0(a,a'))\Big\}\\ 
    =&  \sqrt{\frac{h^{d_A}}{n}}\sum^L_{\ell=1}\sum_{i\in I_\ell} \Big\{m(O_i; \hat{\alpha}_\ell,\hat{\lambda}_\ell, \hat{\gamma}_\ell, \hat{\psi}_\ell(a,a'))-m(O_i;\hat{\alpha}_\ell,\hat{\lambda}_\ell,\hat{\gamma}_\ell,\psi_0(a,a'))\\
    &+m(O_i;\hat{\alpha}_\ell,\hat{\lambda}_\ell,\hat{\gamma}_\ell,\psi_0(a,a'))- m(O_i;\alpha,\lambda,\gamma,\psi_0(a,a'))\Big\}\\ 
    =& -\sqrt{nh^{d_A}} (\hat{\psi}^{MR}(a,a')- \psi_0(a,a') )\\ &+\sqrt{\frac{h^{d_A}}{n}}\sum^L_{\ell=1}\sum_{i\in I_\ell}\Big\{m(O_i;\hat{\alpha}_\ell,\hat{\lambda}_\ell,\hat{\gamma}_\ell,\psi_0(a,a'))- m(O_i;\alpha,\lambda,\gamma,\psi_0(a,a'))\Big\}.
\end{align*}
Since $\frac{1}{n}\sum^L_{\ell=1}\sum_{i\in I_\ell} m(O_i; \hat{\alpha}_\ell,\hat{\lambda}_\ell, \hat{\gamma}_\ell, \hat{\psi}_\ell(a,a'))=0$, we have
\begin{align*}
    &\sqrt{nh^{d_A}} (\hat{\psi}^{MR}(a,a')- \psi_0(a,a') )\\
    =& \sqrt{\frac{h^{d_A}}{n}}\sum^L_{\ell=1}\sum_{i\in I_\ell} \Big\{ m(O_i;\alpha,\lambda,\gamma,\psi_0(a,a'))\Big\}\\
    &+\sqrt{\frac{h^{d_A}}{n}}\sum^L_{\ell=1}\sum_{i\in I_\ell}\Big\{m(O_i;\hat{\alpha}_\ell,\hat{\lambda}_\ell,\hat{\gamma}_\ell,\psi_0(a,a'))- m(O_i;\alpha,\lambda,\gamma,\psi_0(a,a'))\Big\}.
\end{align*}
In order to establish an asymptotically linear representation for our proposed estimator, it suffices to to show that for all $1\le \ell\le L$ we have
\begin{align*}
    \sqrt{\frac{h^{d_A}}{n}} \sum_{i \epsilon I_\ell} \left\{m(O_i;\hat{\alpha}_\ell,\hat{\lambda}_\ell,\hat{\gamma}_\ell,\psi_0(a,a')) - m(O_i;\alpha,\lambda,\gamma,\psi_0(a,a')) \right\} = o_p(1).
\end{align*}

Next, we expand $m(O_i;\hat{\alpha}_\ell,\hat{\lambda}_\ell,\hat{\gamma}_\ell,\psi_0(a,a')) - m(O_i;\alpha,\lambda,\gamma,\psi_0(a,a'))$ into multiple terms and bound each term individually. Note that
\begin{align}
&m(O_i;\hat{\alpha}_\ell,\hat{\lambda}_\ell,\hat{\gamma}_\ell,\psi_0(a,a')) - m(O_i;\alpha,\lambda,\gamma,\psi_0(a,a'))\nonumber\\
=&K_h(A_i - a)\big\{\hat{\lambda}(a,X_i)\frac{\hat{\alpha}(a',M_i,X_i)}{\hat{\alpha}(a,M_i,X_i)}[Y_i - \hat{\gamma}(a,M_i,X_i)]\notag\\
&\quad-\lambda(a,X_i)\frac{\alpha(a',M_i,X_i)}{\alpha(a,M_i,X_i)}[Y_i - \gamma(X_i,M_i,a)]\big\}\label{eq:term1}\\
&\quad+K_h(A_i - a')\big\{\hat{\lambda}(a',X_i)[\hat{\gamma}(a,M_i,X_i) 
- \hat{\eta}(a, a', X_i)]\notag\\
&\quad-\lambda(a',X_i)[\gamma(X_i,M_i,a) - \eta(a, a', X_i)]\big\}\label{eq:term2}\\
&\quad+\hat{\eta}(a, a', X_i)-\eta(a, a', X_i)\tag{R1}.
\end{align}

Defining $R(M_i,X_i) := \frac{\alpha(a^\prime, M_i, X_i)}{\alpha(a, M_i, X_i)}$, terms \eqref{eq:term1} and \eqref{eq:term2} can be expanded additionally. Expanding term \eqref{eq:term1}, we get
\begin{align*}
&K_h(A_i - a)\big\{\hat{\lambda}(a, X_i)\hat{R}(M_i,X_i)\{Y_i - \hat{\gamma}(X_i,M_i,a)\}\\
&\qquad\qquad-\lambda(a,X_i)R(M_i,X_i)\{Y_i - \gamma(X_i,M_i,a)\}\big\}\\
&=-K_h(A_i - a)\big(\hat{R}(M_i,X_i)-R(M_i,X_i)\big)\big(\hat{\lambda}(a, X_i)-\lambda(a,X_i)\big)\big(\hat{\gamma}(X_i,M_i,a)-\gamma(X_i,M_i,a)\big)\tag{CS1}\\
&\quad+K_h(A_i - a)\big(\hat{R}(M_i,X_i)-R(M_i,X_i)\big)\big(\hat{\gamma}(X_i,M_i,a)-\gamma(X_i,M_i,a)\big)\big(Y_i - \gamma(X_i,M_i,a)\big)\tag{CS2}\\
&\quad-K_h(A_i - a)\big(\hat{R}(M_i,X_i)-R(M_i,X_i)\big)\big(\hat{\gamma}(X_i,M_i,a)-\gamma(X_i,M_i,a)\big)\lambda(a,X_i)\tag{CS3}\\
&\quad-K_h(A_i - a)\big(\hat{\lambda}(a, X_i)-\lambda(a,X_i)\big)\big(\hat{\gamma}(X_i,M_i,a)-\gamma(X_i,M_i,a)\big)R(M_i,X_i)\tag{CS4}\\
&\quad+\Big\{K_h(A_i - a)\big(\hat{R}(M_i,X_i)-R(M_i,X_i)\big)\lambda(a,X_i)\big(Y_i - \gamma(X_i,M_i,a)\big)\\
&\quad\quad\quad-\mathbb{E}\big[
K_h(A_i - a)\big(\hat{R}(M_i,X_i)-R(M_i,X_i)\big)\lambda(a,X_i)\big(Y_i - \gamma(X_i,M_i,a)\big) \mid  O^c_{I_\ell} \big]\Big\}\tag{E1}\\
&\quad+\mathbb{E}\big[
K_h(A_i - a)\big(\hat{R}(M_i,X_i)-R(M_i,X_i)\big)\lambda(a,X_i)\big(Y_i - \gamma(X_i,M_i,a)\big) \mid  O^c_{I_\ell}\big]\tag{TR1}\\
&\quad+\Big\{K_h(A_i - a)\big(\hat{\lambda}(a, X_i)-\lambda(a,X_i)\big)R(M_i,X_i)\big(Y_i - \gamma(X_i,M_i,a)\big)\\
&\quad\quad\quad-\mathbb{E}\big[
K_h(A_i - a)\big(\hat{\lambda}(a, X_i)-\lambda(a,X_i)\big)R(M_i,X_i)\big(Y_i - \gamma(X_i,M_i,a)\big) \mid  O^c_{I_\ell}\big]\Big\}\tag{E2}\\
&\quad+\mathbb{E}\big[
K_h(A_i - a)\big(\hat{\lambda}(a, X_i)-\lambda(a,X_i)\big)R(M_i,X_i)\big(Y_i - \gamma(X_i,M_i,a)\big) \mid  O^c_{I_\ell} \big]\tag{TR2}\\
&\quad-K_h(A_i - a)\big(\hat{\gamma}(X_i,M_i,a)-\gamma(X_i,M_i,a)\big)\lambda(a,X_i)R(M_i,X_i)\tag{R2}.
\end{align*}

For term \eqref{eq:term2}, note that
\begin{align*}
&K_h(A_i - a')\big\{\hat{\lambda}(a', X_i)\{\hat{\gamma}(X_i,M_i,a) - \hat{\eta}(a, a', X_i)\}
-\lambda(a', X_i)\{\gamma(X_i,M_i,a) - \eta(a, a', X_i)\}\big\}\\
&=K_h(A_i - a')\big(\hat{\lambda}(a', X_i)-\lambda(a', X_i)\big)\big(\hat{\gamma}(X_i,M_i,a)-\gamma(X_i,M_i,a)\big)\tag{CS5}\\
&\quad-K_h(A_i - a')\big(\hat{\lambda}(a', X_i)-\lambda(a', X_i)\big)\big(\hat{\eta}(a, a', X_i)-\hat{\eta}(a, a', X_i)\big)\tag{CS6}\\
&\quad+K_h(A_i - a')\big(\hat{\lambda}(a', X_i)-\lambda(a', X_i)\big)\gamma(X_i,M_i,a)\tag{R3}\\
&\quad+K_h(A_i - a')\big(\hat{\gamma}(X_i,M_i,a)-\gamma(X_i,M_i,a)\big)\lambda(a', X_i)\tag{R4}\\
&\quad-K_h(A_i - a')\big(\hat{\lambda}(a', X_i)-\lambda(a', X_i)\big)\hat{\eta}(a, a', X_i)\tag{R5}\\
&\quad-K_h(A_i - a')\big(\hat{\eta}(a, a', X_i)-\hat{\eta}(a, a', X_i)\big)\lambda(a', X_i)\tag{R6}.
\end{align*}

Next, we group terms (R1)-(R6) as follows. We pair (R1) with (R6), (R2) with (R4), and (R3) with (R5). Note that every expectation introduced here is only over $O_i$, conditional on $O^c_{I_\ell}$, i.e., $\mathbb{E}(\cdot| O^c_{I_\ell})$, and hence all the terms are random variables. For (R1)+(R6) we have

\begin{align*}
&(R1)+(R6)\\
&=\big(\hat{\eta}(a, a', X_i)-\eta(a, a', X_i)\big)
-K_h(A_i - a')\big(\hat{\eta}(a, a', X_i)-\eta(a, a', X_i)\big)\lambda(a', X_i)\\
&=\big(\hat{\eta}(a, a', X_i)-\eta(a, a', X_i)\big)
-\mathbb{E}\big[\hat{\eta}(a, a', X_i)-\eta(a, a', X_i) \big]\tag{E3}\\
&\quad-\Big\{
K_h(A_i - a')\big(\hat{\eta}(a, a', X_i)-\eta(a, a', X_i)\big)\lambda(a', X_i)\\
&-\mathbb{E}\big[K_h(A_i - a')\big(\hat{\eta}(a, a', X_i) -\eta(a, a', X_i)\big)\lambda(a', X_i) \mid  O^c_{I_\ell} \big]\Big\}\tag{E4}\\
&\quad+\mathbb{E}\big[
\big(\hat{\eta}(a, a', X_i)-\eta(a, a', X_i)\big)\big(1-K_h(A_i - a')\lambda(a', X_i)\big) \mid  O^c_{I_\ell}
\big]\tag{TR3}.
\end{align*}

For (R2)+(R4) we have
 \begin{align*}
 &(R2)+(R4)\\
 &=
 -K_h(A_i - a)\big(\hat{\gamma}(X_i,M_i,a)-\gamma(X_i,M_i,a)\big)\lambda(a,X_i)R(M_i,X_i)
 +\\
 & \hspace{10em}K_h(A_i - a')\big(\hat{\gamma}_a(M_i,X_i)-\gamma(X_i,M_i,a)\big)\lambda(a', X_i)\\
 &=-\Big\{  
 K_h(A_i - a)\big(\hat{\gamma}(X_i,M_i,a)-\gamma(X_i,M_i,a)\big)\lambda(a,X_i)R(M_i,X_i)\\
 &\quad\quad\quad-
 \mathbb{E}\big[
 K_h(A_i - a)\big(\hat{\gamma}(X_i,M_i,a)-\gamma(X_i,M_i,a)\big)\lambda(a,X_i)R(M_i,X_i)
 \mid  O^c_{I_\ell} \big]
 \Big\}\tag{E5}\\
 &\quad+\Big\{  
 K_h(A_i - a')\big(\hat{\gamma}(X_i,M_i,a)-\gamma(X_i,M_i,a)\big)\lambda(a', X_i)\\
 &\quad\quad\quad-
 \mathbb{E}\big[
 K_h(A_i - a')\big(\hat{\gamma}(X_i,M_i,a)-\gamma(X_i,M_i,a)\big)\lambda(a', X_i) \mid  O^c_{I_\ell} \big]
 \Big\}\tag{E6}\\
 &\quad+\mathbb{E}\big[
 \big(\hat{\gamma}_a(M_i,X_i)-\gamma(X_i,M_i,a)\big)
 \big\{
 K_h(A_i - a')\lambda_{a'}(X_i)\\
 &\qquad\qquad
 -K_h(A_i - a)\lambda(a,X_i)R(M_i,X_i)
 \big\} \big|  O^c_{I_\ell} \big]
 \tag{TR4}.
 \end{align*}

For (R3)+(R5) we have
\begin{align*}
&(R3)+(R5)\\
&=
K_h(A_i - a')\big(\hat{\lambda}(a', X_i)-\lambda(a', X_i)\big)\gamma(X_i,M_i,a)\\
&\qquad\qquad
-K_h(A_i - a')\big(\hat{\lambda}(a', X_i)-\lambda(a', X_i)\big)\eta(a, a', X_i)\\
&=\Big\{  
K_h(A_i - a')\big(\hat{\lambda}(a', X_i)-\lambda(a', X_i)\big)\gamma(X_i,M_i,a)\\
&\quad\quad\quad-
\mathbb{E}\big[
K_h(A_i - a')\big(\hat{\lambda}(a', X_i)-\lambda(a', X_i)\big)\gamma(X_i,M_i,a) \mid  O^c_{I_\ell} \big]
\Big\}\tag{E7}\\
&\quad-\Big\{  
K_h(A_i - a')\big(\hat{\lambda}(a', X_i)-\lambda(a', X_i)\big)\eta(a, a', X_i)\\
&\quad\quad\quad-
\mathbb{E}\big[
K_h(A_i - a')\big(\hat{\lambda}(a', X_i)-\lambda(a', X_i)\big)\eta(a, a', X_i) \mid  O^c_{I_\ell} \big]
\Big\}\tag{E8}\\
&\quad+\mathbb{E}\big[
K_h(A_i - a')\big(\hat{\lambda}(a', X_i)-\lambda(a', X_i)\big)
\big\{
\gamma(X_i,M_i,a)-\eta(a, a', X_i)
\big\} \mid  O^c_{I_\ell} \big]\tag{TR5}.
\end{align*}

And so, to prove $\sqrt{\frac{h^{d_A}}{n}} \sum_{i \epsilon I_\ell} \left\{m(O_i;\hat{\alpha}_\ell,\hat{\lambda}_\ell,\hat{\gamma}_\ell,\psi_0(a,a')) - m(O_i;\alpha,\lambda,\gamma,\psi_0(a,a')) \right\} = o_p(1)$, we provide proofs for the convergence of the terms (CS1) - (CS6), (E1) - (E8) and (TR1) - (TR5) in the following sub-sections.    

\vspace{1em}
\noindent\textbf{Proof for Terms (CS1)-(CS6)}
\vspace{1em}

All of these terms contain the product of two or more errors and can be treated similarly. We provide a detailed proof for (CS2), and a similar method can be followed for the rest of the terms. 

For (CS2), write $\Delta_{i\ell}=K_h(A_i - a)\big[\hat{R}(M_i,X_i)-R(M_i,X_i)\big]\big[\hat{\gamma}(X_i,M_i,a)-\gamma(X_i,M_i,a)\big]\big[Y_i - \gamma(X_i,M_i,a)\big]$. Following Lemma \ref{lem:bounding-term}, it suffices to bound $\mathbb{E}\left[|\sqrt{\frac{h^{d_A}}{n}} \sum_{i \in I_\ell} \Delta_{i\ell} | \Big| O^c_{I_\ell}\right]$ as $o_p(1)$ in order to show that $$\sqrt{\frac{h^{d_A}}{n}} \sum_{i \in I_\ell}\Delta_{i\ell} = o_p(1).$$

First, from the triangle inequality, $\mathbb{E}\left[\left|\sqrt{\frac{h^{d_A}}{n}} \sum_{i \in I_\ell} \Delta_{i\ell} \right| \mid O^c_{I_\ell}\right] \leq \frac{1}{L}\sqrt{nh^{d_A}} \mathbb{E}\left[\left|\Delta_{i\ell} \right| \mid O^c_{I_\ell}\right]$, and so it suffices to bound $\sqrt{nh^{d_A}}\mathbb{E} \bigg[\big|\Delta_{i\ell} \big|\bigg| O^c_{I_\ell}\bigg]$. In the interest of space, we introduce the following notation $\Tilde{k}(u) = \prod_{j = 1}^{d_A} k(u_j)$, where $u$ is a vector in $\mathbb{R}^{d_A}$.

{\small

\begin{align*}
    &\sqrt{nh^{d_A}}\mathbb{E} \bigg[\big|\Delta_{i\ell} \big|\bigg| O^c_{I_\ell}\bigg] \\
    =& \sqrt{nh^{d_A}}\int\bigg| K_h(A_i - a)\big[\hat{R}(M_i,X_i)-R(M_i,X_i)\big]\big[\hat{\gamma}(X_i,M_i,a)-\gamma(X_i,M_i,a)\big]\big[Y_i - \gamma(X_i,M_i,a)\big]\bigg| \\
    &\qquad\qquad\times f(Y_i, A_i, M_i, X_i)dO_i\\
    =&\sqrt{nh^{d_A}} \int\bigg|\Tilde{k}(u)\big[\hat{R}(M_i,X_i)-R(M_i,X_i)\big]\big[\hat{\gamma}(X_i,M_i,a)-\gamma(X_i,M_i,a)\big]\left[Y_i - \gamma(X_i,M_i,a)\right]\bigg|\\
    &\qquad\qquad\times f(Y_i, uh+a, M_i, X_i) dudY_idM_idX_i\\
    =&  \sqrt{nh^{d_A}}\int\left\{\int \bigg|\Tilde{k}(u)f(uh+a|M_i,X_i)\bigg\{\int \bigg|\left[Y_i - \gamma(X_i,M_i,a)\right]\bigg| f(Y_i | uh+a, M_i, X_i) dY_i\bigg\}\bigg|du\right\}\\
    &\hspace{5em}\bigg|\big[\hat{R}(M_i,X_i)-R(M_i,X_i)\big]\big[\hat{\gamma}(X_i,M_i,a)-\gamma(X_i,M_i,a)\big]f(M_i,X_i)\bigg|dM_idX_i\\
\end{align*}
}
Next, Assumption 3.1 on the boundedness of $\gamma(X,M,a)$ and Assumption 3.3 on the boundedness of $var(Y_i|a,m,x)$, along with an application of Lemma \ref{lem:kernelexpansion} on $f(a \mid M, X)$, we get 
\begin{align*}
=& O(\sqrt{nh^{d_A}} )\int
\left\{f(a \mid M_i, X_i) + O(h^2) \right\}\\
&\qquad \times \bigg|\big[\hat{R}(M_i,X_i)-R(M_i,X_i)\big]\big[\hat{\gamma}(X_i,M_i,a)-\gamma(X_i,M_i,a)\big]\bigg|f(M_i,X_i)dM_idX_i\\
=& O(\sqrt{nh^{d_A}} )\int
f(a \mid M_i, X_i) \bigg|\big[\hat{R}(M_i,X_i)-R(M_i,X_i)\big]\big[\hat{\gamma}(X_i,M_i,a)-\gamma(X_i,M_i,a)\big]\bigg|\\
&\qquad \times f(M_i,X_i)dM_idX_i\\
&+ O(\sqrt{nh^{d_A + 4}})\int
\bigg|\big[\hat{R}(M_i,X_i)-R(M_i,X_i)\big]\big[\hat{\gamma}(X_i,M_i,a)-\gamma(X_i,M_i,a)\big]\bigg|\\
&\qquad\qquad \times f(M_i,X_i)dM_idX_i\\
\overset{(a)}{\le} & O(\sqrt{nh^{d_A}} )\int
\bigg|\big[\hat{R}(M_i,X_i)-R(M_i,X_i)\big]\big[\hat{\gamma}(X_i,M_i,a)-\gamma(X_i,M_i,a)\big]\bigg|f(M_i,X_i)dM_idX_i\\
&+ O(\sqrt{nh^{d_A + 4}})\int
\bigg|\big[\hat{R}(M_i,X_i)-R(M_i,X_i)\big]\big[\hat{\gamma}(X_i,M_i,a)-\gamma(X_i,M_i,a)\big]\bigg|\\
&\qquad\qquad \times f(M_i,X_i)dM_idX_i\\
\overset{(b)}{\le} & O\bigg(\sqrt{nh^{d_A}} \bigg\{ \int\big[\hat{R}(M_i,X_i)-R(M_i,X_i)\big] ^2f(M_i,X_i)dM_idX_i \\
&\qquad\qquad\int\big[\hat{\gamma}(X_i,M_i,a)-\gamma(X_i,M_i,a)\big] ^2f(M_i,X_i)dM_idX_i\bigg\}^{1/2}\bigg) + o_p(1)\\
    &=  o_p(1).
\end{align*}

Where $(a)$ follows from an application of Holder's inequality combined with Assumption 3.1 on the boundedness of $f(a \mid M, X)$, and (b) and the last equality follows from an application of Cauchy-Schwartz, combined with Assumption 5.1 and $nh^{d_A + 4} \rightarrow C_h$ by Assumption 2.

\vspace{1em}
\noindent\textbf{Proof for Terms (E1)-(E8)}
\vspace{1em}

Terms (E1)-(E8) are normalized terms of the form of a bias times a bounded quantity; they can all be treated similarly. We only provide the proof of the convergence in probability to zero for the term (E2). (E2) is given as

\begin{align*}
&K_h(A_i - a)\big(\hat{\lambda}(a, X_i)-\lambda(a,X_i)\big)R(M_i,X_i)\big(Y_i - \gamma(X_i,M_i,a)\big)\\
-&\mathbb{E}\big[K_h(A_i - a)\big(\hat{\lambda}(a, X_i)-\lambda(a,X_i)\big)R(M_i,X_i)\big(Y_i - \gamma(X_i,M_i,a)\big)\mid  O^c_{I_\ell}\big]
\end{align*}

To prove $\sqrt{nh^{d_A}}$ times (E2) is $o_p(1)$, we set  $\hat{\Delta}_{i\ell}$ as (E2).
By construction, $O^c_{I_\ell}$ and $O_i$ are independent, $i\in I_\ell$, and consequently $\mathbb{E}\left[\hat{\Delta}_{i\ell}|O^c_{I_\ell}\right]=0$ and $\mathbb{E}\left[\hat{\Delta}_{i\ell}\hat{\Delta}_{j\ell}|O^c_{I_\ell}\right]=0$ for $i,j \in I_\ell$ and all $a', a\in \mathcal{A}_0$. Next we note that
\begin{align*}
&h^{d_A}\mathbb{E}\left[\hat{\Delta}^2_{i\ell}|O^c_{I_\ell}\right]\\
    = & h^{d_A}\int K^2_h(A_i-a)  \left[\hat{\lambda}(a, X_i) - \lambda(a,X_i) \right]^2R^2_i(M_i,X_i)\left[Y_i - \gamma(X_i,M_i,a)\right]^2 \\
    &\qquad\qquad \times f(Y_i, A_i, M_i, X_i) dO_i\\
    =& \int \Tilde{k}(u)^2\left[\hat{\lambda}(a, X_i) - \lambda(a,X_i) \right]^2R^2_i(M_i,X_i)\left[Y_i - \gamma(X_i,M_i,a)\right]^2\\
    &\qquad\qquad \times f(Y_i, uh+a, M_i, X_i) dudY_idM_idX_i\\
    =& \iint \Tilde{k}(u)^2f(uh+a|M_i,X_i)\bigg\{\int \left[Y_i - \gamma(X_i,M_i,a)\right]^2 f(Y_i | uh+a, M_i, X_i) dY_i\bigg\}du\\
    &\hspace{5em}\left[\hat{\lambda}(a, X_i) - \lambda(a,X_i) \right]^2R^2_i(M_i,X_i)f(M_i,X_i)dM_idX_i\\
    \overset{(a)}{=} & O\bigg(\int\Tilde{k}(u)^2du \int\left[\hat{\lambda}(a, X_i) - \lambda(a,X_i) \right]^2R^2_i(M_i,X_i)f(M_i,X_i)dM_idX_i\bigg)\\
    \overset{(b)}{=} & O(1)\int\left[\hat{\lambda}(a, X_i) - \lambda(a,X_i) \right]^2R^2_i(M_i,X_i)f(M_i,X_i)dM_idX_i\\
    \overset{(c)}{=} & o_p(1) 
\end{align*}
Where (a) follows from Assumption 3.1 on the boundedness of $f(a \mid M, X)$, along with Assumption 3.1 and Assumption 3.3 combined with the derivation provided below

\begin{align*}
    &\int\left[Y_i - \gamma(X_i,M_i,a)\right]^2 f(Y_i | uh+a, M_i, X_i) dY_i\\
    =& \int\left[Y^2_i + \gamma^2_a(M_i, X_i)-2\gamma(X_i,M_i,a)Y_i\right] f(Y_i | uh+a, M_i, X_i) dY_i\\
    =& \mathbb{E}[Y^2_i | uh+a, M_i, X_i] + \gamma^2_a(M_i, X_i) -2\gamma(X_i,M_i,a) \int_{\mathcal{Y}} Y_i f(Y_i | uh+a, M_i, X_i) dY_i\\
    =&\mathbb{E}[Y^2_i | uh+a, M_i, X_i] + \gamma^2_a(M_i, X_i) -2\gamma(X_i,M_i,a) - 2\gamma(X_i,M_i,a)\gamma_{uh+a}(M_i,X_i)\\
    =& O(1).
\end{align*}
Next, (b) follows from Assumption 2.4, and finally, (c) follows Assumption 3.2 along with Assumption 4.1. Then 
$\mathbb{E} \bigg[\left(\sqrt{h^{d_A}/n} \sum^L_{l=1}\sum_{i\in I_\ell} \hat{\Delta}_{i\ell}\right)^2\bigg| O^c_{I_\ell}\bigg] =h^{d_A} /n \sum^L_{\ell=1}\sum_{i\in I_\ell} \mathbb{E}\left[\hat{\Delta}^2_{i\ell}|O^c_{I_\ell}\right] = h^{d_A}\mathbb{E}\left[\hat{\Delta}^2_{i\ell}|O^c_{I_\ell}\right] = o_p(1).$

Applying Lemma 1 to the above gives $\sqrt{h^{d_A}/n} \sum^L_{l=1}\sum_{i\in I_\ell} \hat{\Delta}_{i\ell}\xrightarrow[]{P} 0$, i.e. $\sqrt{nh^{d_A}}$ times (E2) being $o_p(1)$. 

\vspace{1em}
\noindent\textbf{Proof for Terms (TR1)-(TR5)}
\vspace{1em}

The proofs of the convergence in probability to zero for the terms (TR1)-(TR5) require extra considerations, and we prove them on a case by case basis below. 

Terms (TR1) and (TR2) are similar; we only provide the proof of the convergence in probability to zero for the term (TR2).

To bound TR2, first set $$\hat{\Delta}_{i\ell}=K_h(A_i - a) \left[\hat{\lambda}(a, X_i) - \lambda(a,X_i) \right]R(M_i, X_i)\left[Y_i - \gamma(X_i,M_i,a)\right].$$ Bounding (TR2) amounts to showing $\sqrt{nh^{d_A}}\mathbb{E}[\hat{\Delta}_{i\ell}|O^c_{I_\ell}] = o_p(1)$.

\begin{align*}
&\sqrt{nh^{d_A}}\mathbb{E} \bigg[\hat{\Delta}_{i\ell} \bigg| O^c_{I_\ell}\bigg] \\
    =& \sqrt{nh^{d_A}}\mathbb{E}\bigg\{K_h(A_i - a) 
    \left[\hat{\lambda}(a, X_i) - \lambda(a,X_i) \right]
    R(M_i, X_i)
    \left[Y_i - \gamma(X_i,M_i,a)\right]\bigg| O^c_{I_\ell}\bigg\}\\
    =& \sqrt{nh^{d_A}}\int K_h(A_i-a)  \left[\hat{\lambda}(a, X_i) - \lambda(a,X_i) \right]R(M_i, X_i)\left[Y_i - \gamma(X_i,M_i,a)\right]\\
    &\qquad\qquad \times f(Y_i, A_i, M_i, X_i) dO_i\\
    =&\sqrt{nh^{d_A}}\int \left[ \int K_h(A_i-a)f(A_i \mid Y_i, M_i, X_i) dA_i\right]  \left[\hat{\lambda}(a, X_i) - \lambda(a,X_i) \right]\\
    &\hspace{8em}R(M_i, X_i)\left[Y_i - \gamma(X_i,M_i,a)\right]f(Y_i,M_i, X_i) dY_idM_idX_i\\
    \intertext{Applying Lemma \ref{lem:kernelexpansion} under Assumption 3.1}
    =&\sqrt{nh^{d_A}}\int \left[ f(a \mid Y_i, M_i, X_i) +O(h^2) \right] \\
    &\hspace{4em}\left[\hat{\lambda}(a, X_i) - \lambda(a,X_i) \right]R(M_i, X_i)\left[Y_i - \gamma(X_i,M_i,a)\right]f(Y_i,M_i, X_i) dY_idM_idX_i\\
    \overset{(a)}{=}&\sqrt{nh^{d_A}}\int O(h^2) \left[\hat{\lambda}(a, X_i) - \lambda(a,X_i) \right]\\
    &\hspace{8em}R(M_i, X_i)\left[Y_i - \gamma(X_i,M_i,a)\right]f(Y_i,M_i, X_i) dY_idM_idX_i\\
     \overset{(b)}{=} &O(\sqrt{nh^{d_A + 4}})\int \Big|\left[\hat{\lambda}(a, X_i) - \lambda(a,X_i) \right]R(M_i, X_i)\Big|\\
    &\hspace{8em}\left[\int | Y_i-\gamma(X_i,M_i,a) |  f(Y_i\mid M_i, X_i) dY_i \right]  f(M_i,X_i)dM_idX_i\\
     \overset{(c)}{=}& o_p(1)
\end{align*}

where (a) follows from 
\begin{align*}
    &\int [Y_i -\gamma(X_i,M_i,a)]  f(Y_i\mid a, M_i, X_i) dY_i \\
    =& \int Y_i  f(Y_i\mid a, M_i, X_i) dY_i-\gamma(X_i,M_i,a) = 0,
\end{align*}
(b) is from the exchange of $O(\cdot)$ and integration,
(c) follows from Assumption 2 ($nh^{d_A + 4} \rightarrow C_h$, $h \rightarrow 0$), Assumption 3 and Assumption 4.1, Cauchy-Schwartz combined with the boundedness of $\int |Y_i -\gamma(X_i,M_i,a)| f(Y_i\mid M_i, X_i) dY_i $  derived from Assumption 3.1 shown below
\begin{align*}
    &\int |Y_i -\gamma(X_i,M_i,a)|  f(Y_i\mid M_i, X_i) dY_i \\
    =& \int \bigg[\int |Y_i -\gamma(X_i,M_i,a)| f(Y_i\mid a, M_i, X_i) dY_i \bigg]f(a|M_i,X_i)da \\
    \le & \int \bigg[\text{Var}(Y_i|a,M_i, X_i)\bigg]^{1/2} f(a|M_i,X_i)da < \infty,
\end{align*}
where the last line also comes from the Cauchy-Schwartz inequality.

For Term (TR3), we have
\begin{align*}
&\sqrt{nh^{d_A}}\mathbb{E}\left[
\big(\hat{\eta}(a, a', X_i)-\eta(a, a', X_i)\big)\big(1-K_h(A_i - a')\lambda(a', X_i)\big)
\big| O^c_{I_\ell}\right] \\   
&=\sqrt{nh^{d_A}}\int
\big(\hat{\eta}(a, a', X_i)-\eta(a, a', X_i)\big)\big(1-K_h(A_i - a')\lambda(a', X_i)\big)f(A_i,X_i)dA_idX_i \\
&=\sqrt{nh^{d_A}}\int 
\big(\hat{\eta}(a, a', X_i)-\eta(a, a', X_i)\big)\big(1-\bigg\{\int K_h(A_i - a')f(A_i\mid X_i)dA_i\bigg\}\lambda(a', X_i)\big)\\
&\qquad\qquad \times f(X_i)dX_i \\
&\overset{(a)}{=}\sqrt{nh^{d_A}}\int 
\big(\hat{\eta}(a, a', X_i)-\eta(a, a', X_i)\big)\big(1- f(a'\mid X_i)\lambda_{a'}(X_i)\big)\\
&\qquad\qquad \times f(X_i)dX_i +\sqrt{nh^{d_A}}\int 
\big(\hat{\eta}(a, a', X_i)-\eta(a, a', X_i)\big)O(h^2)\lambda_{a'}(X_i)f(X_i)dX_i\\
&\overset{(b)}{=} o_p(1).
\end{align*}
where $(a)$ follows from Lemma \ref{lem:kernelexpansion},  and (b) follows from the definition of $\lambda_{a^\prime}(X_i)$, $nh^{d_A+4}\rightarrow C_h$, Assumption 4 (convergence of $\hat{\eta}(X_i)$), Assumption 3 (boundedness of $\lambda$) combined with an application of Cauchy-Schwartz inequality.

Demonstrating the bound for (TR4), we have
\begin{align*}
    &\sqrt{nh^{d_A}}\mathbb{E}\big[
\big(\hat{\gamma}(X_i,M_i,a)-\gamma(X_i,M_i,a)\big)\\
&\qquad\qquad \times 
\big\{
K_h(A_i - a')\lambda(a', X_i)
-K_h(A_i - a)\lambda(a,X_i)R(M_i,X_i)
\big\}
\big]\\
&= \sqrt{nh^{d_A}}\mathbb{E}\big[
\big(\hat{\gamma}(X_i,M_i,a)-\gamma(X_i,M_i,a)\big)
\big\{
K_h(A_i - a')\lambda(a', X_i)
\big\}
\big]\tag{TR4-1}\\
 - &\sqrt{nh^{d_A}}\mathbb{E}\big[
\big(\hat{\gamma}(X_i,M_i,a)-\gamma(X_i,M_i,a)\big)
\big\{
K_h(A_i - a)\lambda(a,X_i)R(M_i,X_i)
\big\}
\big]\tag{TR-4-2}
\end{align*}

TR-4-1 can be written as
{\small
\begin{align*}
    &\sqrt{nh^{d_A}}\mathbb{E}\big[\big(\hat{\gamma}(X_i,M_i,a)-\gamma(X_i,M_i,a)\big)\big\{K_h(A_i - a')\lambda(a', X_i)\big\}\big]\\
    &= \sqrt{nh^{d_A}} \int \big(\hat{\gamma}(X_i,M_i,a)-\gamma(X_i,M_i,a)\big)\lambda(a', X_i)\left\{ \int K_h(A_i - a') f(A_i \mid M_i, X_i) dA_i \right\} \\
    &\qquad\qquad \times f(M_i, X_i) dM_i dX_i
\end{align*}
}

An application of Lemma \ref{lem:kernelexpansion} to TR-4-1 gives
\begin{align*}
&\sqrt{nh^{d_A}}\mathbb{E}\big[
\big(\hat{\gamma}(X_i,M_i,a)-\gamma(X_i,M_i,a)\big)
\big\{
K_h(A_i - a')\lambda(a', X_i)
\big\}
\big]\\
    &=\sqrt{nh^{d_A}} \int \big(\hat{\gamma}(X_i,M_i,a)-\gamma(X_i,M_i,a)\big)\lambda(a', X_i)f(a' \mid M_i, X_i)f(M_i, X_i) dM_i dX_i\\
    & + \sqrt{nh^{d_A}}\int \big(\hat{\gamma}(X_i,M_i,a)-\gamma(X_i,M_i,a)\big)\lambda(a', X_i)O(h^2)f(M_i, X_i) dM_i dX_i
\end{align*}

A similar approach applied to TR-4-2 gives
\begin{align*}
&\sqrt{nh^{d_A}}\mathbb{E}\big[
\big(\hat{\gamma}(X_i,M_i,a)-\gamma(X_i,M_i,a)\big)
\big\{
K_h(A_i - a)\lambda(a,X_i)R(M_i,X_i)
\big\}
\big]\\
    =&\sqrt{nh^{d_A}} \int \big(\hat{\gamma}(X_i,M_i,a)-\gamma(X_i,M_i,a)\big)\lambda(a,X_i)R(M_i,X_i)f(a \mid M, X) f(M_i, X_i) dM_i dX_i\\
    &+ \sqrt{nh^{d_A}}\int \big(\hat{\gamma}(X_i,M_i,a)-\gamma(X_i,M_i,a)\big)\lambda(a,X_i)R(M_i,X_i)O(h^2)f(M_i, X_i) dM_i dX_i
\end{align*}

Now, the first terms of TR-4-1 and TR-4-2 cancel out with each other, shown below
{\small
\begin{align*}
    &\int \big(\hat{\gamma}(X_i,M_i,a)-\gamma(X_i,M_i,a)\big)\left\{\lambda(a', X_i)f(a' \mid M, X)
    -  \lambda(a,X_i)R(M_i,X_i)f(a \mid M, X)\right\} \\
    & \qquad \times f(M_i, X_i) dM_i dX_i = 0
\end{align*}
}
This can be seen from
\begin{align*}
    \lambda(a', X_i)f(a' \mid M_i, X_i) = \frac{f(X_i)}{f(a', X_i)}\frac{f(a', M_i, X_i)}{f(M_i, X_i)}
\end{align*}

Along with
\begin{align*}
    \lambda(a,X_i)R(M_i,X_i)f(a \mid M_i, X_i) &= \frac{f(X_i)}{f(a, X_i)}\frac{f(M_i, a', X_i)}{f(a', X_i)}
    \frac{f(a, X_i)}{f(M_i, a, X_i)}\frac{f(a, M_i, X_i)}{f(M_i, X_i)}\\
    &= \frac{f(X_i)}{f(a', X_i)}\frac{f(M_i , a', X_i)}{f(M_i, X_i)}
\end{align*}

Consequently the first terms in TR4-1 and TR4-2 cancel each other out, and this leaves us to bound the remaining terms. 
\[
\sqrt{nh^{d_A}} \int \big(\hat{\gamma}(X_i,M_i,a)-\gamma(X_i,M_i,a)\big)\lambda(a', X_i) O(h^2)f(M, X) dM_i dX_i = o_p(1)
\]
The second term in TR-4-1 and TR-4-2 can be bounded by an application of Cauchy-Schwartz, combined with Assumption 4 (consistency of $\hat{\gamma}$) and boundedness of $\lambda$ in Assumption 3.1.

Finally, for term (TR5), we note that
\begin{align*}
&\sqrt{nh^{d_A}}\mathbb{E}\left[
K_h(A_i - a')\big(\hat{\lambda}(a', X_i)-\lambda(a', X_i)\big)
\big\{
\gamma(X_i,M_i,a)-\eta(X_i)
\big\}\big| O^c_{I_\ell}
\right]\\
&=\sqrt{nh^{d_A}}\int
K_h(A_i - a')\big(\hat{\lambda}(a', X_i)-\lambda(a', X_i)\big)
\big\{
\gamma(X_i,M_i,a)-\eta(a, a', X_i)
\big\}\\
&\qquad\qquad \times f(A_i,M_i,X_i)dA_idM_idX_i\\
&=\sqrt{nh^{d_A}}\int
\Bigg\{\int K_h(A_i - a')f(A_i\mid M_i,X_i)dA_i\Bigg\}
\big(\hat{\lambda}(a', X_i)-\lambda(a', X_i)\big)\\
&\hspace{2in}\times \big\{
\gamma(X_i,M_i,a)-\eta(a, a', X_i)
\big\}f(M_i,X_i)dM_idX_i\\
&\overset{(a)}{=}\sqrt{nh^{d_A}}\int \big(\hat{\lambda}(a', X_i)-\lambda(a', X_i)\big)
\big\{
\gamma(X_i,M_i,a)-\eta(a, a', X_i)
\big\}f(a',M_i,X_i)dM_idX_i\\
&\quad+\sqrt{nh^{d_A}}\int \big(\hat{\lambda}(a', X_i)-\lambda(a', X_i)\big)
\big\{
\gamma(X_i,M_i,a)-\eta(a, a', X_i)
\big\}O(h^2)f(M_i,X_i)dM_idX_i\\
&\overset{(b)}{=} O(\sqrt{nh^{d_A + 4}})\int \Bigg|\big(\hat{\lambda}(a', X_i)-\lambda(a', X_i)\big)
\big\{
\gamma(X_i,M_i,a)-\eta(a, a', X_i)
\big\}\Bigg|f(M_i,X_i)dM_idX_i\\
&\overset{(c)}{=} o_p(1)
\end{align*}

Where $(a)$ follows from an application of Lemma \ref{lem:kernelexpansion}, (b) follows from the definition of $\eta$, and (c) follows from an application of Cauchy-Schwartz combined with the consistency of $\hat{\lambda}$.

\vspace{1em}
\noindent\textbf{Proof of Asymptotic Normality}
\vspace{1em}

The proof for asymptotic normality follows from an application of the Lyapunov Central Limit theorem to the terms $\sqrt{nh^{d_A}}n^{-1} m(O_i;\alpha,\lambda,\gamma,\psi_0(a,a'))$. We first prove the Lyapunov condition holds for $\delta = 1$, i.e. 
\[
\lim_{n \rightarrow \infty} \frac{1}{s^{3}_n} \sum_{i = 1}^n \mathbb{E}\left[\Big|\sqrt{nh^{d_A}}n^{-1} m(O_i;\alpha,\lambda,\gamma,\psi_0(a,a')) - \mu_{i} \Big|^3\right] = 0
\]
Where $\mu_i$ equals $\mathbb{E}\left[\sqrt{nh^{d_A}}n^{-1} m(O_i;\alpha,\lambda,\gamma,\psi_0(a,a'))\right]$ and $s^2_n = \sum_{i = 1}^n \sigma^2_i$ where $\sigma^2_i$ is the variance of of $\sqrt{nh^{d_A}}n^{-1} m(O_i;\alpha,\lambda,\gamma,\psi_0(a,a'))$. To prove the Lyaponuv condition holds, we first derive $\mu_i$ and $\sigma^2_i$.

\vspace{1em}
\noindent\textbf{Calculation for $B(a, a')$ and $\mu_i$}
\vspace{1em}

Given
\begin{align*}
    &m(O_i;\alpha,\lambda,\gamma,\psi_0(a,a')) \\
    =& \frac{K_h(A_i - a)f(M_i \mid A = a', X_i)}{f(M_i \mid A = a, X_i)f(a \mid X_i)}\{Y_i - \mathbb{E}[Y \mid X_i, M_i , A = a]\}\\
    &+ \frac{K_h(A_i - a')}{f(a' \mid X_i)}\{\mathbb{E}[Y \mid X_i, M_i, A = a] - \eta(a, a', X_i)\} + \eta(a, a', X_i) - \psi_0(a,a')
\end{align*}
Since $\mathbb{E}[\eta(a, a', X_i) - \psi_0(a,a')] = 0$, we focus on
$\frac{K_h(A - a)f(M \mid A = a', X)}{f(M \mid A = a, X)f(a \mid X)}\{Y - \mathbb{E}[Y \mid X, M , A = a]\}+ \frac{K_h(A - a')}{f(a' \mid X)}\{\mathbb{E}[Y \mid X, M, A = a] - \eta(a, a', X)\}$.

We start by computing the expectation of each the individual terms one at a time.

\noindent \textbf{Expectation Part 1}
$$\frac{K_h(A - a)f(M \mid A = a', X)}{f(M \mid A = a, X)f(a \mid X)}\{Y - \mathbb{E}[Y \mid X, M , A = a]\}$$

From $\mathbb{E}\{\mathbb{E}[\gamma(X,M,a)|X,M]\} = \mathbb{E}\{\mathbb{E}[Y|X,M]\}$ and $\gamma(X,M,a) = \mathbb{E}(Y|X,M,A)$, expectation of the first term
\begin{align*}
&\mathbb{E}\bigg[\frac{K_h(A - a)f(M \mid A = a', X)}{f(M \mid A = a, X)f(a \mid X)}\{Y - \mathbb{E}[Y \mid X, M , A = a]\}\bigg]   \\
=& \mathbb{E}\Bigg\{\mathbb{E}\bigg[\frac{K_h(A - a)f(M \mid A = a', X)}{f(M \mid A = a, X)f(a \mid X)}\{Y - \mathbb{E}[Y \mid X, M , A = a]\}\bigg| X,M\bigg]\Bigg\}\\
=& \mathbb{E}\Bigg\{\frac{f(M \mid A = a', X)}{f(M \mid A = a, X)f(a \mid X)} \mathbb{E}\bigg[K_h(A - a) (\gamma(X,M,a) - \gamma(X,M,a))\bigg| X,M\bigg]\Bigg\}.
\end{align*}

The inner product further expands as follows,
\begin{align*}
    &\mathbb{E}\bigg[K_h(A - a) (\gamma(X,M,a) - \gamma(X,M,a))\bigg| X,M\bigg]\\
    =& \int K_h(A - a) (\gamma(X,M,a) - \gamma(X,M,a)) f(A|X,M) dA\\
    =& \int \bigg[\prod^{d_A}_{j=1}\frac{1}{h}k\Big(\frac{A_j-a}{h}\Big)\bigg] (\gamma(X,M,a) - \gamma(X,M,a)) f(A|X,M) dA\\
    =& \int \bigg[\prod^{d_A}_{j=1} k(u_j)\bigg] (\gamma(a+uh, X, M) - \gamma(X,M,a)) f(a+uh|X,M) du\\
    =& \int \bigg[\prod^{d_A}_{j=1} k(u_j)\bigg] \bigg(\sum^{d_A}_{j=1} u_j h\partial_{a_j} \gamma(X,M,a) +\sum^{d_A}_{j=1}\sum^{d_A}_{j^\prime=1}\frac{u_ju_{j^\prime}h^2}{2} \partial_{a_j}\partial_{a_{j^\prime}} \gamma(X,M,a) \bigg)\\
    & \times \bigg( f(a|X,M)+\sum^{d_A}_{j=1} u_j h\partial_{a_j} f(a|X,M) +\frac{u^2_j h^2}{2} \partial^2_{a_j} f(a|X,M)\bigg)du_1\cdots du_{d_A} + O(h^3)\\
    =& h^2 \int u^2 k(u)du \bigg(\sum^{d_A}_{j=1} \partial_{a_j} \gamma(X,M,a) \partial_{a_j} f(a|X,M) + \frac{1}{2} \left[\sum^{d_A}_{j=1} \partial^2_{a_j} \gamma(X,M,a)\right]f(a|X,M) \bigg)\\
    &\qquad\qquad + O(h^3)
\end{align*}
for all $X,M$ in respective range. Inserting this back into the original expectation we get,
{\small
\begin{align*}
    &\mathbb{E}\Bigg\{\frac{f(M \mid A = a', X)}{f(M \mid A = a, X)f(a \mid X)} \mathbb{E}\bigg[K_h(A - a) (\gamma(X,M,a) - \gamma(X,M,a))\bigg| X,M\bigg]\Bigg\} \\
    =& h^2 \int u^2 k(u)du \\
    &\times \mathbb{E}\bigg[ \frac{f(M \mid A = a', X)}{f(M \mid A = a, X)f(a \mid X)}\bigg(\sum^{d_A}_{j=1} \partial_{a_j} \gamma(X,M,a) \partial_{a_j} f(a|X,M) \\
    &+ \frac{1}{2} \left[\sum^{d_A}_{j=1} \partial^2_{a_j} \gamma(X,M,a)\right]f(a|X,M) \bigg)\bigg]+ O(h^3).
\end{align*}
}

\noindent \textbf{Expectation Part 2}
$$\frac{K_h(A - a')}{f(a' \mid X)}\{\gamma(X,M,a) - \eta(a, a', X)\}$$
\begin{align*}
    &\mathbb{E}\bigg[\frac{K_h(A - a')}{f(a' \mid X)}\{\gamma(X,M,a) - \eta(a, a', X)\}\bigg]\\
    =&\mathbb{E}\bigg\{\mathbb{E}\bigg[\frac{K_h(A - a')}{f(a' \mid X)}\{\gamma(X,M,a) - \eta(a, a', X)\}\bigg| X,M\bigg]\bigg\} \\
    =&\mathbb{E}\bigg\{\frac{1}{f(a' \mid X)}\mathbb{E}\bigg[K_h(A - a')\{\gamma(X,M,a) - \eta(a, a', X)\}\bigg| X,M\bigg]\bigg\}\\
    =&\mathbb{E}\bigg\{\frac{\gamma(X,M,a) - \eta(a, a', X)}{f(a' \mid X)}\mathbb{E}[K_h(A - a')| X,M]\bigg\}
\end{align*}

The inner expectation can be written as 
\begin{align*}
    &\mathbb{E}\bigg[K_h(A - a')\bigg| X,M\bigg]\\
    =& \int \bigg[\prod^{d_A}_{j=1}\frac{1}{h}k\Big(\frac{A_j-a'}{h}\Big)\bigg]  f(A|X,M) dA \\
    =& \int k(u_1)\cdots k(u_{d_A}) \bigg( f(a'|X,M)+ \sum^{d_A}_{j=1} u_jh\partial_{a_j} f(a'|X,M) \\
    &+\sum^{d_A}_{j=1} \sum^{d_A}_{j^\prime=1}\frac{u_ju_{j^\prime}h^2}{2} \partial_{a_j}\partial_{a_{j^\prime}} f(a^\prime|X,M) \\
    &+ \sum^{d_A}_{j=1} \sum^{d_A}_{j^\prime=1}\sum^{d_A}_{j^{\prime\prime}=1}\frac{u_ju_{j^\prime}u_{j^{\prime\prime}}h^3}{2} \partial_{a_j}\partial_{a_{j^\prime}}\partial_{a_{j^{\prime\prime}}} f(\Bar{a}|X,M)\bigg)du_1\cdots du_{d_A}\\
    =&f(a'|X,M) + \frac{1}{2}h^2\int u^2k(u)du \sum^{d_A}_{j=1}\partial^2_{a_j} f(a'|X,M) + O(h^3)
\end{align*}

Plugging this back into the above expectation
{\small
\begin{align*}
    &\mathbb{E}\Bigg\{ \frac{1}{f(a'|X)}\{\gamma(X,M,a) - \eta(a, a', X)\}\times \bigg(f(a'|X,M) + \frac{1}{2}h^2\int u^2k(u)du \sum_{j = 1}^{h_{d_A}}\partial^2_{a_j} f(a'|X,M) \bigg)\Bigg\}\\
    &\qquad\qquad + O(h^3)\\
    =& \mathbb{E}\Bigg\{ \{\gamma(X,M,a) - \eta(a, a', X)\} \bigg(\frac{f(a'|X,M)}{f(a'|X)} + \frac{1}{2}h^2\int u^2k(u)du \frac{\sum_{j = 1}^{h_{d_A}} \partial^2_{a_j} f(a'|X,M)}{f(a'|X)} \bigg)\Bigg\} + O(h^3)\\
    =& h^2\left[\int u^2k(u)du\right]\mathbb{E}\left[ \{\gamma(X,M,a) - \eta(a, a', X)\}\frac{1}{2}\frac{\sum_{j = 1}^{h_{d_A}} \partial^2_{a_j} f(a'|X,M)}{f(a'|X)}\right]+ O(h^3)\\
\end{align*}
}
from having the first term in this expectation equal to zero, which we prove below
\begin{align*}
    & \mathbb{E}\left[ \{\gamma(X,M,a) - \eta(a, a', X)\}\frac{f(a'|X,M)}{f(a'|X)} \right] \\
    =& \int \{\gamma(X,M,a) - \eta(a, a', X)\}\frac{f(a'|X,M)}{f(a'|X)} f(M,X)dMdX \\
    =& \int \{\gamma(X,M,a) - \eta(a, a', X)\}\frac{f(A=a',X,M)}{f(A = a', X)} f(X)dMdX \\
    =& \int \{\gamma(X,M,a) - \eta(a, a', X)\}f(M|A=a',X) dM f(X)dX \\
    =& \int \{\eta(a, a', X)-\eta(a, a', X)\} f(X)dX=0
\end{align*}

Hence, letting 
\begin{align*}
    &B(a,a') = \\
    & \left[\int u^2 k(u)du\right]\times \mathbb{E}\bigg[ \frac{f(M \mid A = a', X)}{f(M \mid A = a, X)f(a \mid X)}\bigg(\sum^{d_A}_{j=1} \partial_{a_j} \gamma(X,M,a) \partial_{a_j} f(a|X,M) \\
    &+ \frac{1}{2} \left[\sum^{d_A}_{j=1} \partial^2_{a_j} \gamma(X,M,a)\right]f(a|X,M) \bigg)\\
    &+ \{\gamma(X,M,a) - \eta(a, a', X)\}\frac{1}{2}\frac{\sum_{j = 1}^{h_{d_A}} \partial^2_{a_j} f(a'|X,M)}{f(a'|X)}\bigg]+O(h),\\
\end{align*}

we have $\mathbb{E}\left[m(O_i;\alpha,\lambda,\gamma, \psi_0(a,a'))\right] = h^2 B(a,a')$. Additionally from this derivation $\mathbb{E}\left[\sqrt{nh^{d_A}}n^{-1} m(O_i;\alpha,\lambda,\gamma,\psi_0(a,a'))\right] = O(\sqrt{\frac{h^{d_A + 4}}{n}})$. Next, we prove the properties of variance.

\vspace{1em}
\noindent\textbf{Calculation for $V(a, a')$ and $s^2_n$}
\vspace{1em}

From the definition of $s^2_n$, we have
\begin{align*}
    &s^2_n = \sum_{i = 1}^n \sigma^2_i = \sum_{i = 1}^n var\left(\sqrt{nh^{d_A}}n^{-1} m(O_i;\alpha,\lambda,\gamma,\psi_0(a,a'))\right)\\
    =& h^{d_A} var\left(m(O_i;\alpha,\lambda,\gamma,\psi_0(a,a')\right)
\end{align*}

Consequently, we calculate
\begin{align*}
    &h^{d_A}\times var\Bigg\{\frac{K_h(A - a)f(M \mid A = a', X)}{f(M \mid A = a, X)f(a \mid X)}\{Y - \mathbb{E}[Y \mid X, M , A = a]\} \\
    & \qquad + \frac{K_h(A - a')}{f(a' \mid X)}\{\mathbb{E}[Y \mid X, M, A = a] - \eta(a, a', X)\}\\
    & \qquad +  \eta(a, a', X) - \psi_0(a, a^\prime)\Bigg\} 
\end{align*}

Using the property that $var(X) = \mathbb{E}[X^2] - \mathbb{E}[X]^2$ and constant values do not contribute to the variance, the variance term above can be re-written as
\begin{align*}
    &h^{d_A}\mathbb{E}\Bigg\{\bigg[\frac{K_h(A - a)f(M \mid A = a', X)}{f(M \mid A = a, X)f(a \mid X)}\{Y - \mathbb{E}[Y \mid X, M , A = a]\} + \\
    &\qquad \frac{K_h(A - a')}{f(a' \mid X)}\{\mathbb{E}[Y \mid X, M, A = a] - \eta(a, a', X)\} + \eta(a, a', X)\bigg]^2\Bigg\} \\
    &- h^{d_A}\mathbb{E}\Bigg\{\bigg[\frac{K_h(A - a)f(M \mid A = a', X)}{f(M \mid A = a, X)f(a \mid X)}\{Y - \mathbb{E}[Y \mid X, M , A = a]\} + \\
    &\qquad \frac{K_h(A - a')}{f(a' \mid X)}\{\mathbb{E}[Y \mid X, M, A = a] - \eta(a, a', X)\} + \eta(a, a', X)\bigg]\Bigg\}^2\\
\end{align*}

Examining each of the terms above one by one, the first term can be expanded as
\begin{align*}
    =& h^{d_A}\mathbb{E}\Bigg\{\bigg[\frac{K_h(A - a)f(M \mid A = a', X)}{f(M \mid A = a, X)f(a \mid X)}\{Y - \mathbb{E}[Y \mid X, M , A = a]\}  \bigg]^2\Bigg\}\\
    &+ h^{d_A}\mathbb{E}\Bigg\{\bigg[\frac{K_h(A - a')}{f(a' \mid X)}\{\mathbb{E}[Y \mid X, M, A = a] - \eta(a, a', X)\} \bigg]^2\Bigg\} +  h^{d_A}\mathbb{E}\Bigg\{\eta^2(a, a', X)\Bigg\}\\
    &+ 2h^{d_A}\mathbb{E}\Bigg\{\bigg[\frac{K_h(A - a)f(M \mid A = a', X)}{f(M \mid A = a, X)f(a \mid X)}\{Y - \mathbb{E}[Y \mid X, M , A = a]\}  \bigg]\\
    &\qquad 
    \bigg[\frac{K_h(A - a')}{f(a' \mid X)}\{\mathbb{E}[Y \mid X, M, A = a] - \eta(a, a', X)\} \bigg]\Bigg\}\\
    &+ 2h^{d_A}\mathbb{E}\Bigg\{\eta(a, a', X)\bigg[\frac{K_h(A - a)f(M \mid A = a', X)}{f(M \mid A = a, X)f(a \mid X)}\{Y - \mathbb{E}[Y \mid X, M , A = a]\}  \bigg]\Bigg\}\\
    &+ 2h^{d_A}\mathbb{E}\Bigg\{\eta(a, a', X)
    \bigg[\frac{K_h(A - a')}{f(a' \mid X)}\{\mathbb{E}[Y \mid X, M, A = a] - \eta(a, a', X)\} \bigg]\Bigg\}
\end{align*}

We analyze each of these terms part by part.

\vspace{2em}

\noindent \textbf{Variance Part 1}
\begin{align*}
    &h^{d_A}\mathbb{E}\Bigg\{\bigg[\frac{K_h(A - a)f(M \mid A = a', X)}{f(M \mid A = a, X)f(a \mid X)}\{Y - \mathbb{E}[Y \mid X, M , A = a]\}  \bigg]^2\Bigg\}\\
    =&h^{d_A}\mathbb{E}\Bigg\{\mathbb{E}\bigg\{\bigg[\frac{K_h(A - a)f(M \mid A = a', X)}{f(M \mid A = a, X)f(a \mid X)}\{Y - \mathbb{E}[Y \mid X, M , A = a]\}  \bigg]^2\bigg|X,M\bigg\}\Bigg\}\\
    =&h^{d_A}\mathbb{E}\Bigg\{\frac{f(M \mid A = a', X)^2}{f(M \mid A = a, X)^2f(a \mid X)^2}\mathbb{E}\bigg\{K_h(A - a)^2(Y - \mathbb{E}[Y \mid X, M , A = a])^2\bigg|X,M\bigg\}\Bigg\}\\
    =&h^{d_A}\mathbb{E}\Bigg\{\frac{f(M \mid A = a', X)^2}{f(M \mid A = a, X)^2f(a \mid X)^2}\\
    &\qquad \times \mathbb{E}\bigg\{K_h(A - a)^2\mathbb{E}\{(Y - \mathbb{E}[Y \mid X, M , A = a])^2|X,M,A\}\bigg|X,M\bigg\}\Bigg\}\\
    =&h^{d_A}\mathbb{E}\Bigg\{\frac{f(M \mid A = a', X)^2}{f(M \mid A = a, X)^2f(a \mid X)^2}\\
    &\times \mathbb{E}\bigg\{K_h(A - a)^2 \bigg[var(Y|X,M,A) + \gamma(X,M,a)^2 \\
&\qquad -2\gamma(X,M,a)\gamma(X,M,a)+\gamma(X,M,a)^2\bigg]
    \bigg|X,M\bigg\}\Bigg\}\\
    =&h^{d_A}\mathbb{E}\Bigg\{\frac{f(M \mid A = a', X)^2}{f(M \mid A = a, X)^2f(a \mid X)^2}\\
    &\times \mathbb{E}\bigg\{K_h(A - a)^2 \bigg[var(Y|X,M,A) + (\gamma(X,M,a)-\gamma(X,M,a))^2\bigg]
    \bigg|X,M\bigg\}\Bigg\}
\end{align*}
Because $0 < \int u^6 k(u) du < \infty$ from Assumption 2 (3), we also have boundedness of $\int u^6k^2(u)du$. The inner expectation can be written as 
\begin{align*}
     &h^{d_A}\mathbb{E}\bigg\{K_h(A - a)^2 \bigg[var(Y|X,M,A) + (\gamma(X,M,a)-\gamma(X,M,a))^2\bigg]
     \bigg|X,M\bigg\}\\
     = &h^{d_A} \int \bigg[\prod^{d_A}_{j=1}\frac{1}{h^2}k\Big(\frac{A_j-a_j}{h}\Big)^2\bigg]  \bigg\{var(Y|X,M,A) + [\gamma(X,M,a)-\gamma(X,M,a)]^2\bigg\}
     \\
     &\qquad\qquad \times f(A|X,M) dA\\
     = &\int \tilde{k}(u)^2\times \bigg\{var(Y|X,M,a+uh) + [\gamma(a+uh, M, X)-\gamma(X,M,a)]^2\bigg\}
     f(a+uh|X,M) du\\
     = & \int k(u_1)^2\cdots k(u_{d_A})^2 \\
     &\times \bigg\{var(Y|X,M,a)+\sum^{d_A}_{j=1}u_jh\partial_{a_j} var(Y|X,M,a) + \sum^{d_A}_{j=1}\sum^{d_A}_{j^\prime=1}u_ju_{j^\prime}h^2\partial_{a_j}\partial_{a_j^\prime} var(Y|X,M,\Bar{a}_v) \\
     &  \hspace{3em} + \Big[\sum^{d_A}_{j=1} u_jh\partial_{a_j} \gamma(X,M,a) + \sum^{d_A}_{j=1}\sum^{d_A}_{j'=1} u_ju_{j'}h^2\partial_{a_j}\partial_{a_{j'}} \gamma(\Bar{a}_{\gamma},M,X)\Big]^2\bigg\}\\
     &\times \bigg[f(a|X,M)+\sum^{d_A}_{j=1}u_jh\partial_{a_j} f(a|X,M) + \sum^{d_A}_{j=1}\sum^{d_A}_{j^\prime=1}u_ju_{j^\prime}h^2\partial_{a_j}\partial_{a_{j^\prime}} f(\Bar{a}_f|X,M) \bigg]du_1\cdots du_{d_A}\\
     =& \Big[\int \tilde{k}(u)^2 du\Big]\times  var(Y|X,M,a)f(a|X,M) + O(h^2)
\end{align*}
where $\Bar{a}_v,\Bar{a}_\gamma$, and $\Bar{a}_f$ are between $a$ and $a+uh$. Hence, part 1 of the variance 
\begin{align*}
    &h^{d_A}\mathbb{E}\Bigg\{\bigg[\frac{K_h(A - a)f(M \mid a', X)}{f(M \mid a, X)f(a \mid X)}\{Y - \mathbb{E}[Y \mid X, M , a]\}  \bigg]^2\Bigg\}\\
    =&\left[\int k(u)^2 du\right]^{d_A} \mathbb{E}\Bigg\{\frac{f(M \mid a', X)^2}{f(M \mid a, X)^2f(a \mid X)^2}var(Y|X,M,a)f(a|X,M)\Bigg\}+ O(h^2)
\end{align*}
\vspace{5em}

\noindent \textbf{Variance Part 2}
\begin{align*}
    &h^{d_A}\mathbb{E}\Bigg\{\bigg[\frac{K_h(A - a')}{f(a' \mid X)}\Big(\mathbb{E}(Y \mid X, M, A = a) - \eta(a, a', X)\Big) \bigg]^2\Bigg\}\\
    =&h^{d_A}\mathbb{E}\Bigg\{\frac{1}{f(a' \mid X)^2}\mathbb{E}\bigg[K_h(A - a')^2\Big(\gamma(X,M,a) - \eta(a, a', X)\Big)^2\Big|X\bigg]\Bigg\}\\
\end{align*}

The inner expectation can be written as
\begin{align*}
    &h^{d_A}\mathbb{E}\bigg[K_h(A - a')^2\Big(\gamma(X,M,a) - \eta(a, a', X)\Big)^2\Big|X\bigg]\\
    =& h^{d_A}\int \bigg[\prod^{d_A}_{j=1}\frac{1}{h^2}k\Big(\frac{A_j-a'}{h}\Big)^2\bigg]  \Big(\gamma(X,M,a) - \eta(a, a', X)\Big)^2 f(A| M,X)f(M\mid X) dAdM \\
    =&\int \bigg[\prod^{d_A}_{j=1} k(u_j)^2 \bigg]\Big(\gamma(X,M,a) - \eta(a, a', X)\Big)^2 f(a'+uh|M, X)f(M\mid X) dudM\\
    =&\int \bigg[\prod^{d_A}_{j=1} k(u_j)^2 \bigg]\Big(\gamma(X,M,a) - \eta(a, a', X)\Big)^2 \Big\{f(a'|M, X) + \sum_{j = 1}^{d_A}u_jh\partial_{a_j}f(a \mid X, M)\\
    &+ \sum_{j = 1}^{d_A}\sum_{j^\prime = 1}^{d_A}u_ju_{j^\prime}h^2\partial_{a_j}\partial_{a_{j^\prime}}f(\Bar{a} \mid X, M)\Big\}f(M\mid X) dudM\\
    =& \int k^2(u_1)\cdots k^2(u_{d_A})\Big(\gamma(X,M,a) - \eta(a, a', X)\Big)^2f(a'|X, M)f(M\mid X)
    du_1\cdots du_{d_A}dM \\
    &\qquad\qquad + O(h^2)\\
    =& \left[\int k(u)^2 du \right]^{d_A} \times var[E(Y|X,M,a)|X,a']f(a'|X) + O(h^2)\\
\end{align*}

the last equation is from
\begin{align*}
    var[E(Y|X,M,a)|X,a']=& \mathbb{E}\bigg\{\Big[E(Y|X,M,a)-\eta(a,a',X)\Big]^2|X,a'\bigg\}\\
    =& \int \Big[E(Y|X,M,a)-\eta(a,a',X)\Big]^2 f(M|X,a')dM.
\end{align*}

Hence, the part 2 of variance
\begin{align*}
&h^{d_A}\mathbb{E}\Bigg\{\frac{1}{f(a' \mid X)^2}\mathbb{E}\bigg[K_h(A - a')^2\Big(\gamma(X,M,a) - \eta(a, a', X)\Big)^2\Big|X\bigg]\Bigg\}\\
=& \left[\int k(u)^2 du \right]^{d_A} \times \mathbb{E}\bigg\{\frac{1}{f(a'|X)}var[E(Y|X,M,a)|X,a']\bigg\} + O(h^2)
\end{align*}

\noindent \textbf{Variance Part 3}
\begin{align*}
    h^{d_A}\mathbb{E}\left[\eta^2(a, a', X)\right] &= O(h^{d_A})
\end{align*}

This holds because we assume $\eta$ is bounded.

\noindent \textbf{Variance Part 4}
\begin{align*}
    &h^{d_A}\mathbb{E}\Bigg\{\bigg[\frac{K_h(A - a)f(M \mid A = a', X)}{f(M \mid A = a, X)f(a \mid X)}\{Y - \mathbb{E}[Y \mid X, M , A = a]\}  \bigg]\\
    &\qquad \times \bigg[\frac{K_h(A - a')}{f(a' \mid X)}\{\mathbb{E}[Y \mid X, M, A = a] - \eta(a, a', X)\} \bigg]\Bigg\}\\
    =&h^{d_A}\mathbb{E}\Bigg\{\frac{K_h(A-a)K_h(A-a')}{f(a|X)f(a'|X)}\frac{f(M|a',X)}{f(M|a,X)}\Big[Y-\gamma(X,M,a)\Big]\Big[\gamma(X,M,a)-\eta(a,a',X)\Big]\Bigg\}\\
    =&h^{d_A}\mathbb{E}\Bigg\{\frac{1}{f(a|X)f(a'|X)}\frac{f(M|a',X)}{f(M|a,X)}
    \Big[\gamma(X,M,a)-\eta(a,a',X)\Big]\\
    & \qquad\times\mathbb{E}\Big\{K_h(A-a)K_h(A-a')\Big[Y-\gamma(X,M,a)\Big]\Big|X,M\Big\}\Bigg\}\\
    =&h^{d_A}\mathbb{E}\Bigg\{\frac{1}{f(a|X)f(a'|X)}\frac{f(M|a',X)}{f(M|a,X)}
    \Big[\gamma(X,M,a)-\eta(a,a',X)\Big]\\
    & \qquad\times\mathbb{E}\Big\{K_h(A-a)K_h(A-a')\Big[\gamma(X,M,a)-\gamma(X,M,a)\Big]\Big|X,M\Big\}\Bigg\}
\end{align*}

The inner expectation
\begin{align*}
    &h^{d_A}\mathbb{E}\Big\{K_h(A-a)K_h(A-a')\Big[\gamma(X,M,a)-\gamma(X,M,a)\Big]\Big|X,M\Big\}\\
    =& h^{d_A}\int \bigg[\prod^{d_A}_{j=1}\frac{1}{h^2}k\Big(\frac{A_j-a}{h}\Big)k\Big(\frac{A_j-a'}{h}\Big)\bigg]
  \Big[\gamma(X,M,a)-\gamma(X,M,a)\Big] f(A|X,M)dA\\
    =& \int k(u_1)\cdots k(u_{d_A}) k(u_1+\frac{a-a'}{h})\cdots k(u_{d_A}+\frac{a-a'}{h})\\
    &\times \Big[\gamma(uh + a, M, X)-\gamma(X,M,a)\Big] f(uh+ a|X,M)dA\\
    =& \int k(u_1)\cdots k(u_{d_A}) k(u_1+\frac{a-a'}{h})\cdots k(u_{d_A}+\frac{a-a'}{h})\\  &\times \Big[\sum^{d_A}_{j=1} u_j h\partial_{a_j} \gamma(X,M,a)+\frac{u^2_j h^2}{2}\partial^2_{a_j}\gamma(X,M,a)+\frac{u^3_jh^3}{6}\partial^3_{a_j}\gamma(\Bar{a}, M, X)\Big]\\
    &\times \Big[f(a|X,M) +\sum^{d_A}_{j=1} u_jh\partial_{a_j} f(a|X,M)+\frac{u^2_jh^2}{2}\partial^2_{a_j} f(\Bar{a}|X,M)\Big]du_1\cdots du_{d_A}\\
    &= O(h)
\end{align*}

Hence, the part 4 of variance 
\begin{align*}
    &\mathbb{E}\Bigg\{\frac{1}{f(a|X)f(a'|X)}\frac{f(M|a',X)}{f(M|a,X)}
    \Big[\gamma(X,M,a)-\eta(a,a',X)\Big]\\
    & \qquad\times\mathbb{E}\Big\{K_h(A-a)K_h(A-a')\Big[\gamma(X,M,a)-\gamma(X,M,a)\Big]\Big|X,M\Big\}\Bigg\}\\
    &= O(h)
\end{align*}

\noindent \textbf{Variance Part 5}
\begin{align*}
        &2h^{d_A}\mathbb{E}\Bigg\{\eta(a, a', X)\bigg[\frac{K_h(A - a)f(M \mid A = a', X)}{f(M \mid A = a, X)f(a \mid X)}\{Y - \mathbb{E}[Y \mid X, M , A = a]\}  \bigg]\Bigg\}\\
        &= 2h^{d_A}\mathbb{E}\Bigg\{\frac{\eta(a, a', X)f(M \mid A = a', X)}{f(M \mid A = a, X)f(a \mid X)}\mathbb{E}\bigg[\{Y - \mathbb{E}[Y \mid X, M , A = a] \mid X, M\}  \bigg]\Bigg\}\\
\end{align*}
Applying the same expansion in Expectation Part 1, we can write the inner expectation as 
\begin{align*}
   & h^2 \int u^2 k(u)du \bigg(\sum^{d_A}_{j=1} \partial_{a_j} \gamma(X,M,a) \partial_{a_j} f(a|X,M) + \frac{1}{2} \left[\sum^{d_A}_{j=1} \partial^2_{a_j} \gamma(X,M,a)\right]f(a|X,M) \bigg)\\
   &+ O(h^3)
\end{align*}
Inserting this back into the full expectation, combined with the boundedness of $\eta$, $f(M \mid a, X)$ and $f(a \mid X)$, we get
\begin{align*}
    2h^{d_A}\mathbb{E}\Bigg\{\eta(a, a', X)\bigg[\frac{K_h(A - a)f(M \mid A = a', X)}{f(M \mid A = a, X)f(a \mid X)}\{Y - \mathbb{E}[Y \mid X, M , A = a]\}  \bigg]\Bigg\} = O(h^{d_A + 2}) 
\end{align*}

\noindent \textbf{Variance Part 6}
\begin{align*}
    &2h^{d_A}\mathbb{E}\Bigg\{\eta(a, a', X)
    \bigg[\frac{K_h(A - a')}{f(a' \mid X)}\{\mathbb{E}[Y \mid X, M, A = a] - \eta(a, a', X)\} \bigg]\Bigg\}\\
    =2h^{d_A}&\mathbb{E}\bigg\{\mathbb{E}\bigg[\frac{K_h(A - a')\eta(a, a^\prime, X)}{f(a' \mid X)}\{\gamma(X,M,a) - \eta(a, a', X)\}\bigg| X,M\bigg]\bigg\} \\
    =2h^{d_A}&\mathbb{E}\bigg\{\frac{\{\gamma(X,M,a) - \eta(a, a', X)\}\eta(a, a^\prime, X)}{f(a' \mid X)}\mathbb{E}[K_h(A - a')| X,M]\bigg\}
\end{align*}

Using the expansion from Part 2 of the expectation on $\mathbb{E}[K_h(A - a')| X,M]$, we get
\begin{align*}
       \mathbb{E}[K_h(A - a')| X,M] =&f(a'|X,M) + \frac{1}{2}h^2\int u^2k(u)du \sum^{d_A}_{j=1}\partial^2_{a_j} f(a'|X,M) + O(h^3)
\end{align*}
Plugging this back into the full expectation, we get
\begin{align*}
    2h^{d_A}\mathbb{E}\Bigg\{\eta(a, a', X)
    \bigg[\frac{K_h(A - a')}{f(a' \mid X)}\{\mathbb{E}[Y \mid X, M, A = a] - \eta(a, a', X)\} \bigg]\Bigg\} = O(h^{d_A})
\end{align*}

And using the calculation for the bias,
\begin{align*}
    &h^{d_A}\mathbb{E}\Bigg\{\bigg[\frac{K_h(A - a)f(M \mid A = a', X)}{f(M \mid A = a, X)f(a \mid X)}\{Y - \mathbb{E}[Y \mid X, M , A = a]\} + \\
    &\qquad \frac{K_h(A - a')}{f(a' \mid X)}\{\mathbb{E}[Y \mid X, M, A = a] - \eta(a, a', X)\} + \eta(a, a', X)\bigg]\Bigg\}^2\\ 
    &= O(h^{d_A + 4})\\
\end{align*}
Finally, putting the pieces of the variance together, we have 
\begin{align*}
    &h^{d_A}\times var\Bigg\{\frac{K_h(A - a)f(M \mid A = a', X)}{f(M \mid A = a, X)f(a \mid X)}\{Y - \mathbb{E}[Y \mid X, M , A = a]\} \\
    & \qquad + \frac{K_h(A - a')}{f(a' \mid X)}\{\mathbb{E}[Y \mid X, M, A = a] - \eta(a, a', X)\}\\
    & \qquad +  \eta(a, a', X) - \psi_0(a, a^\prime)\Bigg\} = V(a,a') + O(h)
\end{align*}
where the term converges to $V(a, a^\prime)$ as $h\rightarrow 0$ and 
\begin{align*}
    V(a, a^\prime) = &\left[\int k(u)^2  du\right]^{d_A}\times\mathbb{E}\Bigg\{\frac{f(M \mid a', X)^2}{f(M \mid a, X)^2f(a \mid X)^2}var(Y|X,M,a)f(a|X,M) \\
    &+ \frac{1}{f(a'|X)}var[E(Y|X,M,a)|X,a']\Bigg\}.\\
\end{align*}

Having derived the bias and variance terms, we now prove the Lyapunov condition for $\delta = 1$.

\vspace{1em}
\noindent\textbf{Proof for Lyapunov Condition}
\vspace{1em}

We now prove the Lyapunov condition
\[
\lim_{n \rightarrow \infty} \frac{1}{s^{3}_n} \sum_{i = 1}^n \mathbb{E}\left[\Big|\sqrt{nh^{d_A}}n^{-1} m(O_i;\alpha,\lambda,\gamma,\psi_0(a,a')) - \mu_{i} \Big|^3\right] = 0
\]

Note that
$$\Big|\sqrt{nh^{d_A}}n^{-1} m(O_i;\alpha,\lambda,\gamma,\psi_0(a,a')) - \mu_{i} \Big| \leq \Big|\sqrt{nh^{d_A}}n^{-1} m(O_i;\alpha,\lambda,\gamma,\psi_0(a,a'))\Big| +  \Big| \mu_{i} \Big|$$
Since both sides are positive,
\begin{align*}
    \Big|\sqrt{nh^{d_A}}n^{-1} m(O_i;\alpha,\lambda,\gamma,\psi_0(a,a')) - \mu_{i} \Big|^3 \leq &(h^{d_A}n^{-1})^{3/2}\Big| m(O_i;\alpha,\lambda,\gamma,\psi_0(a,a'))\Big|^3 +  \Big| \mu_{i} \Big|^3\\ + & 3(h^{d_A}n^{-1})\Big| m(O_i;\alpha,\lambda,\gamma,\psi_0(a,a'))\Big|^2\Big|\mu_i\Big|\\
    + & 3(h^{d_A}n^{-1})^{1/2}\Big| m(O_i;\alpha,\lambda,\gamma,\psi_0(a,a'))\Big|\Big|\mu_i\Big|^2
\end{align*}
From the monotonicity of the expected value, we have
\begin{align*}
    &\sum_{i = 1}^n \mathbb{E}\left[\Big|\sqrt{nh^{d_A}}n^{-1} m(O_i;\alpha,\lambda,\gamma,\psi_0(a,a')) - \mu_{i} \Big|^3\right] \leq\\
    &\qquad \sum_{i = 1}^n \mathbb{E}\left[(h^{d_A}n^{-1})^{3/2}\Big| m(O_i;\alpha,\lambda,\gamma,\psi_0(a,a'))\Big|^3\right] +  \sum_{i = 1}^n\Big| \mu_{i} \Big|^3\\ 
    &\qquad +  \sum_{i = 1}^n 3h^{d_A}n^{-1}\mathbb{E}\left[\Big| m(O_i;\alpha,\lambda,\gamma,\psi_0(a,a'))\Big|^2\right]\Big|\mu_i\Big|\\
    &\qquad +  \sum_{i = 1}^n 3(h^{d_A}n^{-1})^{1/2}\mathbb{E}\left[\Big| m(O_i;\alpha,\lambda,\gamma,\psi_0(a,a'))\Big|\right]\Big|\mu_i\Big|^2
\end{align*}

Since $$\sum_{i = 1}^n\Big| \mu_{i} \Big|^3 = O(h^{(d_A + 4)3/2}n^{-1/2}) = o(1),$$  $$\sum_{i = 1}^n 3h^{d_A}n^{-1}\mathbb{E}\left[\Big| m(O_i;\alpha,\lambda,\gamma,\psi_0(a,a'))\Big|^2\right]\Big|\mu_i\Big| = O(\sqrt{\frac{h^{d_A + 4}}{n}}) = o(1),\text{ and}$$  $$\sum_{i = 1}^n 3(h^{d_A}n^{-1})^{1/2}\mathbb{E}\left[\Big| m(O_i;\alpha,\lambda,\gamma,\psi_0(a,a'))\Big|\right]\Big|\mu_i\Big|^2 = o(1),$$it suffices to prove the following condition
\[
\lim_{n \rightarrow \infty} \frac{1}{s^{3/2}_n} \sum_{i = 1}^n \mathbb{E}\left[\Big|\sqrt{nh^{d_A}}n^{-1} m(O_i;\alpha,\lambda,\gamma,\psi_0(a,a')) \Big|^3\right] = 0
\]

Following a similar derivation as in the proof for consistency of $\hat{V}(a, a')$, from the assumption $\mathbb{E}\{[Y-\gamma(X,M,a)]^3 | A = a', M = m, X = x\}$ over any $(a,a',m,x)\in \mathcal{A}\times\mathcal{A}\times \mathcal{M} \times \mathcal{X}$, along with $\int^{\infty}_{-\infty} k(u)^{c_1}k(u+\tilde{c})^{c_2} du < \infty$ and $\int^{\infty}_{-\infty} u^2k(u)^{c_1}k(u+\tilde{c})^{c_2} du < \infty$ for $\tilde{c} \in \mathcal{R}$ and $c_1+c_2 \in \{2,3\}$ for $c_1,c_2\in\{0,1,2,3\}$, we can bound $\mathbb{E}\left[\Big|m(O_i;\alpha,\lambda,\gamma,\psi_0(a,a')) \Big|^3\right] = O(\frac{1}{h^{2d_A}})$. Hence
\[
\sum_{i = 1}^n \mathbb{E}\left[\Big|\sqrt{nh^{d_A}}n^{-1} m(O_i;\alpha,\lambda,\gamma,\psi_0(a,a')) \Big|^3\right] = O((nh^{d_A})^{-1/2}) = o(1)
\]

Combining this with $s^2_n = V(a, a') + o(1)$ proves the Lyapunov condition. Hence,
$$
\frac{1}{s_n} \sum_{i = 1}^n \left(\sqrt{\frac{h^{d_A}}{n}}m(O_i;\alpha,\lambda,\gamma,\psi_0(a,a')) - \mu_i\right) \xrightarrow[]{d} \mathcal{N}(0, 1)
$$

An application of Slutsky's theorem provides the desired result that
$$\sqrt{nh^{d_A}} (\hat{\psi}^{MR}(a,a')- \psi_0(a,a') -B(a,a')) \xrightarrow[]{d} N(0, V(a,a'))$$ 

\vspace{1em}
\subsection{Proof of Proposition 1}

Following a similar breakdown as that in Theorem 1, $\hat{\psi}^{MR}(a, a') - \psi(a, a')$ can be expanded as
\begin{align*}
    & \hat{\psi}^{MR}(a,a')- \psi_0(a,a')\\
    =& \frac{1}{\sqrt{nh^{d_A}}} \times \sqrt{\frac{h^{d_A}}{n}}\sum^L_{\ell=1}\sum_{i\in I_\ell} \Big\{ m(O_i;\alpha,\lambda,\gamma,\psi_0(a,a'))\Big\}\\
    &+\frac{1}{\sqrt{nh^{d_A}}} \times \sqrt{\frac{h^{d_A}}{n}}\sum^L_{\ell=1}\sum_{i\in I_\ell}\Big\{m(O_i;\hat{\alpha}_\ell,\hat{\lambda}_\ell,\hat{\gamma}_\ell,\psi_0(a,a'))- m(O_i;\alpha,\lambda,\gamma,\psi_0(a,a'))\Big\}.
\end{align*}

Following the result in Theorem 1 on asymptotic normality by application of the Lyapunov CLT, the term $\sqrt{\frac{h^{d_A}}{n}}\sum^L_{\ell=1}\sum_{i\in I_\ell} \Big\{ m(O_i;\alpha,\lambda,\gamma,\psi_0(a,a'))\Big\} = O_p(1)$. Since $\frac{1}{\sqrt{nh^{d_A}}} = o_p(1)$, the following holds 

\[\frac{1}{\sqrt{nh^{d_A}}} \times \sqrt{\frac{h^{d_A}}{n}}\sum^L_{\ell=1}\sum_{i\in I_\ell} \Big\{ m(O_i;\alpha,\lambda,\gamma,\psi_0(a,a'))\Big\} = o_p(1).
\]

The remainder of the proof demonstrates the remaining the remaining terms are $o_p(1)$, i.e.
\[
\frac{1}{\sqrt{nh^{d_A}}} \times \sqrt{\frac{h^{d_A}}{n}}\sum^L_{\ell=1}\sum_{i\in I_\ell}\Big\{m(O_i;\hat{\alpha}_\ell,\hat{\lambda}_\ell,\hat{\gamma}_\ell,\psi_0(a,a'))- m(O_i;\alpha,\lambda,\gamma,\psi_0(a,a'))\Big\}.
\]

To see this, we first expand these terms identically as the proof of Theorem 1, and provide proofs for the convergence of the terms (CS1) - (CS6), (E1) - (E8) and (TR1) - (TR5) under the assumption that any two out of three nuisance models are correctly specified in the following sub-sections.    

\vspace{1em}
\noindent\textbf{Proof for Terms (CS1)-(CS6)}
\vspace{1em}

All of these terms contain the product of two or more errors and can be treated similarly. We provide a detailed proof for (CS2), and a similar method can be followed for the rest of the terms. 

For (CS2), write $\Delta_{i\ell}=K_h(A_i - a)\big[\hat{R}(M_i,X_i)-R(M_i,X_i)\big]\big[\hat{\gamma}(X_i,M_i,a)-\gamma(X_i,M_i,a)\big]\big[Y_i - \gamma(X_i,M_i,a)\big]$. Following Lemma \ref{lem:bounding-term}, it suffices to bound $\mathbb{E}\left[|\frac{1}{\sqrt{nh^{d_A}}} \times \sqrt{\frac{h^{d_A}}{n}} \sum_{i \in I_\ell} \Delta_{i\ell} | \Big| O^c_{I_\ell}\right]$ as $o_p(1)$ in order to show that $$\frac{1}{\sqrt{nh^{d_A}}} \times \sqrt{\frac{h^{d_A}}{n}} \sum_{i \in I_\ell}\Delta_{i\ell} = o_p(1).$$

First, from the triangle inequality, $\mathbb{E}\left[\left|\frac{1}{\sqrt{nh^{d_A}}} \times \sqrt{\frac{h^{d_A}}{n}} \sum_{i \in I_\ell} \Delta_{i\ell} \right| \mid O^c_{I_\ell}\right] \leq \frac{1}{L}\mathbb{E}\left[\left|\Delta_{i\ell} \right| \mid O^c_{I_\ell}\right]$, and so it suffices to bound $\mathbb{E} \bigg[\big|\Delta_{i\ell} \big|\bigg| O^c_{I_\ell}\bigg]$. 

{\small

\begin{align*}
    &\mathbb{E} \bigg[\big|\Delta_{i\ell} \big|\bigg| O^c_{I_\ell}\bigg] \\
    =&  \int\bigg| K_h(A_i - a)\big[\hat{R}(M_i,X_i)-R(M_i,X_i)\big]\big[\hat{\gamma}(X_i,M_i,a)-\gamma(X_i,M_i,a)\big]\big[Y_i - \gamma(X_i,M_i,a)\big]\bigg| \\
    &\qquad\qquad \times f(Y_i, A_i, M_i, X_i)dO_i\\
    =& \int\bigg|\Tilde{k}(u)\big[\hat{R}(M_i,X_i)-R(M_i,X_i)\big]\big[\hat{\gamma}(X_i,M_i,a)-\gamma(X_i,M_i,a)\big]\left[Y_i - \gamma(X_i,M_i,a)\right]\bigg|\\
    &\qquad\qquad \times f(Y_i, uh+a, M_i, X_i) dudY_idM_idX_i\\
    =&\int\left\{\int \bigg|\Tilde{k}(u)f(uh+a|M_i,X_i)\bigg\{\int \bigg|\left[Y_i - \gamma(X_i,M_i,a)\right]\bigg| f(Y_i | uh+a, M_i, X_i) dY_i\bigg\}\bigg|du\right\}\\
    &\hspace{5em}\bigg|\big[\hat{R}(M_i,X_i)-R(M_i,X_i)\big]\big[\hat{\gamma}(X_i,M_i,a)-\gamma(X_i,M_i,a)\big]f(M_i,X_i)\bigg|dM_idX_i\\
\end{align*}
}
Next, Assumption 3.1 on the boundedness of $\gamma(X,M,a)$ and Assumption 3.3 on the boundedness of $var(Y_i|a,m,x)$, along with an application of Lemma \ref{lem:kernelexpansion} on $f(a \mid M, X)$, we get 
{\small
\begin{align*}
=& O(1)\int
\left\{f(a \mid M_i, X_i) + O(h^2) \right\}\bigg|\big[\hat{R}(M_i,X_i)-R(M_i,X_i)\big]\big[\hat{\gamma}(X_i,M_i,a)-\gamma(X_i,M_i,a)\big]\bigg|\\
&\qquad \times f(M_i,X_i)dM_idX_i\\
=& O(1)\int
f(a \mid M_i, X_i) \bigg|\big[\hat{R}(M_i,X_i)-R(M_i,X_i)\big]\big[\hat{\gamma}(X_i,M_i,a)-\gamma(X_i,M_i,a)\big]\bigg|f(M_i,X_i)dM_idX_i\\
&+ O(h^2)\int
\bigg|\big[\hat{R}(M_i,X_i)-R(M_i,X_i)\big]\big[\hat{\gamma}(X_i,M_i,a)-\gamma(X_i,M_i,a)\big]\bigg|f(M_i,X_i)dM_idX_i\\
\overset{(a)}{\le} & O(1)\int
\bigg|\big[\hat{R}(M_i,X_i)-R(M_i,X_i)\big]\big[\hat{\gamma}(X_i,M_i,a)-\gamma(X_i,M_i,a)\big]\bigg|f(M_i,X_i)dM_idX_i\\
&+ O(h^2)\int
\bigg|\big[\hat{R}(M_i,X_i)-R(M_i,X_i)\big]\big[\hat{\gamma}(X_i,M_i,a)-\gamma(X_i,M_i,a)\big]\bigg|f(M_i,X_i)dM_idX_i
\end{align*}
}
As long as either Assumption 4.2 or Assumption 4.3 hold, then combined with Assumption 3.2, 
(CS2) will be $o_p(1)$. A similar approach can be used to bound the remaining CS terms.

\vspace{1em}
\noindent\textbf{Proof for Terms (E1)-(E8)}
\vspace{1em}

Terms (E1)-(E8) are normalized terms of the form of a bias times a bounded quantity; they can all be treated similarly. We only provide the proof of the convergence in probability to zero for the term (E2). (E2) is given as

\begin{align*}
&K_h(A_i - a)\big(\hat{\lambda}(a, X_i)-\lambda(a,X_i)\big)R(M_i,X_i)\big(Y_i - \gamma(X_i,M_i,a)\big)\\
-&\mathbb{E}\big[K_h(A_i - a)\big(\hat{\lambda}(a, X_i)-\lambda(a,X_i)\big)R(M_i,X_i)\big(Y_i - \gamma(X_i,M_i,a)\big)\mid  O^c_{I_\ell}\big]
\end{align*}

To prove this, we set  $\hat{\Delta}_{i\ell}$ as (E2). By construction, $O^c_{I_\ell}$ and $O_i$ are independent, $i\in I_\ell$, and consequently $\mathbb{E}\left[\hat{\Delta}_{i\ell}|O^c_{I_\ell}\right]=0$ and $\mathbb{E}\left[\hat{\Delta}_{i\ell}\hat{\Delta}_{j\ell}|O^c_{I_\ell}\right]=0$ for $i,j \in I_\ell$ and all $a', a\in \mathcal{A}_0$. Next we note that
\begin{align*}
&\mathbb{E}\left[\hat{\Delta}^2_{i\ell}|O^c_{I_\ell}\right]\\
    = & \int K^2_h(A_i-a)  \left[\hat{\lambda}(a, X_i) - \lambda(a,X_i) \right]^2R^2_i(M_i,X_i)\left[Y_i - \gamma(X_i,M_i,a)\right]^2 \\
    &\qquad\qquad \times f(Y_i, A_i, M_i, X_i) dO_i\\
    =& \frac{1}{h^{d_A}}\int \Tilde{k}(u)^2\left[\hat{\lambda}(a, X_i) - \lambda(a,X_i) \right]^2R^2_i(M_i,X_i)\left[Y_i - \gamma(X_i,M_i,a)\right]^2\\
    &\qquad \times f(Y_i, uh+a, M_i, X_i) dudY_idM_idX_i\\
    =& \frac{1}{h^{d_A}}\int \int \Tilde{k}(u)^2f(uh+a|M_i,X_i)\bigg\{\int \left[Y_i - \gamma(X_i,M_i,a)\right]^2 f(Y_i | uh+a, M_i, X_i) dY_i\bigg\}du\\
    &\hspace{5em}\left[\hat{\lambda}(a, X_i) - \lambda(a,X_i) \right]^2R^2_i(M_i,X_i)f(M_i,X_i)dM_idX_i\\
    \overset{(a)}{=} & O\bigg(\frac{1}{h^{d_A}}\int\Tilde{k}(u)^2du \int\left[\hat{\lambda}(a, X_i) - \lambda(a,X_i) \right]^2R^2_i(M_i,X_i)f(M_i,X_i)dM_idX_i\bigg)\\
    \overset{(b)}{=} & O(\frac{1}{h^{d_A}})\int\left[\hat{\lambda}(a, X_i) - \lambda(a,X_i) \right]^2R^2_i(M_i,X_i)f(M_i,X_i)dM_idX_i\\
    \overset{(c)}{=} & O_p(\frac{1}{h^{d_A}}) 
\end{align*}
Where (a) follows from Assumption 3.1 on the boundedness of $f(a \mid M, X)$, along with Assumption 3.1 and Assumption 3.3 combined with the derivation provided below
\begin{align*}
    &\int\left[Y_i - \gamma(X_i,M_i,a)\right]^2 f(Y_i | uh+a, M_i, X_i) dY_i\\
    =& \int\left[Y^2_i + \gamma^2_a(M_i, X_i)-2\gamma(X_i,M_i,a)Y_i\right] f(Y_i | uh+a, M_i, X_i) dY_i\\
    =& \mathbb{E}[Y^2_i | uh+a, M_i, X_i] + \gamma^2_a(M_i, X_i) -2\gamma(X_i,M_i,a) \int Y_i f(Y_i | uh+a, M_i, X_i) dY_i\\
    =&\mathbb{E}[Y^2_i | uh+a, M_i, X_i] + \gamma^2_a(M_i, X_i) -2\gamma(X_i,M_i,a) - 2\gamma(X_i,M_i,a)\gamma_{uh+a}(M_i,X_i)\\
    =& O(1).
\end{align*}
Next, (b) follows from Assumption 2.4, and finally, (c) follows Assumption 3.2 along with Assumption 4.1. Then 
$\mathbb{E} \bigg[\left(\frac{1}{\sqrt{nh^{d_A}}} \times \sqrt{h^{d_A}/n} \sum_{i\in I_\ell} \hat{\Delta}_{i\ell}\right)^2\bigg| O^c_{I_\ell}\bigg] =\frac{1}{n^2} \sum_{i\in I_\ell} \mathbb{E}\left[\hat{\Delta}^2_{i\ell}|O^c_{I_\ell}\right] = O(\frac{1}{n})\mathbb{E}\left[\hat{\Delta}^2_{i\ell}|O^c_{I_\ell}\right] = O_p(\frac{1}{nh^{d_A}}) = o_p(1)$.

Applying Lemma 1 to the above gives $\frac{1}{\sqrt{nh^{d_A}}} \times \sqrt{h^{d_A}/n} \sum^L_{l=1}\sum_{i\in I_\ell} \hat{\Delta}_{i\ell}\xrightarrow[]{P} 0$.

\vspace{1em}
\noindent\textbf{Proof for Terms (TR1)-(TR5)}
\vspace{1em}

The proofs of the convergence in probability to zero for the terms (TR1)-(TR5) follows a similar outline as Theorem 1, and we prove convergence on a case by case below. 

Terms (TR1) and (TR2) are similar; we only provide the proof of the convergence in probability to zero for the term (TR2).

To bound TR2, first set $\hat{\Delta}_{i\ell}=K_h(A_i - a) \left[\hat{\lambda}(a, X_i) - \lambda(a,X_i) \right]R(M_i, X_i)\left[Y_i - \gamma(X_i,M_i,a)\right]$. Bounding (TR2) amounts to showing $\mathbb{E}[\hat{\Delta}_{i\ell}|O^c_{I_\ell}] = o_p(1)$.

\begin{align*}
&\mathbb{E} \bigg[\hat{\Delta}_{i\ell} \bigg| O^c_{I_\ell}\bigg] \\
    =& \mathbb{E}\bigg\{K_h(A_i - a) 
    \left[\hat{\lambda}(a, X_i) - \lambda(a,X_i) \right]
    R(M_i, X_i)
    \left[Y_i - \gamma(X_i,M_i,a)\right]\bigg| O^c_{I_\ell}\bigg\}\\
    \intertext{Following identical steps as the proof for Theorem 1, gives}
    =&\int \left[ f(a \mid Y_i, M_i, X_i) +O(h^2) \right] \\
    &\quad \left[\hat{\lambda}(a, X_i) - \lambda(a,X_i) \right]R(M_i, X_i)\left[Y_i - \gamma(X_i,M_i,a)\right]f(Y_i,M_i, X_i) dY_idM_idX_i\\
    \overset{(a)}{=}&\int O(h^2) \left[\hat{\lambda}(a, X_i) - \lambda(a,X_i) \right]\\
    &\hspace{8em}R(M_i, X_i)\left[Y_i - \gamma(X_i,M_i,a)\right]f(Y_i,M_i, X_i) dY_idM_idX_i\\
     \overset{(b)}{=} &O(h^2)\int \Big|\left[\hat{\lambda}(a, X_i) - \lambda(a,X_i) \right]R(M_i, X_i)\Big|\\
    &\hspace{8em}\left[\int | Y_i-\gamma(X_i,M_i,a) |  f(Y_i\mid M_i, X_i) dY_i \right]  f(M_i,X_i)dM_idX_i\\
     \overset{(c)}{=}& o_p(1)
\end{align*}

where the equalities follow identically as in the proof of Theorem 1, and the final equality follows from $h \rightarrow 0$ along with the boundedness assumptions in Assumption 3.

For Term (TR3), we have
\begin{align*}
&\mathbb{E}\left[
\big(\hat{\eta}(a, a', X_i)-\eta(a, a', X_i)\big)\big(1-K_h(A_i - a')\lambda(a', X_i)\big)
\big| O^c_{I_\ell}\right] \\   
\intertext{Following a similar approach as used in the proof for Theorem 1, we have}
&{=}\int 
\big(\hat{\eta}(a, a', X_i)-\eta(a, a', X_i)\big)\big(1- f(a'\mid X_i)\lambda_{a'}(X_i)\big)f(X_i)dX_i \\
&\quad+\int 
\big(\hat{\eta}(a, a', X_i)-\eta(a, a', X_i)\big)O(h^2)\lambda_{a'}(X_i)f(X_i)dX_i\\
&\overset{(b)}{=} o_p(1).
\end{align*}
where (b) follows from the definition of $\lambda_{a^\prime}(X_i)$, $h\rightarrow 0$, Assumption 3 (boundedness of $\lambda$, $\hat{\eta}$ and $\eta$).

Demonstrating the bound for (TR4), we have
\begin{align*}
    &\mathbb{E}\big[
\big(\hat{\gamma}(X_i,M_i,a)-\gamma(X_i,M_i,a)\big)
\big\{
K_h(A_i - a')\lambda(a', X_i)
-K_h(A_i - a)\lambda(a,X_i)R(M_i,X_i)
\big\}
\big]\\
&= \mathbb{E}\big[
\big(\hat{\gamma}(X_i,M_i,a)-\gamma(X_i,M_i,a)\big)
\big\{
K_h(A_i - a')\lambda(a', X_i)
\big\}
\big]\tag{TR4-1}\\
 - &\mathbb{E}\big[
\big(\hat{\gamma}(X_i,M_i,a)-\gamma(X_i,M_i,a)\big)
\big\{
K_h(A_i - a)\lambda(a,X_i)R(M_i,X_i)
\big\}
\big]\tag{TR-4-2}
\end{align*}

TR-4-1 can be written as
\begin{align*}
    &\mathbb{E}\big[\big(\hat{\gamma}(X_i,M_i,a)-\gamma(X_i,M_i,a)\big)\big\{K_h(A_i - a')\lambda(a', X_i)\big\}\big]\\
    &= \int \big(\hat{\gamma}(X_i,M_i,a)-\gamma(X_i,M_i,a)\big)\lambda(a', X_i)\left\{ \int K_h(A_i - a') f(A_i \mid M_i, X_i) dA_i \right\} \\
    &\qquad \qquad \times f(M_i, X_i) dM_i dX_i
\end{align*}

An application of Lemma \ref{lem:kernelexpansion} to TR-4-1 gives
\begin{align*}
&\mathbb{E}\big[
\big(\hat{\gamma}(X_i,M_i,a)-\gamma(X_i,M_i,a)\big)
\big\{
K_h(A_i - a')\lambda(a', X_i)
\big\}
\big]\\
    &= \int \big(\hat{\gamma}(X_i,M_i,a)-\gamma(X_i,M_i,a)\big)\lambda(a', X_i)f(a' \mid M_i, X_i)f(M_i, X_i) dM_i dX_i\\
    & + \int \big(\hat{\gamma}(X_i,M_i,a)-\gamma(X_i,M_i,a)\big)\lambda(a', X_i)O(h^2)f(M_i, X_i) dM_i dX_i
\end{align*}

A similar approach applied to TR-4-2 gives
\begin{align*}
&\mathbb{E}\big[
\big(\hat{\gamma}(X_i,M_i,a)-\gamma(X_i,M_i,a)\big)
\big\{
K_h(A_i - a)\lambda(a,X_i)R(M_i,X_i)
\big\}
\big]\\
    =& \int \big(\hat{\gamma}(X_i,M_i,a)-\gamma(X_i,M_i,a)\big)\lambda(a,X_i)R(M_i,X_i)f(a \mid M, X) f(M_i, X_i) dM_i dX_i\\
    &+ \int \big(\hat{\gamma}(X_i,M_i,a)-\gamma(X_i,M_i,a)\big)\lambda(a,X_i)R(M_i,X_i)O(h^2)f(M_i, X_i) dM_i dX_i
\end{align*}

Now, the first terms of TR-4-1 and TR-4-2 cancel out with each other, with an identical proof to that used in the proof of Theorem 1.

Consequently only the remaining terms must be bounded. 
\[
\int \big(\hat{\gamma}(X_i,M_i,a)-\gamma(X_i,M_i,a)\big)\lambda(a', X_i) O(h^2)f(M, X) dM_i dX_i = o_p(1)
\]
The second term in TR-4-1 and TR-4-2 can be bounded from $h \rightarrow 0$, combined with the boundedness assumptions in Assumption 3.2.

Finally, for term (TR5), we note that
\begin{align*}
&\sqrt{nh^{d_A}}\mathbb{E}\left[
K_h(A_i - a')\big(\hat{\lambda}(a', X_i)-\lambda(a', X_i)\big)
\big\{
\gamma(X_i,M_i,a)-\eta(X_i)
\big\}\big| O^c_{I_\ell}
\right]\\
&=\sqrt{nh^{d_A}}\int
K_h(A_i - a')\big(\hat{\lambda}(a', X_i)-\lambda(a', X_i)\big)
\big\{
\gamma(X_i,M_i,a)-\eta(a, a', X_i)
\big\}\\
&\qquad \times f(A_i,M_i,X_i)dA_idM_idX_i\\
&=\sqrt{nh^{d_A}}\int
\Bigg\{\int K_h(A_i - a')f(A_i\mid M_i,X_i)dA_i\Bigg\}
\big(\hat{\lambda}(a', X_i)-\lambda(a', X_i)\big)\\
&\qquad \times \big\{
\gamma(X_i,M_i,a)-\eta(a, a', X_i)
\big\}f(M_i,X_i)dM_idX_i\\
&\overset{(a)}{=}\sqrt{nh^{d_A}}\int \big(\hat{\lambda}(a', X_i)-\lambda(a', X_i)\big)
\big\{
\gamma(X_i,M_i,a)-\eta(a, a', X_i)
\big\}f(a',M_i,X_i)dM_idX_i\\
&\quad+\sqrt{nh^{d_A}}\int \big(\hat{\lambda}(a', X_i)-\lambda(a', X_i)\big)
\big\{
\gamma(X_i,M_i,a)-\eta(a, a', X_i)
\big\}O(h^2)f(M_i,X_i)dM_idX_i\\
&\overset{(b)}{=} 0 +O(\sqrt{nh^{d_A + 4}})\int \Bigg|\big(\hat{\lambda}(a', X_i)-\lambda(a', X_i)\big)
\big\{
\gamma(X_i,M_i,a)-\eta(a, a', X_i)
\big\}\Bigg|\\
&\qquad \qquad \times f(M_i,X_i)dM_idX_i\\
&\overset{(c)}{=} o_p(1)
\end{align*}

Where $(a)$ follows from an application of Lemma \ref{lem:kernelexpansion}, (b) follows from the definition of $\eta$, and (c) follows from an application of Cauchy-Schwartz combined with the consistency of $\hat{\lambda}$.

\vspace{1em}
\subsection{Proof of Proposition 2}

Recall that
\begin{equation*}
\begin{aligned}
m(O;\alpha,\lambda,\gamma,\psi(a,a')) 
&= K_h(A - a)\lambda(a,X)R(M,X)\{Y - \gamma(X,M,a)\}\\
&\quad+ K_h(A - a')\lambda(a',X)\{\gamma(X,M,a) - \eta(a, a', X)\} \\
&\quad + \eta(a, a', X) - \psi(a,a'),
\end{aligned}
\end{equation*}

To prove consistency of $\widehat{V}(a, a')$, we first prove propositions \textbf{(I)}, \textbf{(II)} and \textbf{(III)}, which together prove the desired result. 

\noindent \textbf{(I)} $h^{d_A}n^{-1} \sum_{i\in I_\ell} m^2(O_i;\alpha,\lambda,\gamma,\psi(a,a')) - V(a,a') = o_p(1)$

To simplify notation, denote $m(O_i;\alpha,\lambda,\gamma,\psi(a,a'))$ as $m_i(a,a')$. From the proof of Theorem 1, we have $h^{d_A} \mathbb{E}[m^2_i(a,a')] = V(a,a') + o_p(1)$. 

We write
\begin{align*}
    &U_1(a,a') = K_h(A - a)\lambda(a,X)R(M,X)\{Y - \gamma(X,M,a)\},\\
    &U_2(a,a') = K_h(A - a')\lambda(a',X)\{\gamma(X,M,a) - \eta(a, a', X)\},\\
    &U_3(a,a') =\eta(a, a', X) - \psi(a,a').
\end{align*}
Then,
\begin{equation*}
    \begin{aligned}
        &\mathbb{E} (m^4_i) = \mathbb{E}[(U_1+U_2+U_3)^4]\\
        =&  \mathbb{E}(U_1^4) + 4\mathbb{E}(U_1^3U_2) + 4\mathbb{E}(U_1^3U_3) + 6\mathbb{E}(U_1^2U_2^2) + 12\mathbb{E}(U_1^2U_2U_3)+\\
        &6\mathbb{E}(U_1^2U_3^2) + 4\mathbb{E}(U_1U_2^3) + 12\mathbb{E}(U_1U_2^2U_3) + 12\mathbb{E}(U_1U_2U_3^2) + 4\mathbb{E}(U_1U_3^3) + \\
        &\mathbb{E}(U_2^4) + 4\mathbb{E}(U_2^3U_3) + 6\mathbb{E}(U_2^2U_3^2) + 4\mathbb{E}(U_2U_3^3) + \mathbb{E}(U_3^4) 
    \end{aligned}
\end{equation*}

We only need to investigate the terms $\mathbb{E}(U_1^{c_1} U_2^{c_2}U_3^{c_3})$ for any $c_1 \ge 0$, $c_2 \ge 0$ and $c_3 \ge 0$ with $c_1 + c_2 + c_3 = 4$. To be specific, dropping the terms with power index being zero, we will be studying $\mathbb{E}(U_1^{c_1})$, $\mathbb{E}(U_2^{c_2})$, $\mathbb{E}(U_1^{c_1} U_2^{c_2})$, $\mathbb{E}(U_1^{c_1} U_3^{c_3})$, $\mathbb{E}(U_2^{c_2} U_3^{c_3})$, and $\mathbb{E}(U_1^{c_1} U_2^{c_2}U_3^{c_3})$ for positive $c_1, c_2,$ and $c_3$.

\begin{enumerate}
    \item $\mathbb{E}(U_1^{c_1})$: By the assumed boundedness of $\lambda(a,X)$, $R(M,X)$, and $\mathbb{E}\{[Y-\gamma(X,M,a)]^4 | A = a', M = m, X = x\}$ over any $(a,a',m,x)\in \mathcal{A}\times\mathcal{A}\times \mathcal{M} \times \mathcal{X}$ from Assumption 7, 

    \begin{align*}
        &\mathbb{E}(U_1^{c_1}) = \int \bigg\{K_h(A - a)\lambda(a,X)R(M,X)[Y - \gamma(X,M,a)]\bigg\}^{c_1} f(Y,A, M,X) dO\\
        =& O(\frac{1}{h^{(c_1-1)d_A}})\int \tilde{k}(u)^{c_1}  \bigg\{\int |Y - \gamma(X,M,a)|^{c_1} f(Y|A = uh+a, M, X)dY \bigg\}\\
        & \qquad \times f(uh+a,M,X) dudMdX\\
        =& O(\frac{1}{h^{(c_1-1)d_A}})\int \tilde{k}(u)^{c_1} \mathbb{E}\{ |Y - \gamma(X,M,a)|^{c_1} |A = uh+a, M, X\}\\
        &\qquad\qquad f(uh+a,M,X) dudMdX\\
        =& O(\frac{1}{h^{(c_1-1)d_A}})\int \tilde{k}(u)^{c_1} f_{MX}(M,X)\\ &\bigg\{f(a|M,X)+ \sum^{d_A}_{j=1}u_jh\frac{\partial}{\partial a} f(a|M,X) + \sum^{d_A}_{j=1}\sum^{d_A}_{j'=1}u_ju_{j'}h^2\frac{\partial^2}{\partial a_j\partial a_{j'}} f(\bar{a}|M,X) \bigg\}dudMdX\\
        =& O(\frac{1}{h^{(c_1-1)d_A}})\int \tilde{k}(u)^{c_1} du + o(\frac{1}{h^{(c_1-1)d_A}}) =  O(\frac{1}{h^{(c_1-1)d_A}}).
    \end{align*}
    where $\bar{a}$ is between $a$ and $a+uh$.
    \item $\mathbb{E}(U_2^{c_2})$: \\
    From the boundedness of $\lambda(a',X)$, $\gamma(X,M,a)$ and $\eta(a, a', X)$ over any $(a,a',a'',m,x)\in \mathcal{A}^3\times \mathcal{M} \times \mathcal{X}$,
    \begin{align*}
    &\mathbb{E}(U_2^{c_2}) = \int \bigg\{K_h(A - a')\lambda(a',X)[\gamma(X,M,a) - \eta(a, a', X)]\bigg\}^{c_2} f(A, M,X) dO\\
    =& O(\frac{1}{h^{(c_2-1)d_A}})\int \tilde{k}(u)^{c_2}f_{MX}(M,X)   f(uh+a'|M,X) dudMdX\\
    =& O(\frac{1}{h^{(c_2-1)d_A}})\int \tilde{k}(u)^{c_2}f_{MX}(M,X) \\
    &\bigg\{f(a'|M,X)+ \sum^{d_A}_{j=1}u_jh\frac{\partial}{\partial a'} f(a'|M,X) \\
    &\qquad \qquad + \sum^{d_A}_{j=1}\sum^{d_A}_{j'=1}u_ju_{j'}h^2\frac{\partial^2}{\partial a'_j\partial a'_{j'}} f(\bar{a}|M,X) \bigg\}dudMdX\\
    =& O(\frac{1}{h^{(c_2-1)d_A}})\int \tilde{k}(u)^{c_2} du + o(\frac{1}{h^{(c_2-1)d_A}}) =  O(\frac{1}{h^{(c_2-1)d_A}}).
\end{align*}
where $\bar{a}$ is between $a'$ and $a'+uh$.
\item $\mathbb{E}(U_1^{c_1} U_2^{c_2})$
\begin{align*}
    &\mathbb{E}(U_1^{c_1} U_2^{c_2})\\
    =& \int \bigg\{K_h(A - a)\lambda(a,X)\frac{\alpha(a',M,X)}{\alpha(a,M,X)}[Y - \gamma(X,M,a)]\bigg\}^{c_1}\\
    &\qquad\times \bigg\{K_h(A - a')\lambda(a',X)[\gamma(X,M,a) - \eta(a, a', X)]\bigg\}^{c_2} f(Y,A, M,X) dO\\
    =& O(\frac{1}{h^{(c_1+c_2-1)d_A}})\int \tilde{k}(u)^{c_1}  \bigg\{\int |Y - \gamma(X,M,a)|^{c_1} f(Y|A = uh+a, M, X)dY \bigg\}\\
    &\qquad\times \tilde{k}(u+\frac{a-a'}{h})^{c_2} f_{MX}(M, X) f(uh+a|M,X) dudMdX \\
    =& O(\frac{1}{h^{(c_1+c_2-1)d_A}})\int \bigg[\prod^{d_A}_{j=1} k(u_j)^{c_1} k(u_j + \frac{a_j-a'_j}{h})^{c_2}\bigg] \\
    & \quad \mathbb{E}\{ |Y - \gamma(X,M,a)|^{c_1} |A = uh+a, M, X\}\times  f_{MX}(M, X) f(uh+a|M,X) dudMdX \\
    =& O(\frac{1}{h^{(c_1+c_2-1)d_A}})\int \bigg[\prod^{d_A}_{j=1} k(u_j)^{c_1} k(u_j + \frac{a_j-a'_j}{h})^{c_2}\bigg] f_{MX}(M, X) \\
&\bigg\{f(a|M,X)+ \sum^{d_A}_{j=1}u_jh\frac{\partial}{\partial a} f(a|M,X) + \sum^{d_A}_{j=1}\sum^{d_A}_{j'=1}u_ju_{j'}h^2\frac{\partial^2}{\partial a_j\partial a_{j'}} f(\bar{a}|M,X) \bigg\}dudMdX\\
    =& O(\frac{1}{h^{(c_1+c_2-1)d_A}})\int \bigg[\prod^{d_A}_{j=1} k(u_j)^{c_1} k(u_j + \frac{a_j-a'_j}{h})^{c_2}\bigg] du + o(\frac{1}{h^{(c_1+c_2-1)d_A}})\\
     =& O(\frac{1}{h^{(c_1+c_2-1)d_A}}).
\end{align*}
where $\bar{a}$ is between $a$ and $a+uh$.

\item $\mathbb{E}(U_1^{c_1} U_3^{c_3})$
\begin{align*}
    &\mathbb{E}(U_1^{c_1} U_3^{c_3})\\
    =& \int \bigg\{K_h(A - a)\lambda(a,X)\frac{\alpha(a',M,X)}{\alpha(a,M,X)}[Y - \gamma(X,M,a)]\bigg\}^{c_1}\bigg\{\eta(a, a', X) - \psi(a,a') \bigg\}^{c_2} \\
    &\qquad\times f(Y,A, M,X) dO\\
    =&O(1) \int K_h(A - a)^{c_1}|Y - \gamma(X,M,a)|^{c_1}f(Y,A, M,X) dO\\
    =&O(1) \int K_h(A - a)^{c_1}\mathbb{E}\bigg[|Y - \gamma(X,M,a)|^{c_1} \bigg| A, M, X \bigg]f(A, M,X) dO\\
    =&O(1) \int K_h(A - a)^{c_1}f(A \mid M,X) dA f_{MX}(M, X) dMdX\\
    =&O(\frac{1}{h^{(c_1 - 1)d_A}}) \int \tilde{k}(u)^{c_1}f_{MX}(M, X)\\
    &\bigg\{f(a|M,X)+ \sum^{d_A}_{j=1}u_jh\frac{\partial}{\partial a} f(a|M,X) + \sum^{d_A}_{j=1}\sum^{d_A}_{j'=1}u_ju_{j'}h^2\frac{\partial^2}{\partial a_j\partial a_{j'}} f(\bar{a}|M,X) \bigg\}dudMdX\\
    =& O(\frac{1}{h^{(c_1 - 1)d_A}})
\end{align*}
where $\bar{a}$ is between $a$ and $a+uh$, the second equality is from from the boundedness of $\lambda$, $\eta$, $\alpha$ and $\psi$, and the fourth equality comes from the assumed boundedness of $\mathbb{E}[|Y - \gamma|^4 \mid A, M, X]$.

\item $\mathbb{E}(U_2^{c_2} U_3^{c_3})$
\begin{align*}
    &\mathbb{E}(U_2^{c_2} U_3^{c_3})\\
    =& \int \bigg\{K_h(A - a')\lambda(a',X)[\gamma(X,M,a) - \eta(a, a', X)]\bigg\}^{c_2}\bigg\{\eta(a, a', X) - \psi(a,a')\bigg\}^{c_3}\\
    &\qquad\times f(Y,A, M,X) dO\\
    =& O(\frac{1}{h^{(c_2-1)d_A}})\int \tilde{k}(u+\frac{a-a'}{h})^{c_2}   f_{MX}(M, X) f_{A|X}(uh+a|M,X) dudMdX \\
    =& O(\frac{1}{h^{(c_2-1)d_A}})\int \tilde{k}(u+\frac{a-a'}{h})^{c_2}  f_{MX}(M, X)  \\
   &\bigg\{f(a'|M,X)+ \sum^{d_A}_{j=1}u_jh\frac{\partial}{\partial a'} f(a'|M,X) + \sum^{d_A}_{j=1}\sum^{d_A}_{j'=1}u_ju_{j'}h^2\frac{\partial^2}{\partial a'_j\partial a'_{j'}} f(\bar{a}|M,X) \bigg\}\\
   &\qquad\qquad dudMdX\\
    =& O(\frac{1}{h^{(c_2-1)d_A}})\int \tilde{k}(u+\frac{a-a'}{h})^{c_2}   du + o(\frac{1}{h^{(c_2-1)d_A}})\\
     =& O(\frac{1}{h^{(c_2-1)d_A}}).
\end{align*}
where $\bar{a}$ is between $a'$ and $a'+uh$.

\item $\mathbb{E}(U_1^{c_1} U_2^{c_2} U_3^{c_3} )$
\begin{align*}
    &\mathbb{E}(U_1^{c_1} U_2^{c_2}U_3^{c_3})\\
    =& \int \bigg\{K_h(A - a)\lambda(a,X)\frac{\alpha(a',M,X)}{\alpha(a,M,X)}[Y - \gamma(X,M,a)]\bigg\}^{c_1}\\
    &\bigg\{K_h(A - a')\lambda(a',X)[\gamma(X,M,a) - \eta(a, a', X)]\bigg\}^{c_2}\bigg\{\eta(a, a', X) - \psi(a,a')\bigg\}^{c_3} \\
&\qquad\times f(Y,A, M,X) dO\\
    =& O(\frac{1}{h^{(c_1+c_2-1)d_A}})\int \tilde{k}(u)^{c_1}  \bigg\{\int |Y - \gamma(X,M,a)|^{c_1} f(Y|A = uh+a, M, X)dY \bigg\}\\
    &\hspace{10em}\tilde{k}(u+\frac{a-a'}{h})^{c_2}   f_{MX}(M, X) f(uh+a|M,X) dudMdX\\
     =& O(\frac{1}{h^{(c_1+c_2-1)d_A}}).
\end{align*}
where the last equality is obtained as in the calculation for ${E}(U_1^{c_1} U_2^{c_2})$.
\end{enumerate}

Combining all the terms, we obtain $\mathbb{E} (m^4_i) = O(h^{-3d_A})$. Then by Markov inequality, for any $\epsilon > 0$,
\begin{align*}
   & P(|h^{d_A}n^{-1} \sum_{i\in I_\ell} m^2_i(a,a') - V(a,a') | > \epsilon) \le \frac{1}{\epsilon^2} \mathbb{E}\bigg\{\big[h^{d_A}n^{-1} \sum_{i\in I_\ell} m^2_i(a,a') - V(a,a') \big]^2\bigg\}\\
   =& \frac{1}{\epsilon^2} \mathbb{E}\bigg\{\big[h^{d_A}n^{-1} \sum_{i\in I_\ell} m^2_i(a,a') - h^{d_A}\mathbb{E}[m^2_i(a,a')] + o_p(1) \big]^2\bigg\}\\
    =& \frac{h^{2d_A}}{n^2\epsilon^2} \mathbb{E}\bigg\{\big[\sum_{i\in I_\ell} m^2_i - \mathbb{E}(\sum_{i\in I_\ell} m^2_i)\big]^2\bigg\} +  o_p(1)\\
    =& \frac{h^{2d_A}}{n^2\epsilon^2} var(\sum_{i\in I_\ell} m^2_i) +  o_p(1)\\
    =& \frac{h^{2d_A}}{n\epsilon^2} var( m^2_i) +  o_p(1)\\
    =& O(\frac{1}{nh^{d_A}} ) = o_p(1),
\end{align*}
where the equality in the last row comes from $var( m^2_i) = O(\mathbb{E} (m^4_i)) = O(h^{-3d_A})$.

\noindent \textbf{(II)}: $$h^{d_A}|I_\ell|^{-1}\sum_{i\in I_\ell}\mathbb{E}[ m^2(O_i;\hat{\alpha}_\ell,\hat{\lambda}_\ell,\hat{\gamma}_\ell,\hat{\psi}_\ell(a,a')) - m^2(O_i;\alpha,\lambda,\gamma,\psi(a,a')) \mid O^c_{I_\ell}] = o_p(1)$$

For simplicity in notation, we ignore the subscripts $\ell$ below for nuisance parameters estimated from $O^c_{I_\ell}$.
First, we analyze $h^{d_A}\mathbb{E}[m^2(O_i;\hat{\alpha},\hat{\lambda},\hat{\gamma},\hat{\psi}(a,a'))\mid O^c_{I_\ell}]$ as follows. 
We write
\begin{align*}
    &\hat{U}_1(a,a') = K_h(A - a)\hat{\lambda}(a,X)\hat{R}(M,X)\{Y - \hat{\gamma}(X,M,a)\},\\
    &\hat{U}_2(a,a') = K_h(A - a')\hat{\lambda}(a',X)\{\hat{\gamma}(X,M,a) - \hat{\eta}(a, a', X)\},\\
    &\hat{U}_3(a,a') =\hat{\eta}(a, a', X) - \hat{\psi}(a,a').
\end{align*}
Denote $m(O_i;\hat{\alpha},\hat{\lambda},\hat{\gamma},\hat{\psi}(a,a'))$ as $\hat{m}_i$. Then,
\begin{equation*}
    \begin{aligned}
        &\mathbb{E} (\hat{m}^2_i|O^c_{I_\ell}) = \mathbb{E}[(\hat{U}_1+\hat{U}_2+\hat{U}_3)^2|O^c_{I_\ell}]\\
        =&  \mathbb{E}(\hat{U}^2_1|O^c_{I_\ell}) + \mathbb{E}(\hat{U}^2_2|O^c_{I_\ell})+\mathbb{E}(\hat{U}^2_3|O^c_{I_\ell}) + 2\mathbb{E}(\hat{U}_1\hat{U}_2|O^c_{I_\ell})+2\mathbb{E}(\hat{U}_2\hat{U}_3|O^c_{I_\ell})+2\mathbb{E}(\hat{U}_1\hat{U}_3|O^c_{I_\ell})
    \end{aligned}
\end{equation*}
\begin{enumerate}
    \item $h^{d_A}\mathbb{E}(\hat{U}_1^{2} \mid  O^c_{I_\ell})$ \begin{align*}
    &h^{d_A}\mathbb{E}\Bigg(\bigg\{\frac{K_h(A - a)\hat{f}(M \mid A = a', X)}{\hat{f}(M \mid A = a, X)\hat{f}(a \mid X)}[Y - \hat{\gamma}(X,M,a)]  \bigg\}^2 \Bigg| O^c_{I_\ell}\Bigg)\\
    =&h^{d_A}\mathbb{E}\Bigg[\mathbb{E}\bigg(\bigg\{\frac{K_h(A - a)\hat{f}(M \mid A = a', X)}{\hat{f}(M \mid A = a, X)\hat{f}(a \mid X)}[Y - \hat{\gamma}(X,M,a)]  \bigg\}^2\bigg|X,M, O^c_{I_\ell}\bigg)\Bigg| O^c_{I_\ell}\Bigg]\\
    =&h^{d_A}\mathbb{E}\Bigg\{\frac{\hat{f}(M \mid A = a', X)^2}{\hat{f}(M \mid A = a, X)^2\hat{f}(a \mid X)^2}\times\\
    &\mathbb{E}\bigg[K_h(A - a)^2\mathbb{E}\{[Y - \hat{\gamma}(X,M,a)]^2|X,M,A, O^c_{I_\ell}\}\bigg|X,M, O^c_{I_\ell}\bigg] \Bigg| O^c_{I_\ell} \Bigg\}
\end{align*}
After adding and subtracting $\mathbb{E}[Y \mid X, M, A]$, the middle expectation can be written as 
\begin{align*}
    &h^{d_A}\mathbb{E}\bigg(K_h(A - a)^2 \bigg\{var(Y|X,M,A) + [\gamma(X,M,a)-\hat{\gamma}(X,M,a)]^2\bigg\}
     \bigg|X,M, O^c_{I_\ell} \bigg)\\
     = &\int \bigg[\prod^{d_A}_{j=1} k(u_j)^2 \bigg]\times \bigg\{var(Y|X,M,a+uh) + [\gamma(a+uh,M,X)-\hat{\gamma}(X,M,a)]^2\bigg\}
     \\
     &\times f(a+uh|X,M) du_1 \dots du_{d_A}\\
     = & \int k(u_1)^2\cdots k(u_{d_A})^2
     \times \bigg\{var(Y|X,M,a)+\sum^{d_A}_{j=1}u_jh\partial_{a_j} var(Y|X,M,a) +\\
    &\sum^{d_A}_{j=1}\sum^{d_A}_{j^\prime=1}u_ju_{j^\prime}h^2\partial_{a_j}\partial_{a_j^\prime} var(Y|X,M,\Bar{a}_v)  + \\
     &\bigg[ \gamma(X,M,a)-\hat{\gamma}(X,M,a)+\sum^{d_A}_{j=1} u_jh\partial_{a_j} \gamma(X,M, a)\\
     &\qquad\qquad + \sum^{d_A}_{j=1}\sum^{d_A}_{j'=1} u_ju_{j'}h^2\partial_{a_j}\partial_{a_{j'}} \gamma(X,M,\Bar{a}_{\gamma})\bigg]^2\bigg\}\\
     &\times \bigg[f(a|X,M)+\sum^{d_A}_{j=1}u_jh\partial_{a_j} f(a|X,M) \\
     &\qquad\qquad + \sum^{d_A}_{j=1}\sum^{d_A}_{j^\prime=1}u_ju_{j^\prime}h^2\partial_{a_j}\partial_{a_{j^\prime}} f(\Bar{a}_f|X,M) \bigg]du_1\cdots du_{d_A}\\
    \overset{(a)}{=}& \Big[\int \tilde{k}(u)^2 du\Big] \times  \Bigg\{ var(Y|X,M,a) + [\gamma(X,M,a) - \hat{\gamma}(X,M,a)]^2\Bigg\}f(a|X,M)\\
    &\qquad\qquad + O(h^2)\\
     =& \Big[\int \tilde{k}(u)^2 du\Big] \times  \mathbb{E}\Big\{[Y - \hat{\gamma}(X,M,a)]^2\Big|X,M, a, O^c_{I_\ell}\Big\}f(a|X,M) + O(h^2)
\end{align*}
where $\bar{a}_v, \bar{a}_{\gamma}$, and $\bar{a}_f$ are between $a$ and $a+h$. Equality (a) comes from the boundedness of $\int u^6k^2(u)du$, which is true because we assume $0 < \int u^6 k(u) du < \infty$. Plugging this back into the original expectation gives
\begin{align*}
    &\mathbb{E}\Bigg\{\frac{\hat{f}(M \mid A = a', X)^2}{\hat{f}(M \mid A = a, X)^2\hat{f}(a \mid X)^2}
    \Big[\int \tilde{k}(u)^2 du\Big] \\
    &\qquad\times  \mathbb{E}\{[Y - \hat{\gamma}(X,M,a)]^2|X,M, a, O^c_{I_\ell}\}\Bigg| O^c_{I_\ell} \Bigg\} + o_p(1)
\end{align*}

\item $h^{d_A}\mathbb{E}(\hat{U}_2^{2} \mid  O^c_{I_\ell})$
\begin{align*}
    &h^{d_A}\mathbb{E}\Bigg\{\bigg[\frac{K_h(A - a')}{\hat{f}(a' \mid X)}\Big(\hat{\gamma}(X,M,a) - \hat{\eta}(a, a', X)\Big) \bigg]^2 \Bigg| O^c_{I_{\ell}}\Bigg\}\\
    =&h^{d_A}\mathbb{E}\Bigg\{\frac{1}{\hat{f}(a' \mid X)^2}\mathbb{E}\bigg[K_h(A - a')^2\Big(\hat{\gamma}(X,M,a) - \hat{\eta}(a, a', X)\Big)^2\Big|X, O^c_{I_{\ell}}\bigg] \Bigg| O^c_{I_{\ell}}\Bigg\}
\end{align*}

The inner expectation can be written as
\begin{align*}
    &h^{d_A}\mathbb{E}\bigg[K_h(A - a')^2\Big(\hat{\gamma}(X,M,a) - \hat{\eta}(a, a', X)\Big)^2\Big|X,  O^c_{I_{\ell}} \bigg]\\
    \intertext{Following a similar kernel expansion as before}
    =& \int k^2(u_1)\cdots k^2(u_{d_A})\Big(\hat{\gamma}(X,M,a) - \hat{\eta}(a, a', X)\Big)^2f(a'|X, M)f(M\mid X) \\
    &\qquad du_1\cdots du_{d_A}dM+ O(h^2)
\end{align*}
Plugging this back into the original expectation leads to
\begin{align*}
& \int \tilde{k}(u)^2 du \times \mathbb{E}\bigg\{\frac{1}{\hat{f}(a'|X)^2}\Big(\hat{\gamma}(X,M,a) - \hat{\eta}(a, a', X)\Big)^2 \bigg| O^c_{I_\ell}\bigg\} + o_p(1)
\end{align*}

\item $h^{d_A}\mathbb{E}(\hat{U}_3^{2} \mid  O^c_{I_\ell})$ 
\begin{align*}
    h^{d_A}\mathbb{E}\Bigg\{[\hat{\eta}(a, a', X)-\hat{\psi}(a,a')]^2 \bigg| O^c_{I_\ell}\Bigg\} & = o_p(1)
\end{align*}
This holds because we assume the nuisance estimators are bounded, and following a similar calculation as the variance it can be seen that $h^{d_A}\mathbb{E}[\hat{\psi}^2(a, a') \mid  O^c_{I_\ell}] = o_p(1)$, which combined with Jensen's inequality can be used to obtain the desired result.

\item $h^{d_A}\mathbb{E}(\hat{U}_1\hat{U}_2|O^c_{I_\ell})$
\begin{align*}
    &h^{d_A}\mathbb{E}\Bigg\{\bigg[\frac{K_h(A - a)\hat{f}(M \mid A = a', X)}{\hat{f}(M \mid A = a, X)\hat{f}(a \mid X)}\{Y - \hat{\gamma}(X,M,a)\}  \bigg]\\
    &\qquad \times \bigg[\frac{K_h(A - a')}{\hat{f}(a' \mid X)}\{\hat{\gamma}(X,M,a) - \hat{\eta}(a, a', X)\} \bigg] \bigg|  O^c_{I_\ell} \Bigg\}\\
\end{align*}

\begin{align*}
    &h^{d_A}\mathbb{E}\Bigg\{\bigg[\frac{K_h(A - a)\hat{f}(M \mid A = a', X)}{\hat{f}(M \mid A = a, X)\hat{f}(a \mid X)}\{Y - \hat{\gamma}(X,M,a)\}  \bigg]\\
    &\qquad \times \bigg[\frac{K_h(A - a')}{\hat{f}(a' \mid X)}\{\hat{\gamma}(X,M,a) - \hat{\eta}(a, a', X)\} \bigg]\Bigg\}\\
    =&h^{d_A}\mathbb{E}\Bigg\{\frac{K_h(A-a)K_h(A-a')}{\hat{f}(a|X)\hat{f}(a'|X)}\frac{\hat{f}(M|a',X)}{\hat{f}(M|a,X)}\\
    &\qquad \qquad\times \Big[Y-\hat{\gamma}(X,M,a)\Big]\Big[\hat{\gamma}(X,M,a)-\hat{\eta}(a,a',X)\Big]\Bigg\}\\
    =&h^{d_A}\mathbb{E}\Bigg\{\frac{1}{\hat{f}(a|X)\hat{f}(a'|X)}\frac{\hat{f}(M|a',X)}{\hat{f}(M|a,X)} \Big[\hat{\gamma}(X,M,a)-\hat{\eta}(a,a',X)\Big]\\
    & \qquad\times\mathbb{E}\Big\{K_h(A-a)K_h(A-a')\Big[Y-\hat{\gamma}(X,M,a)\Big]\Big|X,M\Big\}\Bigg\}\\
    =&h^{d_A}\mathbb{E}\Bigg\{\frac{1}{\hat{f}(a|X)\hat{f}(a'|X)}\frac{\hat{f}(M|a',X)}{\hat{f}(M|a,X)}
    \Big[\hat{\gamma}(X,M,a) - \hat{\eta}(a,a',X)\Big]\\
    & \qquad\times\mathbb{E}\Big\{K_h(A-a)K_h(A-a')\Big[\gamma(X,M,A)-\hat{\gamma}(X,M,a)\Big]\Big|X,M\Big\}\Bigg\}
\end{align*}

The inner expectation
\begin{align*}
    &h^{d_A}\mathbb{E}\Big\{K_h(A-a)K_h(A-a')\Big[\gamma(X,M,A)-\hat{\gamma}(X,M,a)\Big]\Big|X,M\Big\}\\
    =& h^{d_A}\int \bigg[\prod^{d_A}_{j=1}\frac{1}{h^2}k\Big(\frac{A_j-a}{h}\Big)k\Big(\frac{A_j-a'}{h}\Big)\bigg]
  \Big[\gamma(X,M,A)-\hat{\gamma}(X,M,a)\Big] f(A|X,M)dA\\
    =& \int \tilde{k}(u) \tilde{k}(u+\frac{a-a'}{h})\Big[\gamma(X,M, uh + a)-\hat{\gamma}(X,M,a)\Big] f(uh+ a|X,M)du\\
    =& \int k(u_1)\cdots k(u_{d_A}) k(u_1+\frac{a-a'}{h})\cdots k(u_{d_A}+\frac{a-a'}{h})\\  
    &\times \Big[(\gamma(X,M,a) - \hat{\gamma}(X,M,a)) + \sum^{d_A}_{j=1} u_j h\partial_{a_j} \gamma(X,M,a)+\\
    &\qquad \frac{u^2_j h^2}{2}\partial^2_{a_j}\gamma(X,M,a)+\frac{u^3_jh^3}{6}\partial^3_{a_j}\gamma(X,M,\Bar{a}_\gamma)\Big]\\
    &\times \Big[f(a|X,M) +\sum^{d_A}_{j=1} u_jh\partial_{a_j} f(a|X,M)+\frac{u^2_jh^2}{2}\partial^2_{a_j} f(\Bar{a}_f|X,M)\Big]du_1\cdots du_{d_A}
\end{align*}

where $\bar{a}_{\gamma}$ and $\bar{a}_f$ are between $a$ and $a+h$. Inserting this back into the full expectation combined with Assumption 4 bounds this term as $o_p(1)$.

\item $h^{d_A}\mathbb{E}(\hat{U}_1\hat{U}_3|O^c_{I_\ell})$
\begin{align*}
        &2h^{d_A}\mathbb{E}\Bigg\{\bigg[\hat{\eta}(a, a', X) - \hat{\psi}(a,a')\bigg]\bigg[\frac{K_h(A - a)\hat{f}(M \mid A = a', X)}{\hat{f}(M \mid A = a, X)\hat{f}(a \mid X)}\{Y - \hat{\gamma}(X,M,a)\}  \bigg] \bigg|  O^c_{I_\ell} \Bigg\} \\
        &= o_p(1)
\end{align*}

Expanding this into two terms,
\begin{align*}
    &2h^{d_A}\mathbb{E}\Bigg\{\hat{\eta}(a, a', X) \bigg[\frac{K_h(A - a)\hat{f}(M \mid A = a', X)}{\hat{f}(M \mid A = a, X)\hat{f}(a \mid X)}\{Y - \hat{\gamma}(X,M,a)\}  \bigg] \bigg|  O^c_{I_\ell} \Bigg\} \\
    &- 2h^{d_A}\mathbb{E}\Bigg\{\hat{\psi}(a,a') \bigg[\frac{K_h(A - a)\hat{f}(M \mid A = a', X)}{\hat{f}(M \mid A = a, X)\hat{f}(a \mid X)}\{Y - \hat{\gamma}(X,M,a)\}  \bigg] \bigg|  O^c_{I_\ell} \Bigg\} 
\end{align*}

The first term can be bounded as $o_p(1)$ using a similar approach used above, and for the second term, from the i.i.d assumption on the data we can re-write it as
 \begin{align*}
    &2h^{d_A}\mathbb{E}\Bigg\{\hat{\psi}(a,a') \bigg[\frac{K_h(A - a)\hat{f}(M \mid A = a', X)}{\hat{f}(M \mid A = a, X)\hat{f}(a \mid X)}\{Y - \hat{\gamma}(X,M,a)\}  \bigg] \bigg|  O^c_{I_\ell} \Bigg\} \\
        &= 2h^{d_A} |I_\ell|^{-1}\Bigg(\sum_{i \in I_\ell}\mathbb{E}\Bigg\{\bigg[\frac{K_h(A_i - a)\hat{f}(M_i \mid A = a', X_i)}{\hat{f}(M_i \mid A = a, X_i)\hat{f}(a \mid X_i)}\{Y_i - \hat{\gamma}(X_i,M_i,a)\}  \bigg]^2  \bigg|  O^c_{I_\ell} \Bigg\}\\
        &+ \sum_{i \in I_\ell}\mathbb{E}\Bigg\{\bigg[\frac{K_h(A_i - a)\hat{f}(M_i \mid A = a', X_i)}{\hat{f}(M_i \mid A = a, X_i)\hat{f}(a \mid X_i)}\{Y_i - \hat{\gamma}(X_i,M_i,a)\}  \bigg]\times \\
        &\qquad \bigg[\frac{K_h(A - a')}{\hat{f}(a' \mid X)}\Big(\hat{\gamma}(X,M,a) - \hat{\eta}(a, a', X)\Big) \bigg] \bigg|  O^c_{I_\ell} \Bigg\}\\
        &+ \sum_{i \in I_\ell}\mathbb{E}\Bigg\{\bigg[\frac{K_h(A_i - a)\hat{f}(M_i \mid A = a', X_i)}{\hat{f}(M_i \mid A = a, X_i)\hat{f}(a \mid X_i)}\{Y_i - \hat{\gamma}(X_i,M_i,a)\}  \bigg]\hat{\eta}(a, a', X) \bigg|  O^c_{I_\ell} \Bigg\}
\end{align*}

By the boundedness of$$\mathbb{E}\{[Y - \hat{\gamma}(X,M,a)]^2|X,M,A, O^c_{I_\ell}\} = var(Y|X,M,A) + [\gamma(X,M,a)-\hat{\gamma}(X,M,a)]^2$$ from Assumption 3, and following the results in the first part of (II), we know that $h^{d_A}\mathbb{E}\bigg[K_h(A - a)^2[Y - \hat{\gamma}(X,M,a)]^2\bigg|X,M, O^c_{I_\ell}\bigg]$ is bounded.
Thus, the first term is $O(|I_\ell|^{-1}) = o_p(1)$ from the law of total expectation. Because $\hat{f}(a'|X)$, $\hat{\gamma}$, and $\hat{\eta}$ are bounded by assumptions, the boundedness of $h^{d_A}\mathbb{E}\bigg[K_h(A_i - a)K_h(A_j-a')[Y - \hat{\gamma}(X,M,a)]\bigg|X,M, O^c_{I_\ell}\bigg]$ can be obtained similar to the third part of (I). Hence, the second term also has  $O(|I_\ell|^{-1}) = o_p(1)$.  From the boundedness of $h^{d_A/2}\mathbb{E}\bigg[K_h(A - a)[Y - \hat{\gamma}(X,M,a)]\bigg|X,M, O^c_{I_\ell}\bigg]$ based on Jensen's inequality and the boundedness of $\hat{\eta}$, the third term satisfies $O(h^{d_A/2}|I_\ell|^{-1}) = o_p(1)$. As a result, $h^{d_A}\mathbb{E}(\hat{U}_1\hat{U}_3|O^c_{I_\ell}) = o_p(1)$.

\item $h^{d_A}\mathbb{E}(\hat{U}_2\hat{U}_3|O^c_{I_\ell})$
\begin{align*}
    &h^{d_A}\mathbb{E}\Bigg\{\bigg[\hat{\eta}(a, a', X) - \hat{\psi}(a,a')\bigg]
    \bigg[\frac{K_h(A - a')}{\hat{f}(a' \mid X)}\{\hat{\gamma}(X,M,a) - \hat{\eta}(a, a', X)\} \bigg] \bigg|  O^c_{I_\ell} \Bigg\}\\
\end{align*}
From the boundedness of $\hat{\gamma}$, $\hat{\eta}$, and $\hat{f}(a'|X)$, there is $$h^{d_A}\mathbb{E}\Bigg\{\hat{\eta}(a, a', X)
    \bigg[\frac{K_h(A - a')}{\hat{f}(a' \mid X)}\{\hat{\gamma}(X,M,a) - \hat{\eta}(a, a', X)\} \bigg] \bigg|  O^c_{I_\ell} \Bigg\} = O(h^{d_A}) = o_p(1).$$ 
A similar proof as the fifth part of (II) above can show that the second term is also $o_p(1)$.
\end{enumerate}

Combining all the six parts, we have $h^{d_A}\mathbb{E}[m^2(O_i;\hat{\alpha},\hat{\lambda},\hat{\gamma},\hat{\psi}(a,a'))]$ equal to
\begin{align*}
    &\Big[\int \tilde{k}(u)^2 du\Big]\mathbb{E}\Bigg\{\frac{\hat{f}(M \mid A = a', X)^2}{\hat{f}(M \mid A = a, X)^2\hat{f}(a \mid X)^2} \times  \mathbb{E}\{(Y - \hat{\gamma}(X,M,a))^2|X,M,a, O^c_{I_\ell}\}\Bigg| O^c_{I_\ell} \Bigg\} \\
+ &\Big[\int \tilde{k}(u)^2 du\Big]  \times \mathbb{E}\bigg\{\frac{1}{\hat{f}(a'|X)^2}\Big(\hat{\gamma}(X,M,a) - \hat{\eta}(a, a', X)\Big)^2 \Bigg| O^c_{I_\ell}\bigg\} + o_p(1).
\end{align*}

Next, $h^{d_A}\mathbb{E}[m^2(O_i;\alpha,\lambda,\gamma,\psi(a,a'))]$ can be written as
\begin{align*}
    &\Big[\int \tilde{k}(u)^2 du\Big]\mathbb{E}\Bigg\{\frac{f(M \mid A = a', X)^2}{f(M \mid A = a, X)^2f(a \mid X)^2} \times  \mathbb{E}\{(Y - \gamma(X,M,a))^2|X,M,a, O^c_{I_\ell}\}\Bigg| O^c_{I_\ell} \Bigg\} \\
+ &\Big[\int \tilde{k}(u)^2 du\Big] \times \mathbb{E}\bigg\{\frac{1}{f(a'|X)^2}\Big(\gamma(X,M,a) - \eta(a, a', X)\Big)^2 \Bigg| O^c_{I_\ell}\bigg\} + o_p(1)
\end{align*}

Define $\int \tilde{k}(u)^2 du = R^2_{d_A}$, 
\begin{align*}
\omega_1 = &R^2_{d_A}\mathbb{E}\Bigg\{\frac{\hat{f}(M \mid A = a', X)^2}{\hat{f}(M \mid A = a, X)^2\hat{f}(a \mid X)^2} \times  \mathbb{E}\{(Y - \hat{\gamma}(X,M,a))^2|X,M, a, O^c_{I_\ell}\}\Bigg| O^c_{I_\ell} \Bigg\} \\
 & \quad - R^2_{d_A}\mathbb{E}\Bigg\{\frac{f(M \mid A = a', X)^2}{f(M \mid A = a, X)^2f(a \mid X)^2} \times  \mathbb{E}\{(Y - \gamma(X,M,a))^2|X,M, a, O^c_{I_\ell}\}\Bigg| O^c_{I_\ell} \Bigg\}, \end{align*}
 and
 \begin{align*}
\omega_2= &R^2_{d_A}\mathbb{E}\bigg\{\frac{1}{\hat{f}(a'|X)^2}\Big(\hat{\gamma}(X,M,a) - \hat{\eta}(a, a', X)\Big)^2 \Bigg| O^c_{I_\ell}\bigg\}\\
 &\quad -R^2_{d_A}\mathbb{E}\bigg\{\frac{1}{f(a'|X)^2}\Big(\gamma(X,M,a) - \eta(a, a', X)\Big)^2 \Big| O^c_{I_\ell}\Big\}.
\end{align*}
Then $h^{d_A}\mathbb{E}[m^2(O_i;\hat{\alpha},\hat{\lambda},\hat{\gamma},\hat{\psi}(a,a'))]- h^{d_A}\mathbb{E}[m^2(O_i;\alpha,\lambda,\gamma,\psi(a,a'))] = \omega_1 + \omega_2 + o_p(1).$ First, we focus on simplifying $\omega_2$, which equals
\begin{align*}
    R^2_{d_A}\mathbb{E}\bigg\{\frac{1}{\hat{f}(a'|X)^2}\Big(\hat{\gamma}(X,M,a) - \hat{\eta}(a, a', X)\Big)^2 - \frac{1}{f(a'|X)^2}\Big(\gamma(X,M,a) - \eta(a, a', X)\Big)^2 \Big| O^c_{I_\ell} \Big\}.
\end{align*}
From expressing $\frac{1}{\hat{f}(a'|X)}\Big(\hat{\gamma}(X,M,a) - \hat{\eta}(a, a', X)\Big)$ as 
\begin{align*}
& \frac{1}{f(a'|X)}\Big(\gamma(X,M,a) - \eta(a, a', X)\Big) 
    + \frac{1}{f(a'|X)}\Big(\hat{\gamma}(X,M,a) - \gamma(X,M,a)\Big)\\
    &+ \frac{1}{f(a'|X)}\Big(\eta(a, a', X) - \hat{\eta}(a, a', X)\Big) + \Big(\hat{\gamma}(X,M,a) - \hat{\eta}(a, a', X)\Big)\Big(\frac{1}{\hat{f}(a'|X)} - \frac{1}{f(a'|X)}\Big),
\end{align*}
there is
\begin{align*}
    &\omega_2\\
    =&R^2_{d_A}\mathbb{E}\bigg\{\frac{1}{\hat{f}(a'|X)^2}\Big(\hat{\gamma}(X,M,a) - \hat{\eta}(a, a', X)\Big)^2 - \frac{1}{f(a'|X)^2}\Big(\gamma(X,M,a) - \eta(a, a', X)\Big)^2 \Big| O^c_{I_\ell} \Big\}\\
    =& 
  R^2_{d_A}\mathbb{E}\bigg\{  \frac{1}{f^2(a'|X)}\Big(\hat{\gamma}(X,M,a) - \gamma(X,M,a)\Big)^2
    + \frac{1}{f^2(a'|X)}\Big(\eta(a, a', X)^2 - \hat{\eta}(a, a', X)\Big)^2\\
    &+ \Big(\hat{\gamma}(X,M,a) - \hat{\eta}(a, a', X)\Big)^2\Big(\frac{1}{\hat{f}(a'|X)} - \frac{1}{f(a'|X)}\Big)^2 \\
    &+ 2\frac{1}{f(a'|X)}\Big(\gamma(X,M,a) - \eta(a, a', X)\Big)\frac{1}{f(a'|X)}\Big(\hat{\gamma}(X,M,a) - \gamma(X,M,a)\Big)\\
    &+ 2\frac{1}{f(a'|X)}\Big(\gamma(X,M,a) - \eta(a, a', X)\Big)\frac{1}{f(a'|X)}\Big(\eta(a, a', X) - \hat{\eta}(a, a', X)\Big)\\
    &+ 2\frac{1}{f(a'|X)}\Big(\gamma(X,M,a) - \eta(a, a', X)\Big)\Big(\hat{\gamma}(X,M,a) - \hat{\eta}(a, a', X)\Big)\Big(\frac{1}{\hat{f}(a'|X)} - \frac{1}{f(a'|X)}\Big)\\
    &+ 2\frac{1}{f(a'|X)}\Big(\hat{\gamma}(X,M,a) - \gamma(X,M,a)\Big)\frac{1}{f(a'|X)}\Big(\eta(a, a', X) - \hat{\eta}(a, a', X)\Big)\\
    &+ 2\frac{1}{f(a'|X)}\Big(\hat{\gamma}(X,M,a) - \gamma(X,M,a)\Big)\Big(\hat{\gamma}(X,M,a) - \hat{\eta}(a, a', X)\Big)\Big(\frac{1}{\hat{f}(a'|X)} - \frac{1}{f(a'|X)}\Big)\\
    &+2\frac{1}{f(a'|X)}\Big(\eta(a, a', X) - \hat{\eta}(a, a', X)\Big)\Big(\hat{\gamma}(X,M,a) - \hat{\eta}(a, a', X)\Big)\Big(\frac{1}{\hat{f}(a'|X)} - \frac{1}{f(a'|X)}\Big)\\
    &\qquad \Big| O^c_{I_\ell} \Big\} 
\end{align*}

We show each of these terms are $o_p(1)$ as follows. Because $f^2(a'| X)$ is bounded away from $0$ based on Assumption 3 (ii) and the consistency of $\hat{\gamma}$ from Assumption 4(iii), there is $$\mathbb{E}\bigg\{ \frac{1}{f^2(a'|X)}\Big[\hat{\gamma}(X,M,a) - \gamma(X,M,a)\Big]^2\Big| O^c_{I_\ell} \Big\} = o_p(1).$$
Under a similar argument and from Assumption 4(iv), $$\mathbb{E}\bigg\{ \frac{1}{f^2(a'|X)}\Big[\eta(a, a', X) - \hat{\eta}(a, a', X)\Big]^2\Big| O^c_{I_\ell} \Big\} = o_p(1).$$
Based on the boundedness of nuisance estimators from Assumption 3(ii) and the consistency of $\hat{f}(a'|X)$ from Assumption 4(i), there is
$$\mathbb{E}\bigg\{ \Big[\hat{\gamma}(X,M,a) - \hat{\eta}(a, a', X)\Big]^2\Big[\frac{1}{\hat{f}(a'|X)} - \frac{1}{f(a'|X)}\Big]^2\Big| O^c_{I_\ell} \Big\} = o_p(1).$$
Each of the remaining cross terms is a product of a term that is $o_p(1)$ from the estimator's consistency and a term that is bounded. Hence, we have $\omega_2 = o_p(1)$. 

Next, we employ a similar derivation to simplify $\omega_1$.
Note that
\begin{align*}
   & \frac{\hat{f}(M \mid A = a', X)}{\hat{f}(M \mid A = a, X)\hat{f}(a \mid X)}(Y - \hat{\gamma}(X,M,a)) \\
   =&
    \frac{f(M \mid A = a', X)}{f(M \mid A = a, X)f(a \mid X)}(Y - \gamma(X,M,a))\\
     &+\Big(\frac{\hat{f}(M \mid A = a', X)}{\hat{f}(M \mid A = a, X)} - \frac{f(M \mid A = a', X)}{f(M \mid A = a, X)}\Big)\frac{1}{f(a \mid X)}(Y - \gamma(X,M,a))\\
    &+ \frac{\hat{f}(M \mid A = a', X)}{\hat{f}(M \mid A = a, X)}\Big(\frac{1}{\hat{f}(a \mid X)} - \frac{1}{f(a \mid X)}\Big)(Y - \gamma(X,M,a))\\
    &+  \frac{\hat{f}(M \mid A = a', X)}{\hat{f}(M \mid A = a, X)\hat{f}(a \mid X)}(\gamma(X,M,a) - \hat{\gamma}(X,M,a)).
\end{align*}
Hence, 
\begin{align*}
\omega_1=&R^2_{d_A}\mathbb{E}\Bigg\{\mathbb{E}\bigg\{\frac{\hat{f}(M \mid A = a', X)^2}{\hat{f}(M \mid A = a, X)^2\hat{f}(a \mid X)^2}(Y - \hat{\gamma}(X,M,a))^2\bigg|X,M, a, O^c_{I_\ell}\bigg\} \\
& \qquad -  \mathbb{E}\bigg\{\frac{f(M \mid A = a', X)^2}{f(M \mid A = a, X)^2f(a \mid X)^2}(Y - \gamma(X,M,a))^2 \bigg|X,M, a, O^c_{I_\ell}\bigg\}\Bigg| O^c_{I_\ell} \Bigg\}\\
=& R^2_{d_A} \mathbb{E}\Bigg\{\mathbb{E}\bigg\{\bigg[\Big(\frac{\hat{f}(M \mid A = a', X)}{\hat{f}(M \mid A = a, X)} - \frac{f(M \mid A = a', X)}{f(M \mid A = a, X)}\Big)\frac{1}{f(a \mid X)}(Y - \gamma(X,M,a))\\
    &\quad+ \frac{\hat{f}(M \mid A = a', X)}{\hat{f}(M \mid A = a, X)}\Big(\frac{1}{\hat{f}(a \mid X)} - \frac{1}{f(a \mid X)}\Big)(Y - \gamma(X,M,a))\\
    &\quad+  \frac{\hat{f}(M \mid A = a', X)}{\hat{f}(M \mid A = a, X)\hat{f}(a \mid X)}(\gamma(X,M,a) - \hat{\gamma}(X,M,a))\bigg]^2\bigg|X,M, a, O^c_{I_\ell}\bigg\}\Bigg| O^c_{I_\ell} \Bigg\}\\
    +&2 R^2_{d_A} \mathbb{E}\Bigg\{\mathbb{E}\bigg\{\bigg[\Big(\frac{\hat{f}(M \mid A = a', X)}{\hat{f}(M \mid A = a, X)} - \frac{f(M \mid A = a', X)}{f(M \mid A = a, X)}\Big)\frac{1}{f(a \mid X)}(Y - \gamma(X,M,a))\\
    &\quad+ \frac{\hat{f}(M \mid A = a', X)}{\hat{f}(M \mid A = a, X)}\Big(\frac{1}{\hat{f}(a \mid X)} - \frac{1}{f(a \mid X)}\Big)(Y - \gamma(X,M,a))\\
    &\quad+  \frac{\hat{f}(M \mid A = a', X)}{\hat{f}(M \mid A = a, X)\hat{f}(a \mid X)}(\gamma(X,M,a) - \hat{\gamma}(X,M,a))\bigg]\times\\
    &\quad  \bigg[\frac{f(M \mid A = a', X)}{f(M \mid A = a, X)f(a \mid X)}(Y - \gamma(X,M,a))\bigg]\bigg|X,M, a, O^c_{I_\ell}\bigg\}\Bigg| O^c_{I_\ell} \Bigg\}
\end{align*}

After further expansions, we can show that the squared terms contain a component that is bounded based on Assumption 3 and another component that is $o_p(1)$ from Assumption 4. The $(Y - \gamma(X,M,a))^2$ in some squared terms is integrated out as a bounded component due to $var(Y|X,M,a)$ being bounded as assumed in Assumption 3(3). For interaction terms, those containing $(Y - \gamma(X,M,a))$ equals zero because $\int (Y - \gamma(X,M,a))f(Y|X,M,a)dY = 0$. All of the interaction terms contain a bounded component and a $o_p(1)$ component. Consequently, $\omega_1 = o_p(1)$, leading to $h^{d_A}\mathbb{E}[m^2(O_i;\hat{\alpha},\hat{\lambda},\hat{\gamma},\hat{\psi}(a,a'))]- h^{d_A}\mathbb{E}[m^2(O_i;\alpha,\lambda,\gamma,\psi(a,a'))] = o_p(1).$

\noindent \textbf{(III)} $h^{d_A}|I_\ell|^{-1} \sum_{i\in I_\ell} \Delta_{i\ell} = o_p(1)$, where 

\begin{align*}
&\Delta_{i\ell} = m^2(O_i;\hat{\alpha},\hat{\lambda},\hat{\gamma},\hat{\psi}(a,a')) - m^2(O_i;\alpha,\lambda,\gamma,\psi(a,a')) \\
&\qquad - \mathbb{E}\Big\{m^2(O_i;\hat{\alpha},\hat{\lambda},\hat{\gamma},\hat{\psi}(a,a')) - m^2(O_i;\alpha,\lambda,\gamma,\psi(a,a'))\Big| O^c_{I_\ell}\Big\}.
\end{align*}

By Lemma \ref{lem:bounding-term}, it suffices to bound $\mathbb{E}\left[\Big(h^{d_A}|I_\ell|^{-1} \sum_{i\in I_\ell} \Delta_{i\ell} \Big)^2 \Big| O^c_{I_\ell}\right] =h^{2d_A}|I_\ell|^{-1} \mathbb{E}\left[\Delta^2_{i\ell}  \Big| O^c_{I_\ell}\right]  $ as $o_p(1)$. Note that $\mathbb{E}[\Delta_{i\ell}] = 0$ and interaction terms are zero due to conditional independence. We start with analyzing $\mathbb{E}[\Delta^2_{i\ell} \mid O^c_{I_\ell}]$ as follows. For simplicity of notation, we adopt the notation definitions in parts (I) and (II), ignoring the subscripts $\ell$ for nuisance estimators. We have 
\begin{align*}
    \mathbb{E}[\Delta^2_{i\ell} \mid O^c_{I_\ell}] =&  \mathbb{E}\bigg\{ \big[\hat{m}^2_{i} - m^2_i - \mathbb{E}(\hat{m}^2_{i} - m^2_i\mid O^c_{I_\ell}) \big]^2\mid O^c_{I_\ell}\bigg\}\\
    =& \mathbb{E}\bigg[\big(\hat{m}^2_{i} - m^2_i \big)^2\Big| O^c_{I_\ell}\bigg] - \mathbb{E} \left(\hat{m}^2_{i} - m^2_i \mid O^c_{I_\ell}\right)^2
\end{align*}
From (II), we know that $\mathbb{E} \left(\hat{m}^2_{i} - m^2_i \mid O^c_{I_\ell}\right)^2 = o_p( h^{-2d_A})$. To bound $\mathbb{E}\bigg[\big(\hat{m}^2_{i} - m^2_i \big)^2\Big| O^c_{I_\ell}\bigg]$, by $\hat{m}_{i} = \hat{U}_1 + \hat{U}_2 + \hat{U}_3$ and $m_i = U_1 + U_2 + U_3$, we can rewrite the term as
\begin{align*}
    &\mathbb{E}\bigg[\big(\hat{m}^2_{i} - m^2_i \big)^2\Big| O^c_{I_\ell}\bigg]=\mathbb{E}\bigg\{ \big[(\hat{U}_1 + \hat{U}_2 + \hat{U}_3)^2 - (U_1+U_2+U_3)^2 \big]^2\Big| O^c_{I_\ell}\bigg\}\\
    =& \mathbb{E}\bigg\{ \big[(\hat{U}^2_1 + \hat{U}^2_2 + \hat{U}^2_3 +2\hat{U}_1\hat{U}_2 +2\hat{U}_2\hat{U}_3+2\hat{U}_1\hat{U}_3) \\
    &\qquad - (U^2_1+U^2_2+U^2_3+2U_1U_2 +2U_2U_3 + 2U_1U_3) \big]^2\Big| O^c_{I_\ell}\bigg\}\\
     =& \mathbb{E}\bigg\{ \big[(\hat{U}^2_1-U^2_1) + (\hat{U}^2_2-U^2_2) + (\hat{U}^2_3-U^2_3) \\
     &\qquad+2(\hat{U}_1\hat{U}_2 -U_1U_2+\hat{U}_2\hat{U}_3-U_2U_3+\hat{U}_1\hat{U}_3 - U_1U_3) \big]^2\Big|O^c_{I_\ell}\bigg\}\\
      =& \mathbb{E}\bigg\{ (\hat{U}^2_1-U^2_1)^2 + (\hat{U}^2_2-U^2_2)^2 + (\hat{U}^2_3-U^2_3)^2 +\sum_{\bar{c} \in \mathcal{W}} c_0\hat{U}^{c_1}_1\hat{U}^{c_2}_2\hat{U}^{c_3}_3 U^{c_4}_1U^{c_5}_2U^{c_6}_3 \Big| O^c_{I_\ell}\bigg\}
\end{align*}
where $\bar{c} = (c_1, \ldots, c_6)$ and $\mathcal{W}$ represents the possible combinations of $\bar{c}$ from the decomposition. We will prove that $\mathbb{E}\bigg\{\hat{U}^{c_1}_1\hat{U}^{c_2}_2\hat{U}^{c_3}_3 U^{c_4}_1U^{c_5}_2U^{c_6}_3 \Big|O^c_{I_\ell}\bigg\} = O(h^{-(c_1 + c_2 + c_4 + c_5 - 1)d_A})$. Note that
$$\mathbb{E}\bigg\{\hat{U}^{c_1}_1\hat{U}^{c_2}_2\hat{U}^{c_3}_3 U^{c_4}_1U^{c_5}_2U^{c_6}_3 \Big| O^c_{I_\ell}\bigg\} = \iint \hat{U}^{c_1}_1\hat{U}^{c_2}_2\hat{U}^{c_3}_3 U^{c_4}_1U^{c_5}_2U^{c_6}_3 f(Y, A, M, X \mid O^c_{I_\ell})dYdMdAdX.$$ By the boundedness of nuisance parameters and their estimates (Assumption 3(ii)), the above term equals
{\small $$ O\Big(\iint K_h(A - a)^{c_1 + c_4}K_h(A - a')^{c_2 + c_5}\Big|[Y - \hat{\gamma}(X,M,a)]^{c_{1}}[Y - \gamma(X,M,a)]^{c_{4}}\Big| f(Y, A, M, X \mid O^c_{I_\ell})dYdMdAdX \Big).$$}
The possible combinations of $c_1, c_4$ in $\bar{c}$ are $$\{(c_1, c_4): (1, 1), (2, 0), (0, 2), (2, 1), (1, 2), (3, 0), (0, 3)\}.$$ Similar to the derivation in part (I), we will prove that the rate is $O(h^{-(c_1 + c_2 + c_4 + c_5 - 1)d_A})$ case-by-case. For the terms with $c_1 = 0$, the boundedness of $\mathbb{E}[|Y - \gamma|^4 \mid X, M, A]$ from Assumption 7 provides the boundedness of lower moments by separately considering the regions on which $|Y-\gamma|^{c_4}$ is $\ge$ or $< 1$. Next, we prove for the remaining terms.
\begin{enumerate}
    \item $c_1 > 0$ and $c_4 = 0$. The integral can be written as 
    \begin{align*}
        &\iint K_h^{c_1 + c_4}(A - a)K_h^{c_2 + c_5}(A - a')\big|Y - \hat{\gamma}(X,M,a)\big|^{c_1} f(Y, A, M, X \mid O^c_{I_\ell})dYdMdAdX \\
        =& \iint K_h^{c_1 + c_4}(A - a)K_h^{c_2 + c_5}(A - a')\mathbb{E}\left[\big|Y - \hat{\gamma}(X,M,a)\big|^{c_1} \Big| A, M, X, O^c_{I_\ell}\right]\\
        &\qquad \times f(A, M, X \mid O^c_{I_\ell})dMdAdX \end{align*}
    The inner expectation $\mathbb{E}\left[\big|Y - \hat{\gamma}(X,M,a)\big|^{c_1} \big| A, M, X, O^c_{I_\ell}\right]$ can be bounded as follows,
    \begin{align*}
        &\mathbb{E}\left[\big|Y - \hat{\gamma}(X,M,a)\big|^{c_1} \Big| A, M, X, O^c_{I_\ell}\right]\\
        =& \mathbb{E}\left[\big|Y - \gamma(X,M,a) + \gamma(X,M,a) - \hat{\gamma}(X,M,a)\big|^{c_1} \Big| A, M, X, O^c_{I_\ell}\right]\\
        \leq& \mathbb{E}[|Y - \gamma(X,M,a)|^{c_1} \mid A, M, X] + \mathbb{E}[|\gamma(X,M,a) - \hat{\gamma}(X,M,a)|^{c_1} \mid A, M, X, O^c_{I_\ell}]\\
        &\quad +  \sum^{c_1-1}_{k=1} \binom{c_1}{k} \mathbb{E}[|Y - \gamma(X,M,a)|^{k}|\gamma(X,M,a) - \hat{\gamma}(X,M,a)|^{c_1-k} \big| A, M, X, O^c_{I_\ell}].
    \end{align*}
Each of the terms in the expansion can be bounded from Assumption 3(ii) combined with the boundedness of $\mathbb{E}\left[(Y - \gamma(X, M , a))^4 \big| X, M, A\right]$ from Assumption 7. Hence, the original integral equals
\begin{align*}
    &O\left(\iint K_h(A - a)^{c_1 + c_4}K_h(A - a')^{c_2 + c_5} f(A, M, X \mid O^c_{I_\ell}) dMdAdX\right)\\
    &= O(h^{-(c_1 + c_2 + c_4 + c_5 - 1)d_A})
    \end{align*}
    where the last equality holds from the boundedness of the integrals of the kernels. 

    \item $c_1 > 0$ and $c_4 > 0$. The integral is 
    \begin{align*}
        &\iint K_h(A - a)^{c_1 + c_4}K_h(A - a')^{c_2 + c_5}\Big|[Y - \hat{\gamma}(X,M,a)]^{c_1}[Y - \gamma(X,M,a)]^{c_4}\Big| \\
        &\qquad \times f(Y, A, M, X \mid O^c_{I_\ell})dYdAdMdX\\
        =& \iint K_h(A - a)^{c_1 + c_4}K_h(A - a')^{c_2 + c_5}\\
        &\mathbb{E}\left[\big|[Y - \hat{\gamma}(X,M,a)]^{c_1}[Y - \gamma(X,M,a)]^{c_4}\big| \Big| A, M, X, O^c_{I_\ell}\right] f(A, M, X \mid O^c_{I_\ell})dAdMdX
        \end{align*}
        The inner expectation $\mathbb{E}\left[\big|[Y - \hat{\gamma}(X,M,a)]^{c_1}[Y - \gamma(X,M,a)]^{c_4}\big| \Big| A, M, X, O^c_{I_\ell}\right]$ can be bounded with
        \begin{align*}
        &\mathbb{E}\left[\big|[Y - \hat{\gamma}(X,M,a)]^{c_1}[Y - \gamma(X,M,a)]^{c_4}\big| \Big| A, M, X, O^c_{I_\ell}\right]\\
        =& \mathbb{E}\left\{\big|[Y - \gamma(X,M,a) + \gamma(X,M,a) - \hat{\gamma}(X,M,a)]^{c_1}[Y - \gamma(X,M,a)]^{c_4}\big| \Big| A, M, X, O^c_{I_\ell}\right\}\\
        \leq& \mathbb{E}\{|Y - \gamma(X,M,a)|^{c_1+c_4} \mid A, M, X\} \\
        & \quad + \mathbb{E}[|\gamma(X,M,a) - \hat{\gamma}(X,M,a)|^{c_1} |Y - \gamma(X,M,a)|^{c_4} \mid A, M, X, O^c_{I_\ell}]\\
        &\quad +\sum^{c_1-1}_{k=1} \binom{c_1}{k} \mathbb{E}[|\gamma(X,M,a) - \hat{\gamma}(X,M,a)|^{c_1-k}|Y - \gamma(X,M,a)|^{k+c_4} \mid A, M, X, O^c_{I_\ell}]
        \end{align*}
        The first term is the conditional variance, which is bounded by Assumption 3(3). The second and third terms can be bounded from Assumption 3(ii) and Assumption 7. 
   The bound of the integral follows similarly as before. 
\end{enumerate}

The remaining terms to bound are $\mathbb{E}[\hat{U}^2_1-U^2_1)^2\mid O^c_{I_\ell}]$, $ \mathbb{E}[(\hat{U}^2_2-U^2_2)^2 \mid O^c_{I_\ell}]$, and $ \mathbb{E}[(\hat{U}^2_3-U^2_3)^2 \mid O^c_{I_\ell}]$. First, $\mathbb{E}[(\hat{U}^2_3-U^2_3)^2 \mid O^c_{I_\ell}]$ can be bounded from Assumption 3(ii). Next, we demonstrate the boundedness of $\mathbb{E}[(\hat{U}^2_2-U^2_2)^2 \mid O^c_{I_\ell}]$; a similar derivation applies to $\mathbb{E}[\hat{U}^2_1-U^2_1)^2\mid O^c_{I_\ell}]$. To start with, we re-express the term
    \begin{align*}
    &\mathbb{E}[(\hat{U}^2_2-U^2_2)^2 \mid O^c_{I_\ell}]\\
        =&\mathbb{E}\Big(K^4_h(A - a')\Big\{\frac{1}{\hat{f}(a'|X)^2}\big[\hat{\gamma}(X,M,a) - \hat{\eta}(a, a', X)\big]^2 \\
        &\qquad - \frac{1}{f(a'|X)^2}\big[\gamma(X,M,a) - \eta(a, a', X)\big]^2\Big\}^2 \Big| O^c_{I_\ell}\Big)\\
        = &\mathbb{E}\Big(\Big\{\frac{1}{\hat{f}(a'|X)^2}\big[\hat{\gamma}(X,M,a) - \hat{\eta}(a, a', X)\big]^2 - \frac{1}{f(a'|X)^2}\big[\gamma(X,M,a) - \eta(a, a', X)\big]^2\Big\}^2\\
        &\quad \times\mathbb{E}[K^4_h(A - a') \mid X, M]\Big| O^c_{I_\ell}\Big).
        \end{align*}
     From Assumption 7, 
     \begin{align*}
         \mathbb{E}[K^4_h(A - a') \mid X, M] = & \int \prod^{d_A}_{j=1}\Big[ \frac{1}{h^4}k(\frac{A_j-a'_j}{h})^4\Big]f(A|X,M)dA\\
         =&h^{-3d_A} \int \tilde{k}(u)^4 f(uh+a'|X,M)du
         = O(h^{-3d_A}).
     \end{align*} 
     Hence,
     \begin{align*}
       &\mathbb{E}[(\hat{U}^2_2-U^2_2)^2 \mid O^c_{I_\ell}]\\ =&O(h^{-3d_A})\mathbb{E}\Big(\Big\{\frac{1}{\hat{f}(a'|X)^2}\big[\hat{\gamma}(X,M,a) - \hat{\eta}(a, a', X)\big]^2 \\
       &\qquad- \frac{1}{f(a'|X)^2}\big[\gamma(X,M,a) - \eta(a, a', X)\big]^2\Big\}^2\Big| O^c_{I_\ell}\Big).
\end{align*}
From the expansion of $\omega_2$ in proving term 6 of the part (II), we can express $\frac{1}{\hat{f}(a'|X)^2}\big[\hat{\gamma}(X,M,a) - \hat{\eta}(a, a', X)\big]^2 - \frac{1}{f(a'|X)^2}\big[\gamma(X,M,a) - \eta(a, a', X)\big]^2$ as a summation of 9 components, i.e.
    \begin{align*}
&\frac{1}{\hat{f}(a'|X)^2}\Big(\hat{\gamma}(X,M,a) - \hat{\eta}(a, a', X)\Big)^2 - \frac{1}{f(a'|X)^2}\Big(\gamma(X,M,a) - \eta(a, a', X)\Big)^2 \\
    =& 
   \frac{1}{f^2(a'|X)}\Big(\hat{\gamma}(X,M,a) - \gamma(X,M,a)\Big)^2
    + \frac{1}{f^2(a'|X)}\Big(\eta(a, a', X)^2 - \hat{\eta}(a, a', X)\Big)^2\\
    &+ \Big(\hat{\gamma}(X,M,a) - \hat{\eta}(a, a', X)\Big)^2\Big(\frac{1}{\hat{f}(a'|X)} - \frac{1}{f(a'|X)}\Big)^2 \\
    &+ 2\frac{1}{f(a'|X)}\Big(\gamma(X,M,a) - \eta(a, a', X)\Big)\frac{1}{f(a'|X)}\Big(\hat{\gamma}(X,M,a) - \gamma(X,M,a)\Big)\\
    &+ 2\frac{1}{f(a'|X)}\Big(\gamma(X,M,a) - \eta(a, a', X)\Big)\frac{1}{f(a'|X)}\Big(\eta(a, a', X) - \hat{\eta}(a, a', X)\Big)\\
    &+ 2\frac{1}{f(a'|X)}\Big(\gamma(X,M,a) - \eta(a, a', X)\Big)\Big(\hat{\gamma}(X,M,a) - \hat{\eta}(a, a', X)\Big)\Big(\frac{1}{\hat{f}(a'|X)} - \frac{1}{f(a'|X)}\Big)\\
    &+ 2\frac{1}{f(a'|X)}\Big(\hat{\gamma}(X,M,a) - \gamma(X,M,a)\Big)\frac{1}{f(a'|X)}\Big(\eta(a, a', X) - \hat{\eta}(a, a', X)\Big)\\
    &+ 2\frac{1}{f(a'|X)}\Big(\hat{\gamma}(X,M,a) - \gamma(X,M,a)\Big)\Big(\hat{\gamma}(X,M,a) - \hat{\eta}(a, a', X)\Big)\Big(\frac{1}{\hat{f}(a'|X)} - \frac{1}{f(a'|X)}\Big)\\
    &+2\frac{1}{f(a'|X)}\Big(\eta(a, a', X) - \hat{\eta}(a, a', X)\Big)\Big(\hat{\gamma}(X,M,a) - \hat{\eta}(a, a', X)\Big)\Big(\frac{1}{\hat{f}(a'|X)} - \frac{1}{f(a'|X)}\Big)
\end{align*}
For the multiplication of any two of the nine components chosen with replacement, the corresponding conditional expectation $\mathbb{E}(\cdot|O^c_{I_\ell})$ is a construct of a subcomponent that is $o_p(1)$ from the consistency of nuisance parameters multiplied by other subcomponents that are bounded from Assumption 3. As a consequence,  we obtained that $\mathbb{E}[(\hat{U}^2_2-U^2_2)^2 \mid O^c_{I_\ell}] = o_p(h^{-3d_A})$. A similar argument can be used to prove $\mathbb{E}[\hat{U}^2_1-U^2_1)^2\mid O^c_{I_\ell}] = o_p(h^{-3d_A})$ by utilizing the boundedness of $\mathbb{E}[(Y - \gamma)^4 \mid X, M, A]$ from Assumption 7 (i).

Because $ c_1+c_2+c_4+c_5 \le 4$, $O(1) \le O(h^{-(c_1 + c_2 + c_4 + c_5 - 1)d_A}) \le O(h^{-3d_A})$. 
We conclude that $\mathbb{E}\left[\Delta^2_{i\ell}  \Big| O^c_{I_\ell}\right]  = O(h^{-3d_A})$ and
$$\mathbb{E}\left[\Big(h^{d_A}|I_\ell|^{-1} \sum_{i\in I_\ell} \Delta_{i\ell} \Big)^2 \Big| O^c_{I_\ell}\right] =h^{2d_A}|I_\ell|^{-1} \mathbb{E}\left[\Delta^2_{i\ell}  \Big| O^c_{I_\ell}\right]  =O\left([nh^{d_A}]^{-1}\right) = o_p(1).$$ 

\section{Practical Implications of Irregularity}
\label{sec:irregular}

The pointwise continuous-treatment mediation functional considered here is irregular because it involves point evaluation at the treatment values $a$ and $a'$.
In a nonparametric model, a parameter is regular, or pathwise differentiable, if its pathwise derivative can be represented as a continuous linear functional of the score, equivalently by a square-integrable influence function \citep{bickel1993efficient,van2000asymptotic}. This property is what underlies classical root-$n$ efficient influence function-based estimation.

The exact pointwise mediation functional
\(
\psi_0(a,a')
=
\int \mu(a,m,x) f_{M\mid A,X}(m\mid a',x) f_X(x)\,d\nu(m)\,d\nu(x)
\)
does not have this property in the nonparametric continuous-treatment model. To see why, consider the part of the formal influence function corresponding to the outcome regression. If $A$ were discrete, the corresponding residual term would contain the ordinary inverse-probability factor
\[
\frac{\mathbbm{1}(A=a)}{\Pr(A=a\mid X)}
\frac{f_{M\mid A,X}(M\mid a',X)}
     {f_{M\mid A,X}(M\mid a,X)}
\{Y-\mu(a,M,X)\},
\]
which is square-integrable under positivity and boundedness conditions. When $A$ is continuous, however, $\Pr(A=a\mid X)=0$. The analogous expression for the exact pointwise intervention would require the generalized weight
\[
\frac{\delta_a(a)}{f_{A\mid X}(a\mid X)}
\frac{f_{M\mid A,X}(M\mid a',X)}
     {f_{M\mid A,X}(M\mid a,X)}
\{Y-\mu(a,M,X)\},
\]
where $\delta_a(\cdot)$ denotes a Dirac mass at $a$. This object is not an ordinary square-integrable random variable. Analogous Dirac-mass terms arise for the mediator-density component evaluated at $a'$. Hence, the exact pointwise functional $\psi_0(a,a')$ is not pathwise differentiable in the nonparametric continuous-treatment model. Consequently, unlike in the binary-treatment case, there is no finite-variance influence function for the exact pointwise target and no model-free root-$n$ regular estimator based on a classical efficient influence function.

The proposed estimator addresses this nonregularity by replacing the Dirac mass with a kernel approximation. That is, terms involving point evaluation at $A=a$ are regularized by weights of the form
\(
K_h(A-a)=h^{-d_A}K\left(\frac{A-a}{h}\right).
\)
For fixed $h>0$, this corresponds to a regularized, smoothed functional. However, as $h\to 0$, the variance of the kernel-weighted term increases. In particular, under standard regularity conditions,
\[
E\{K_h(A-a)^2\mid X\}
\approx
h^{-d_A} f_{A\mid X}(a\mid X)\int K(u)^2\,du,
\]
so the squared $L_2$ norm of the localized influence-function component is of order $h^{-d_A}$. Therefore the stochastic error of the estimator is of order $(nh^{d_A})^{-1/2}$ rather than $n^{-1/2}$. This is the sense in which the effective sample size is local and of order $nh^{d_A}$.

This irregularity has direct practical implications for the proposed estimators. First, standard errors and confidence intervals must be based on the local rate $(nh^{d_A})^{-1/2}$ rather than the usual root-$n$ rate. Second, bandwidth selection is an inferentially important bias--variance trade-off: smaller bandwidths better approximate the exact pointwise target but increase variance, whereas larger bandwidths improve stability but correspond to a more heavily smoothed target. For inference centered at the exact pointwise parameter $\psi_0(a,a')$, the smoothing bias, of order $h^2$ under the smoothness conditions used here, must be negligible relative to $(nh^{d_A})^{-1/2}$, for example through undersmoothing or explicit bias correction. If this condition is not imposed, the estimator is more naturally interpreted as targeting a smoothed version of the mediation functional, with possible residual smoothing bias.
Third, the irregularity makes overlap near the target treatment values especially important. When few observations fall in the local neighborhoods of $a$ or $a'$, the effective sample size can be small and inverse-density weights can be unstable. This motivates the practical diagnostics and stabilizations used in the application, including inspection of localized weights, effective sample size calculations, H\'ajek-type stabilization, clipping sensitivity analyses, and bandwidth sensitivity analyses. Multiple robustness and cross-fitting reduce sensitivity to nuisance-model misspecification and nuisance-estimation error, but they do not remove the fundamental nonregularity induced by point evaluation with a continuous treatment.

\section{Stabilization of the inverse-propensity weights in the application}
\label{sec:hajek}

In the application in Section 6 of the main manuscript, two modifications of the inverse-propensity component of the estimator were used to improve finite-sample stability: a H\'ajek-type stabilized propensity density estimator, and clipping of that estimator from below. This section gives the consistency of the former and discusses the properties of the latter. Throughout, $\hat f(a|X)$ denotes the model-based conditional treatment density estimator, written $\hat f_0(a|X)$ in the main manuscript.

\subsection{Consistency of the H\'ajek-type propensity estimator in cross-fitting}

Given a consistent estimator of propensity score at treatment value $a$ for person $i$ in cross validation fold $I_\ell$, $\hat{f}(a|X_i)$, we define the corresponding Hajek-type stabilized weighted propensity score as follows, $$\hat{f}(a|X_i) \times \frac{1}{|I_{-\ell}|}\sum^{|I_{-\ell}|}_{j\in I_{-\ell}} \frac{K_h (A_j-a)}{\hat{f}(a|X_j)}. $$

The goal here is to prove that
$$\lim_{|I_{-\ell}| \rightarrow \infty}\hat{f}(a|X_i) \times \frac{1}{|I_{-\ell}|}\sum^{|I_{-\ell}|}_{j\in I_{-\ell}} \frac{K_h (A_j-a)}{\hat{f}(a|X_j)}\stackrel{p}{=} f(a|X_i). $$

\begin{enumerate}
    \item 
    \begin{align*}
        &\lim_{|I_{-\ell}| \rightarrow \infty}\frac{1}{|I_{-\ell}|}\sum^{|I_{-\ell}|}_{j\in I_{-\ell}} \frac{K_h (A_j-a)}{\hat{f}(a|X_j)} \stackrel{p}{=} \iint \frac{K_h(A-a)}{f(a|X)} f(A,X) dAdX \\
        =&\int \bigg\{ \int K_h(A-a) \frac{f(A,X)}{f(a|X)}dA\bigg\}dX\\
        & \text{by Lemma 2}\\
        =& \iint \prod^{d_A}_{k=1} k(u_k)\bigg\{ \frac{f(a,X)}{f(a|X)} + \sum^{d_A}_{k=1} u_k h \frac{\partial_{a_k}f(a,X)}{f(a|X)} + \frac{1}{2} \sum^{d_A}_{k=1} \sum^{d_A}_{k'=1} u_k u_{k'}h^2\frac{\partial_{a_k}\partial_{a_{k'}}f(a,X)}{f(a|X)}\bigg|_{\bar{a}}\bigg\} \\
        &\qquad du_1\ldots du_{d_A} dX\\
        & \text{assume } \int\partial_{a_k}\partial_{a_{k'}}f(a,X)|_{\bar{a}}dX < \infty \text{ for } \bar{a} \text{ between } a \text{ and } a+uh, \text{ then}\\
        =& \int f(X) dX + O(h^2) = 1 + O(h^2)
    \end{align*} 
    \item  $\lim_{|I_{-\ell}| \rightarrow \infty}\hat{f}(a|X_i) = f(a|X_i)$ from the consistency of the propensity estimator $\hat{f}$
    \item Combining the first two bullets, we get
    \begin{align*}
        &\lim_{|I_{-\ell}| \rightarrow \infty}\hat{f}(a|X_i) \times  \frac{1}{|I_{-\ell}|}\sum^{|I_{-\ell}|}_{j\in I_{-\ell}} \frac{K_h (A_j-a)}{\hat{f}(a|X_j)}\\
        =& \lim_{|I_{-\ell}| \rightarrow \infty}\hat{f}(a|X_i) \times
        \lim_{|I_{-\ell}| \rightarrow \infty}\frac{1}{|I_{-\ell}|}\sum^{|I_{-\ell}|}_{j\in I_{-\ell}} \frac{K_h (A_j-a)}{\hat{f}(a|X_j)}\\
        \stackrel{p}{=}& f(a|X_i)
    \end{align*}
\end{enumerate}

The stabilization rescales the localized inverse-propensity weights and reduces sensitivity to random fluctuations in the denominator of the weighted estimator.

\subsection{Clipping of the stabilized propensity density}

Because the localized continuous-treatment weights can be unstable when $\hat f_H(a|X_i)$ is close to zero, especially near treatment values with limited empirical overlap, in the application we lower bounded the stabilized propensity density by a small constant $c$,
\[
\hat f_{H,c}(a|X_i)=\max\{\hat f_H(a|X_i),c\},
\qquad c=0.01,
\]
and used $1/\hat f_{H,c}(a|X_i)$ in place of $1/\hat f_H(a|X_i)$, where $\hat f_H(a|X_i)$ denotes the H\'ajek-type stabilized estimator above. This clipping step is a finite-sample bias--variance trade-off: it can reduce the influence of observations with extremely large inverse-propensity weights, but it is not a substitute for the positivity assumption. If the clipping threshold is fixed, the resulting estimator should be interpreted as a stabilized finite-sample implementation; asymptotically, clipping is innocuous when it is inactive with high probability, or when a threshold sequence tending to zero is used and the true density is bounded away from zero as in Assumption 3 of the main manuscript.

\section{Application Data Summary}
The table below presents summary statistics of the following variables: outcome $Y$, mediator $M$, treatment $A$, and confounders $X$. Missing values in confounders are addressed by including the indicators of missingness as covariates. 

\begin{table}
\centering
\setlength{\tabcolsep}{1.5pt} 
\renewcommand{\arraystretch}{0.9}
\resizebox{\linewidth}{!}{%
\begin{tabular}{@{}lrrrrrrrrr@{}} 
\hline
        \textbf{} & \textbf{Missing(\%)} & \textbf{Median (IQR)} & \textbf{0} & \textbf{1} & \textbf{2} & \textbf{3} & \textbf{4} & \textbf{5} & \textbf{6} \\ \hline
        female & ~ & ~ & ~ & 43.55 & ~ & ~ & ~ & ~ & ~ \\ 
        age & ~ & 18 (17-20) & ~ & ~ & ~ & ~ & ~ & ~ & ~ \\  
        white & ~ & ~ & ~ & 24.88 & ~ & ~ & ~ & ~ & ~ \\  
        black & ~ & ~ & ~ & 50.35 & ~ & ~ & ~ & ~ & ~ \\  
        Hispanic & ~ & ~ & ~ & 17.2 & ~ & ~ & ~ & ~ & ~ \\  
        years of education & 1.38 & 10 (9-11) & ~ & ~ & ~ & ~ & ~ & ~ & ~ \\  
        GED diploma & ~ & ~ & ~ & 4.12 & ~ & ~ & ~ & ~ & ~ \\  
        high school diploma & ~ & ~ & ~ & 18.22 & ~ & ~ & ~ & ~ & ~ \\  
        native English & ~ & ~ & ~ & 84.6 & ~ & ~ & ~ & ~ & ~ \\  
        divorced & ~ & ~ & ~ & 0.75 & ~ & ~ & ~ & ~ & ~ \\  
        separated & ~ & ~ & ~ & 1.27 & ~ & ~ & ~ & ~ & ~ \\  
        cohabiting & ~ & ~ & ~ & 3.15 & ~ & ~ & ~ & ~ & ~ \\  
        married & ~ & ~ & ~ & 1.62 & ~ & ~ & ~ & ~ & ~ \\  
        has children & ~ & ~ & ~ & 17.77 & ~ & ~ & ~ & ~ & ~ \\  
        ever worked & ~ & ~ & ~ & 14.4 & ~ & ~ & ~ & ~ & ~ \\  
        average weekly earnings in USD & ~ & 0 (0-0) & ~ & ~ & ~ & ~ & ~ & ~ & ~ \\  
        is household head & ~ & ~ & ~ & 10.38 & ~ & ~ & ~ & ~ & ~ \\  
        household size & 1.4 & 3 (2-5) & ~ & ~ & ~ & ~ & ~ & ~ & ~ \\  
        designated for nonresidential slot & ~ & ~ & ~ & 17.22 & ~ & ~ & ~ & ~ & ~ \\  
        total household gross income & 37.3 & ~ & ~ & 24.52 & 20.14 & 12.04 & 10.01 & 7.89 & 7.22 \\  
        total personal gross income & 55.65 & ~ & ~ & 92.78 & 5.02 & 1.24 & 0.51 & 0.28 & 0.06 \\  
        mum's years of education & 18.42 & 12 (11-12) & ~ & ~ & ~ & ~ & ~ & ~ & ~ \\  
        dad's years of education & 37.35 & 12 (11-12) & ~ & ~ & ~ & ~ & ~ & ~ & ~ \\  
        dad did not work at 14 & ~ & ~ & ~ & 5.2 & ~ & ~ & ~ & ~ & ~ \\  
        received AFDC per month & ~ & ~ & ~ & 23.08 & ~ & ~ & ~ & ~ & ~ \\  
        received public assistance per month & ~ & ~ & ~ & 20.03 & ~ & ~ & ~ & ~ & ~ \\  
        received food stamps & ~ & ~ & ~ & 42.95 & ~ & ~ & ~ & ~ & ~ \\  
        welfare receipt during childhood & 6.85 & ~ & ~ & 46.78 & 20.91 & 11.33 & 20.99 & ~ & ~ \\  
        poor/fair health & ~ & ~ & ~ & 12.4 & ~ & ~ & ~ & ~ & ~ \\  
        physical/emotional problems & ~ & ~ & ~ & 4.28 & ~ & ~ & ~ & ~ & ~ \\  
        extent of marijuana use & 63.28 & ~ & 18.58 & 9.73 & 12.87 & 16.88 & 41.93 & ~ & ~ \\  
        extent of hallucinogen use & 94.9 & ~ & 25.49 & 2.45 & 2.94 & 8.33 & 60.78 & ~ & ~ \\  
        ever used other illegal drugs & ~ & ~ & ~ & 0.45 & ~ & ~ & ~ & ~ & ~ \\  
        extent of smoking & 47.9 & ~ & 3.5 & 63.92 & 17.03 & 7.44 & 8.11 & ~ & ~ \\  
        extent of alcohol consumption & 42.35 & ~ & 8.72 & 1.21 & 10.93 & 25.85 & 53.3 & ~ & ~ \\  
        ever arrested & ~ & ~ & ~ & 23.75 & ~ & ~ & ~ & ~ & ~ \\  
        times in prison & ~ & ~ & 94.9 & 3.92 & 0.78 & 0.2 & 0.17 & 0.03 & ~ \\  
        time spent by Job Corps recruiter & 1.95 & ~ & ~ & 32 & 40.54 & 17.62 & 9.84 & ~ & ~ \\  
        extent of recruiter support & 2.22 & ~ & ~ & 64.13 & 26.75 & 1.48 & 0.84 & 6.8 & ~ \\  
        idea about wished training & ~ & ~ & ~ & 84.05 & ~ & ~ & ~ & ~ & ~ \\  
        expected hourly wage after training & 55.02 & 8 (7-10) & ~ & ~ & ~ & ~ & ~ & ~ & ~ \\  
        expected improvement in maths & 2.1 & ~ & ~ & 70.76 & 26.28 & 2.96 & ~ & ~ & ~ \\  
        expected improvement in reading skills & ~ & ~ & 1.7 & 54.77 & 34.88 & 8.65 & ~ & ~ & ~ \\  
        expected improvement in reading skills & ~ & ~ & 1.7 & 61.52 & 26.52 & 10.25 & ~ & ~ & ~ \\  
        expected to be training for a job & 1.95 & ~ & ~ & 96.07 & 3.39 & 0.54 & ~ & ~ & ~ \\  
        worried about training & ~ & ~ & ~ & 36.27 & ~ & ~ & ~ & ~ & ~ \\  
        1st contact with recruiter by phone & ~ & ~ & ~ & 40.75 & ~ & ~ & ~ & ~ & ~ \\  
        1st contact with recruiter in office & ~ & ~ & ~ & 22.8 & ~ & ~ & ~ & ~ & ~ \\  
        expected stay in training & ~ & 0 (0-12) & ~ & ~ & ~ & ~ & ~ & ~ & ~ \\  
        total training hours in yr 1 ($A$) & ~ & 965.71 (404.79-1767.21) & ~ & ~ & ~ & ~ & ~ & ~ & ~ \\  
        proportion of weeks employed in yr 2 ($M$) & ~ & 40.38 (0-80.77) & ~ & ~ & ~ & ~ & ~ & ~ & ~ \\  
        any arrests in yr 4 ($Y$) & ~ & ~ & ~ & 8.7 & ~ & ~ & ~ & ~ & ~ \\ \hline
    \end{tabular}}
    \caption{Descriptive statistics. The table provides [median (interquartile range: 25th - 75th percentiles)] for numeric variables, proportion (in percentage) for each level of categorical variables, and the proportion of 1's (in percentage) for dummy variables. For the total household/personal gross income, the proportion percentages do not add up to 100$\%$ because the last level was eliminated for simplicity of display. The data has a sample size of 4,000.}
     \label{Xtab}
\end{table}

\end{document}